\documentclass[10pt,a4paper]{article}
\usepackage[utf8]{inputenc}
\usepackage[english]{babel}
\usepackage{amsmath}
\usepackage{amsfonts}
\usepackage{amssymb}
\usepackage{graphicx}
\usepackage{tikz}
\usepackage{wrapfig}
\usepackage{caption}
\usepackage{subcaption}
\usepackage{authblk}
\usepackage{algorithm}
\usepackage{algpseudocode}

\usepackage{xcolor}

\usepackage[backend=biber, style=numeric, ]{biblatex}
\title{Reconnection of infinitely thin antiparallel vortices and coherent structures}

\author[1]{Sergei Iakunin}
\author[1,2]{Luis Vega}
\affil[1]{Basque Center for Applied Mathematics (BCAM)}
\date{}
\affil[2]{University of the Basque Country (UPV/EHU)}

\begin{document}

    \maketitle


\begin{abstract}
    One of the characteristic features of turbulent flows is the emergence of many vortices 
    which interact, deform, and intersect, generating a chaotic movement. The evolution of 
    a pair of vortices, e.g. condensation trails of a plane, can be considered as a basic element 
    of a turbulent flow. This simple example nevertheless demonstrates very rich
    behavior which still lacks a complete explanation. In particular, after the reconnection
    of the vortices some coherent structures with the shape of a horseshoe emerge. They have a high level of complexity generated by the interaction
    of waves running from the reconnection region. These structures also turn to be very reminiscent to the ones obtained 
    from the localized induction approximation applied to a polygonal vortex. It can be considered as an evidence that a 
    pair of vortices creates a corner singularity during the reconnection. In this work we focus on a
    study of the reconnection phenomena and the emerged structures. In order to do it we present a new model based on 
    the approximation of an infinitely thin vortex, which allows us to focus on 
    the chaotic movement of the vortex center line. The main advantage of the developed model 
    consists in the ability to go beyond the reconnection time and to see the coherent 
    structures. It is also possible to define the reconnection time by analyzing the fluid impulse.
\end{abstract}
\section{Introduction}\label{sec:introduction}

\begin{figure}
    \centering
    \begin{subfigure}[c]{0.75\textwidth}
        \centering
        \includegraphics[width=\textwidth]{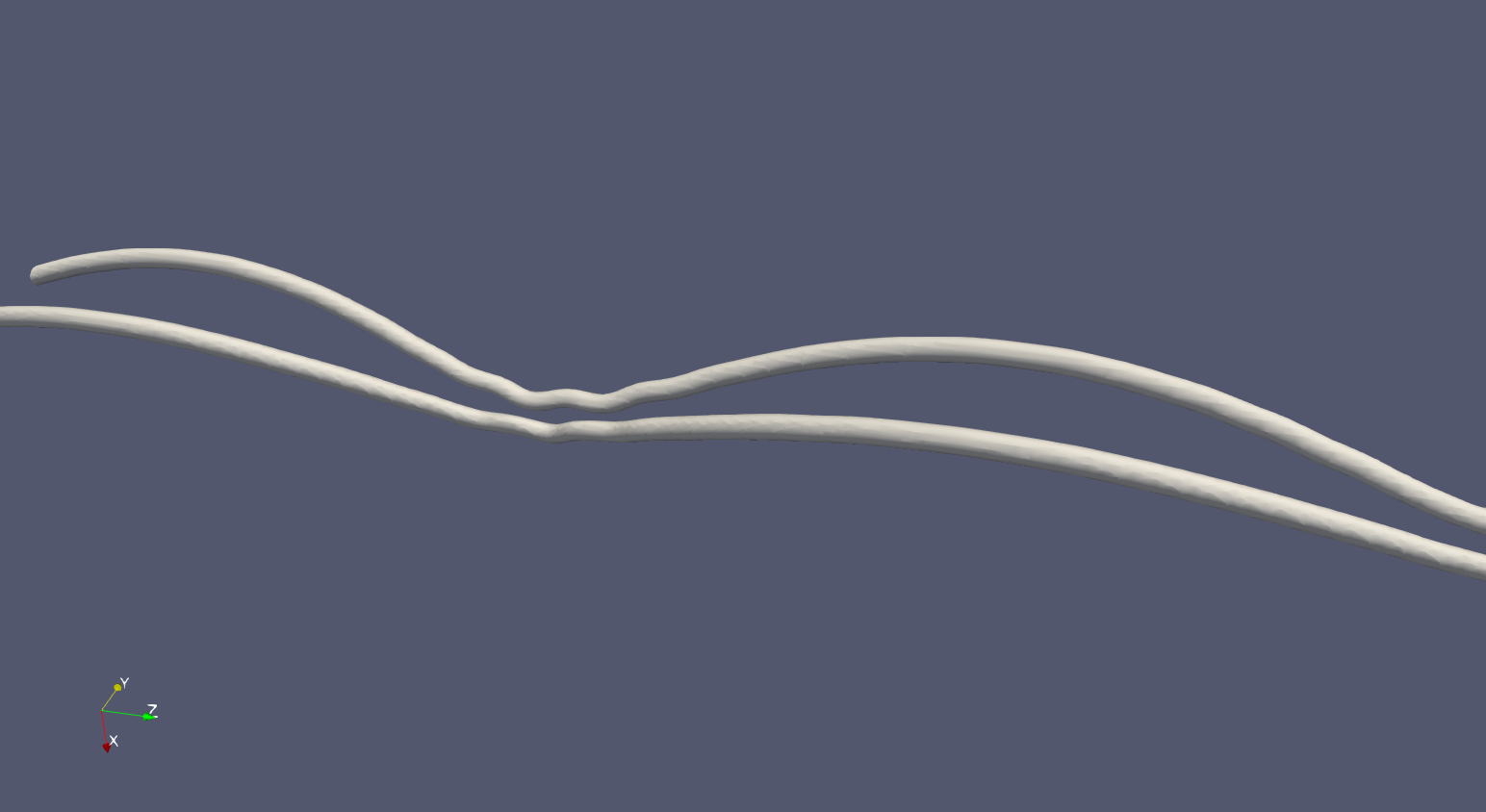}
        \caption{Crow waves (time $t = 5$)}\label{fig:ns_crow}
    \end{subfigure}
    \begin{subfigure}[c]{0.75\textwidth}
        \centering
        \includegraphics[width=\textwidth]{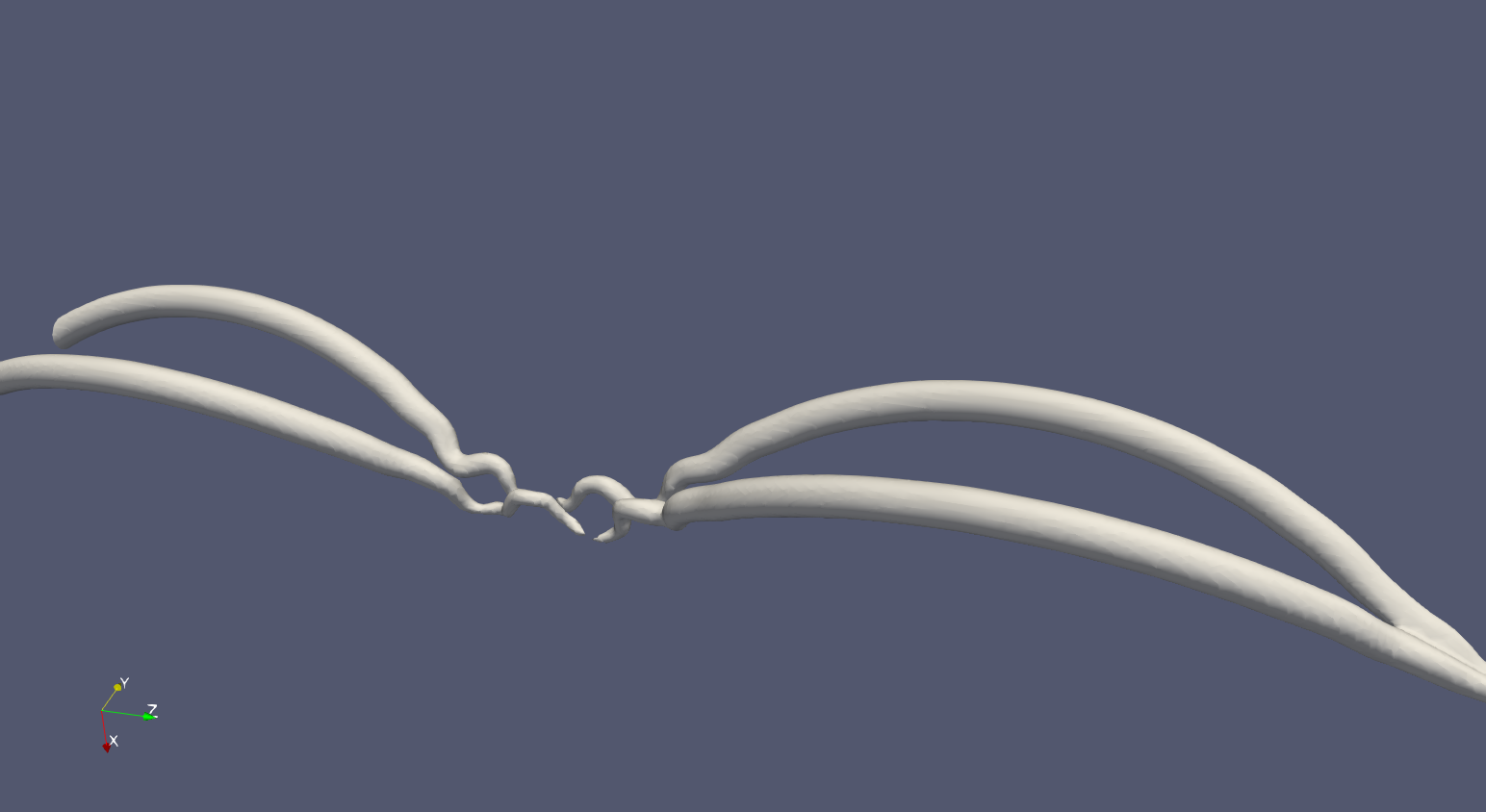}
        \caption{Kelvin waves (time $t = 6$)}\label{fig:ns_kelvin}
    \end{subfigure}
    \begin{subfigure}[c]{0.75\textwidth}
        \centering
        \includegraphics[width=\textwidth]{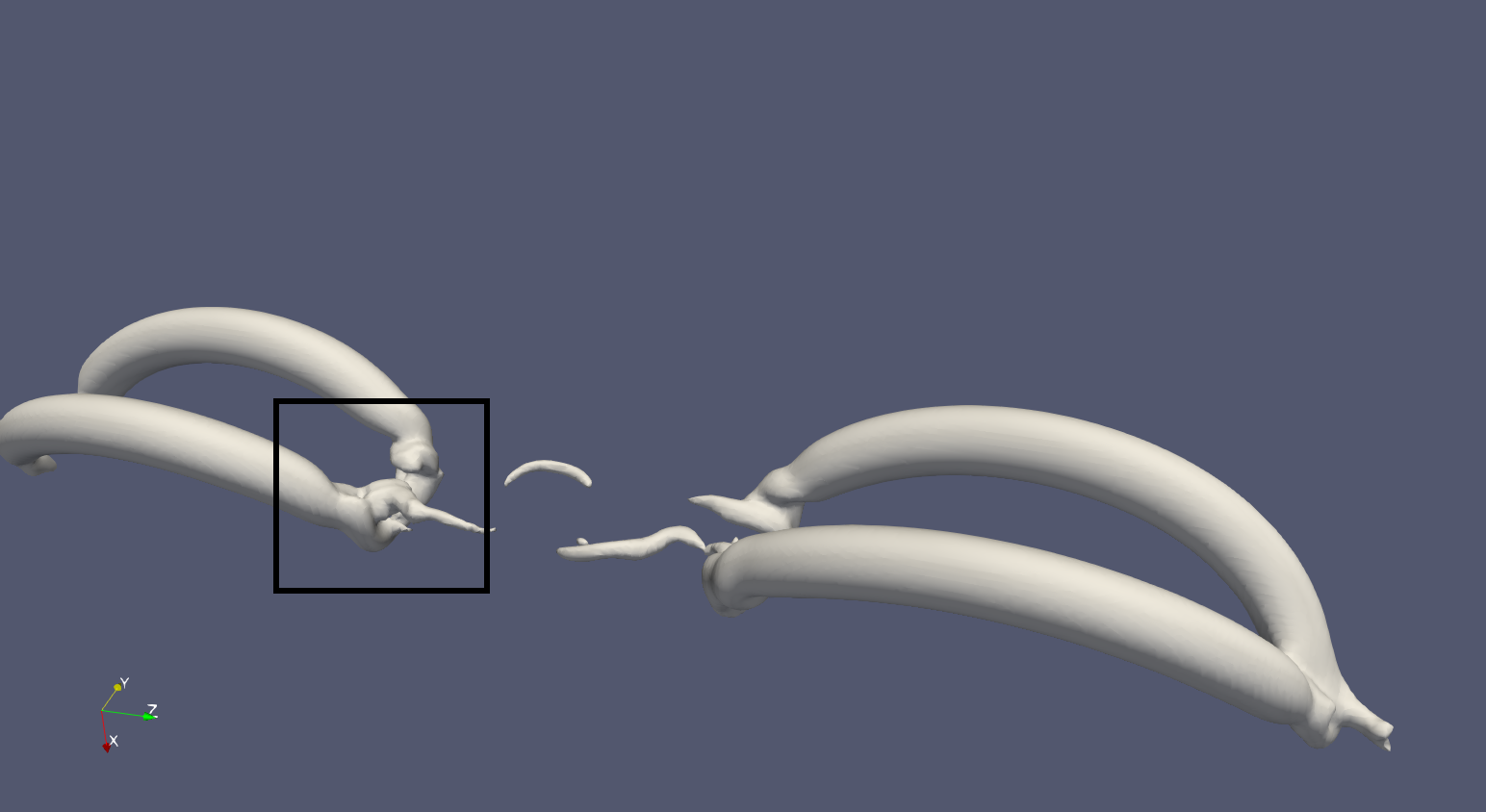}
        \caption{Horseshoes (one in the black rectangle and the symmetric one) and bridge between them (time $t = 7$)}\label{fig:ns_horseshoe}
    \end{subfigure}
    \caption{Simulation of vortex reconnection using LES in OpenFOAM in a periodic domain of nondimensionalized size
    $80\times 80\times 320$ discretized with $112\times 112\times 324$ elements. Vortex strength is $\Gamma = -1500$, vortex core radius is $r_c = 2$,
    and viscosity $\mu = 10^{-5}$, $Re \approx 6\cdot 10^8$. The visualisation is done with the $\lambda_2$ criterium.}\label{fig:navier_stokes}
\end{figure}

One possible way of transition from laminar to turbulent flow is the
interaction between vortices. This is a fascinating process
characterized by the reconnection phenomena when the topology of vortices 
changes producing a cascade of smaller structures which can also reconnect. 
Simple examples of such processes are the interaction of a pair of vortices
like a reconnection of aircraft condensation trials, considered by Crow in~\cite{crow1970}, 
or the collision of vortex rings~\cite{kida1988, lim1992}. In both cases the vortices firstly undergo
long wave deformation and then reconnect generating a series of
smaller rings or eye-shaped structures as those in the figure~\ref{fig:navier_stokes}. This sudden change of the flow structure is quite
impressive and still not completely understood. It is also interesting that
there is a surprising similarity between the statistical behavior of vortex filaments
in turbulent flows for quantum and classical fluids~\cite{bewley2008, nemirovskii2020}. Therefore, 
an explanation of this phenomena may be very useful for the understanding of turbulence. 

Vortices are regions in the fluid were a circular movement of particles happens. These regions
can be recognized by a high magnitude of the vorticity field which is the curl of the velocity.
However, the vorticity field moves with the flow obeying the Navier-Stokes equations, so the vortices
are also moving and deforming as if they were objects embedded into the flow. The fact that 
the vortices are part of a flow and also moved by it leads to that the detection and extraction
of vortices from the solutions of the Navier-Stokes equations is quite challenging, and it is
even more complicated to follow their evolution~\cite{pontin2018}.

We develop a new model of vortex interaction based on the approximation of infinitely thin vortex
embedded into a non-viscous fluid. This model allows us to analyze the behavior of vortices before and after 
the reconnection moment. After the reconnection time the model predicts the emergence of structures very reminiscent to those
obtained from the localized induction approximation (LIA) applied to an isolated eye-shaped vortex. The solution of
LIA is related to Riemann's non-differentiable function (RNDF, see formula~\eqref{eq:rndf}). The ability of the developed model to pass beyond 
the reconnection time and see the 
coherent structures allows to establish a relation between a classical mathematical object such as RNDF and the turbulent flows.

In the simplest cases of the reconnection it is enough to consider a pair of vortex tubes, that is,
cylindrical regions where the vorticity lines are parallel, and the vorticity magnitude is decaying 
far from the cylinder central line. Even in this case we have two different processes:
deformation of the vortex core and deformation of the central line. The first one leads 
to the creation of the helical Kelvin waves seen in many experiments~\cite{leweke2016}, 
numerical simulations~\cite{laporte2002}, and described in detail in~\cite{ledizes2005}. 
The presence of the Kelvin waves leads to a less clean reconnection process what
makes more challenging to understand the phenomena. Their amplitude, however, depends on the size of 
the vortex core so in this article we focus on infinitely thin vortices in order to avoid these waves. 
We will refer to this approximation  as the vortex filament approximation. 
In our model the asymmetric helical structures of Kelvin waves before the reconnection
disappear, but the symmetric ones persist. It can be considered as an indicator that these waves 
have different nature: the ones before are due to the deformation of the core, but the ones after
are due to the emerged singularity. 

Many attempts to describe the vortex reconnection were done in the last 50 years. It started with the 
pioneer paper of Crow~\cite{crow1970} where the evolution of aircraft 
condensation trails is studied. A pair of initially straight vortices undergoes slow,
nearly planar sinusoidal instability. However, the planes containing the deformed vortices
are different and incline one to another by approximately~$\pi / 4$. The amplitude 
of the deformation grows and after some time the vortices reconnect 
forming a train of vortex rings. The stability analysis is performed in~\cite{crow1970} under the assumption that
the vortices are infinitely thin, and the viscosity can be neglected. Suppose 
that $\mathbf{X}_i(s,t) \subset \mathbb{R}^3$ are curves defining central lines of both 
vortices. We can find the velocity $\mathbf{v}$ of the flow in any point $\mathbf{x}$ 
outside the vortices using the Biot-Savart integral:
\begin{equation}
    \mathbf{v}(\mathbf{x}) = \sum_{j = 1}^2 
    \frac{\Gamma_j}{4 \pi}\int_{-\infty}^\infty 
    \frac{(\mathbf{x}- \mathbf{X}_j(s,t))\wedge \frac{\partial}{\partial s}\mathbf{X}_j(s,t)ds}
    {|\mathbf{x}- \mathbf{X}_j(s,t)|^3}.\label{eq:intro_biot_savart}
\end{equation}
Here $t$ is time, $s$ is the parameter of the curve, $i \in \{1,2\}$, $\Gamma_i$ is the vortex strength, the symbol $\wedge$ 
defines the vector product,
the domain is supposed to be $\mathbb{R}^3$, and the vortices are infinite. The integrals~\eqref{eq:intro_biot_savart}
after the introduction of a cutoff can be also applied to a point on the vortex providing the velocity 
$\frac{\partial}{\partial t}\mathbf{X}_i(s,t)$. In~\cite{crow1970} a linear 
perturbation analysis is applied to a pair of initially antiparallel vortices finding that the most unstable mode is the long
symmetric wave, called Crow wave. There are also short and asymmetric waves seen in experiments and
predicted by the model, but their growing rate is much slower, and they probably emerge only under 
certain atmospheric conditions. The wave length predicted by the Crow model for aircraft condensation trails
equals $8.6 b$ where $b$ is the initial distance between vortices (aircraft wingspan). This result
slightly exceeds the wavelength obtained in experiments~\cite{ortega2003} and numerical simulations~\cite{han2000}.
It happens because the finite size core enhance the growth of shorter waves especially in the case of intense turbulence.

Even though the Crow model predicts instability it is hard to use it for the numerical simulation 
or any further analysis. A series of simpler models was proposed in several papers by Klein, Majda, and Damodaran 
\cite{klein1991a},\cite{klein1995}. Since there is only long wave deformation  we can choose an orthonormal basis
$\mathbf{e}_1$, $\mathbf{e}_2$, $\mathbf{e}_3$ and suppose
that the vortices are nearly parallel to $\mathbf{e}_3$:
\begin{equation}
    \mathbf{X}_i(s,t) = s \mathbf{e}_3 + 
    \delta^2 \mathbf{X}_i^{(2)}\left(\frac{s}{\delta}, \frac{t}{\delta^2}\right) + 
    o(\delta^2),\label{eq:intro_klein_assumption}
\end{equation}
where $\delta \ll 1$, and $\mathbf{X}_i^{(2)}(s,t)$ is always orthogonal to $\mathbf{e}_3$, so it can
be considered as a 2D vector.  The allowed wavelength is proportional to $\delta$, but they are
still long comparing to the initial distance between vortices that is proportional to $\delta^2$.
Under this assumption the Biot-Savart integral can be 
approximated up to the leading order in $\delta$:
\begin{equation}
    \frac{\partial}{\partial t}\mathbf{X}_i(s,t) = 
    \mathbf{J}\left(\alpha_i \Gamma_i \frac{\partial^2}{\partial s^2} \mathbf{X}_i(s,t) + 
    \sum_{i \ne j} 2 \Gamma_j \frac{\mathbf{X}_i(s,t) - \mathbf{X}_j(s,t)}
    {|\mathbf{X}_i(s,t) - \mathbf{X}_j(s,t)|^2}\right), 
    \label{eq:intro_klein_model}
\end{equation}
where $\alpha_i$ is a constant which depends on vortex core, $\Gamma_i$ is the vortex strength, and
\begin{equation}
    \mathbf{J} = \begin{pmatrix}
        0 & -1 \\ 1 & 0
    \end{pmatrix}.\nonumber
\end{equation}
There are two terms in~\eqref{eq:intro_klein_model}: the first one is the local self-induction
that is the velocity of the flow generated by the vortex itself, whereas the second one is
the velocity produced by the external flow generated by other vortices. The reconnection after finite time for this model
is proven in~\cite{banica2017}. The equations~\eqref{eq:intro_klein_model} a much easier to use for numerical simulation. However, due to 
the singularity in the second term it is impossible to go beyond the reconnection time.
Furthermore, it is not clear if the assumption~\eqref{eq:intro_klein_assumption} holds true when the 
amplitude of Crow waves is large and consequently the distance between vortices is small. 

Another possible approach to study the reconnection phenomena is to consider the Navier-Stokes
equations. In this case we do not have any problems related to the singularity and can 
include all the details such as compressibility, viscosity, and core deformation. The 
simulation, however, requires more computations and is more difficult. The direct numerical simulation (DNS)
of the incompressible Navier-Stokes equations is done in~\cite{hussain2011}. The visualisation
of vortices with the $\lambda_2$ criterion~\cite{hussain1995} shows the flattering of the vortex core near the reconnection
region and the formation of threads between vortex rings. These threads are stretching and 
may reconnect again if the Reynolds number is big enough. The presence of viscosity makes it difficult
to follow vortex lines due to the~dissipation. Therefore, the topology of the vortices is defined by 
surfaces not by the curves and is much more complicated to analyse. In~\cite{pontin2018} the 
vortex lines are extracted from the solution of the Navier-Stokes equations and classificated into
ones that reconnect, threads, and the additional vortex rings which emerges in the reconnection zone.
It again demonstrates that these phenomena contains a lot of different effects. 
Another attempt is done in~\cite{hussain2020} where the reconnection processes 
is divided in 3 stages: (i) the vortex cores flattering and stretching in the reconnection region; 
(ii) cutting and reconnection of the inner vortex lines that leads to formation of bridges; 
(iii) formation of threads from the rest vortex lines where the energy is dissipating through 
a turbulent cascade. More details are given in the review article~\cite{hussain2022}.
The evolution of threads in the last stage is very complicated and chaotic for high Reynolds
number, that is for thin vortices or almost non-viscous fluids. Thus, if we consider the Navier-Stokes equation 
for small Reynolds number then we have to deal with the deformation of the core and the dissipation due to viscosity
which does not allow us to extract main coherent structures. On the other hand for high Reynolds number we face the chaotic
behavior of the threads after reconnection which is also difficult to filter. We can also note that the 
coherent structures emerge not only in the reconnection region but also far from it as a result of the interaction of running waves.
These structures are clearer, have a distinctive horseshoe shape, and they are the focus of our attention. In particular, we highlight in
this paper that they have a behavior very reminiscent to the evolution of the corner vortex under LIA.

We show in figure~\ref{fig:navier_stokes} the main features that emerge due to the reconnection process. 
We perform the solution of the Navier-Stokes equations using the large eddy simulation (LES) in 
OpenFOAM software. In figure~\ref{fig:ns_crow} the symmetric large length Crow waves emerge. Further, in 
figure~\ref{fig:ns_kelvin} we can notice small length asymmetric Kelvin waves in the region close to the 
reconnection. And finally after the reconnection we can see in figure~\ref{fig:ns_horseshoe} one horseshoe structure 
in black rectangle and the bridge going to another symmetric one. In~\cite{hussain2011, hussain2020, hussain2022, brenner2016} the bridge is called thread
whereas the horseshoes are called bridges but in this article we will follow the introduced terminology. In this paper we are not interested in 
the vortex reconnection cascade \cite{melander1988} but in the evolution of the horseshoe structures and in the complexity of the interaction
of the waves that emerge. The extraction of these structures 
and the definition of the reconnection time is quite challenging due to the finite thickness 
of the vortices in the Navier-Stokes simulation.

The reconnection of vortices also happens in superfluids~\cite{bewley2008,fonda2019}. Even though, there are many differences between 
classical fluids and superfluids some features of the turbulent regime have similarities, e.g. classical 
turbulence flows have a filamentary structure~\cite{nemirovskii2020}. Therefore, the results in
superfluid reconnection are also important for us. Instead of classical vortices, in superfluids there are 
quantum vortices that are topological defects where the density tends to zero. Usually the quantum vortices are studied 
using Gross-Pitaevsky equations~\cite{krstulovic2017} or doing vortex filament approximation~\cite{schwarz1985} since the quantum 
vortices are infinitely thin. There is still a problem on how to model the reconnection and the change of the topology
of the vortices. In~\cite{schwarz1985} it is done by a heuristic way of measuring
the distance between nodes in the discretization of the filaments. The resulting shape demonstrates a self-similar behavior 
with a circular horseshoe and helical waves running along the reconnected vortices. 

\begin{figure}
    \centering
    \begin{subfigure}[c]{0.3\textwidth}
        \centering
        \includegraphics[width=\textwidth, trim={4cm 2cm 15mm 15mm}, clip]{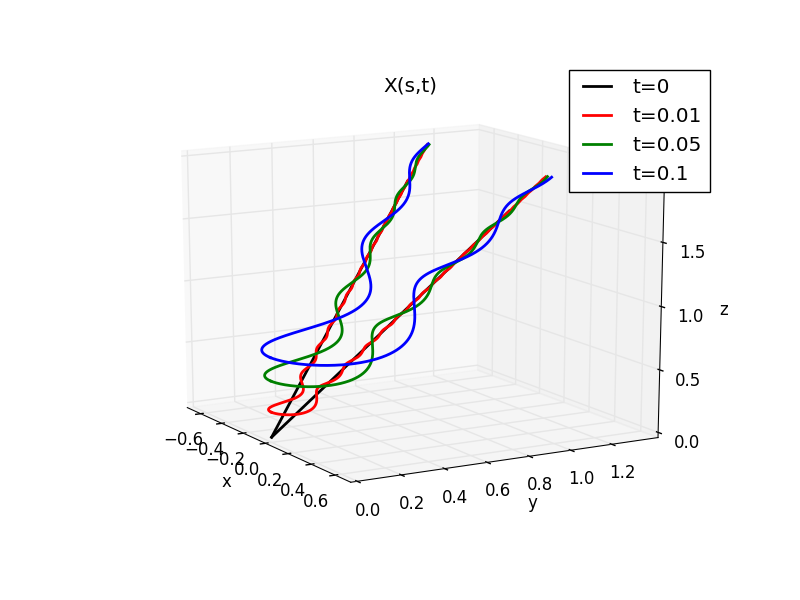}
        \caption{3D view}
    \end{subfigure}
    \begin{subfigure}[c]{0.3\textwidth}
        \centering
        \includegraphics[width=\textwidth, trim={4cm 2cm 15mm 15mm}, clip]{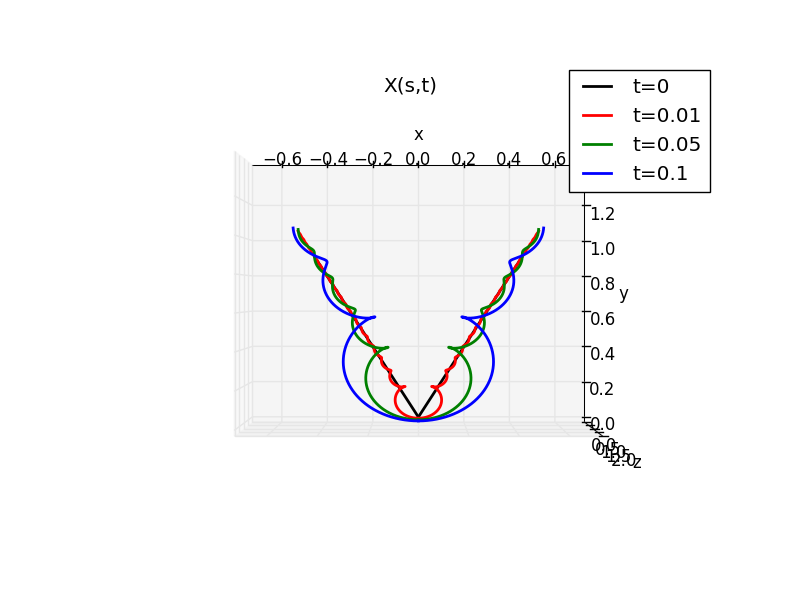}
        \caption{top view}
    \end{subfigure}
    \begin{subfigure}[c]{0.3\textwidth}
        \centering
        \includegraphics[width=\textwidth, trim={4cm 2cm 15mm 15mm}, clip]{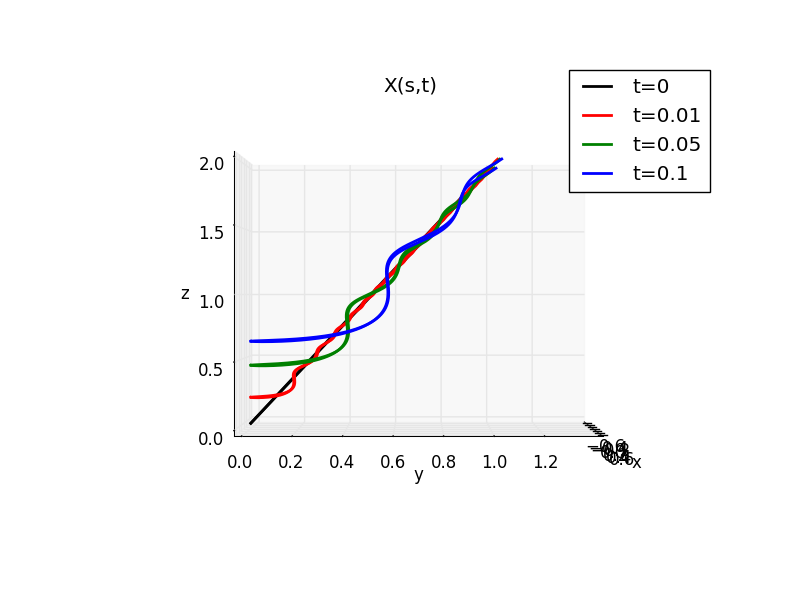}    
        \caption{side view}
    \end{subfigure}    
    \caption{Self-similar evolution of a corner vortex from different points of view}\label{fig:corner}
\end{figure}

The configuration with the horseshoe and helical waves is very reminiscent to the evolution of a vortex filament that moves
according to LIA and at the initial
time is given by two half-lines that meet at point (the corner) with an angle $\theta$~\cite{vega2003, delahoz2018, lipniacki2003}. We will call this
vortex the corner vortex. Its evolution is depicted in figure~\ref{fig:corner}.  The curve is self-similar and has constant curvature 
$c(s,t) = c_0 / \sqrt{2 t}$, and torsion $\tau(s,t) = s / t$, where $s$ is the arclength of the curve, $t$ is time, $c_0$ and $\theta$ satisfy 
the relation $\sin\frac{\theta}{2} = e^{-\pi c_0^2 / 2}$~(see~\cite{vega2003}). 
It is possible to see that the corner is turned into a circular horseshoe quite similar to the one we can
see in figure~\ref{fig:ns_horseshoe}. The parameter $c_0$ and the angle $\theta$ are estimated from the quantum vortex reconnection experiment 
in the superfluid helium $^4 He$ studied in~\cite{fonda2019}. The obtained data is consistent with the analytical results from~\cite{vega2003},
so one would expect the local induction approximation to work reasonably well, at least for quantum vortices.  One of the goals of this work is 
to generalize the model by adding the interaction term and to establish a relation between the self-similar behavior of an infinitely thin corner 
vortex and the reconnection of vortices.

In this article we present a new model based on the vortex filament approximation but with the interaction term similar
to the one in~\cite{klein1995}. This model is able to go beyond the reconnection moment in a natural way, without any heuristic, and generate
coherent structures. It seems to behave as the one showed in the figure~\ref{fig:corner} that we have just explained. The reconnection
time can not be found from the configuration of vortices, so we can not change the topology without introducing an error. However,
if we focus on an integral quantity the reconnection time becomes much more clear. We study the fluid impulse that is an integral of
a cross product of position and vorticity calculated around the reconnection region. Initially this integral changes monotonously but 
at some time starts to oscillate. We consider this time as the reconnection time. Furthermore, the behavior of the fluid impulse
after reconnection is very reminiscent to the RNDF. A similar effect for a polygonal vortex was
discussed in~\cite{delahoz2018}, and we can consider this as an evidence, that the vortices at the reconnection seem to create a corner 
similar to the one described in~\cite{vega2003}.

The paper has the following structure. In the section~\ref{sec:derivation} the derivation of the equations from the Biot-Savart law are shown. 
Next, in the section~\ref{sec:analysis} we discuss 
some properties of the new model and its relation to the previous ones. In the section~\ref{sec:numerics} we 
describe the numerical method we use to solve the derived equations. Finally, in the section~\ref{sec:results} 
we present results of the numerical simulation, and in the section~\ref{sec:eye-shaped} we compare the behavior after the reconnection 
with the one of an isolated vortex which has the shape of an eye. It can be used as an approximation of the vortex which emerges after 
the reconnection (see figure~\ref{fig:ns_horseshoe}) and also as analog of a curvilinear polygon with only two corners. In the section~\ref{sec:conclusions} and we 
make conclusions and discuss possible directions of further research. The appendix~\ref{sec:frenet} is devoted to an alternative formulation and possible simplification,
the appendix~\ref{sec:algorithm} contains the numerical algorithm, and in the appendix~\ref{sec:eye-shaped-details} the evolution of the eye-shaped vortex is described.

\section{Derivation}\label{sec:derivation}

The velocity of the flux produced by a pair 
of infinitely thin antiparallel vortices is 
given by the Biot-Savart law:
\begin{equation}
\mathbf{v}(\mathbf{x})\!=\!-\frac{\Gamma}{4\pi}\left(\!
\int_{-\infty}^{\infty} \frac{(\mathbf{x}\!-\!\mathbf{X}(s))\wedge 
\frac{\partial\mathbf{X}}{\partial s}(s)}
{|\mathbf{x}\!-\!\mathbf{X}(s)|^3} ds\!-\!\int_{-\infty}^{\infty} 
\frac{(\mathbf{x}\!-\!\mathbf{Y}(s))\wedge  \frac{\partial\mathbf{Y}}{\partial s}(s)}
{|\mathbf{x}\!-\!\mathbf{Y}(s)|^3}ds\!\right)\!, \label{eq:biot-savart}
\end{equation}
where $\Gamma$ is the circulation, $\mathbf{X}(s)$ 
and $\mathbf{Y}(s)$ are curves in $\mathbb{R}^3$ defining 
central lines of both vortices,  the symbol $\wedge$ 
defines the vector product.

The vortices are moving by the flow generated by them, therefore $\mathbf{X}(s,t)$ and  $\mathbf{Y}(s,t)$
are also functions of time. For a point $\mathbf{x}$ belonging to the vortex 
$\mathbf{X}(s, t)$ the first integral represents the
velocity due to the local self-induction $\mathbf{v}_{lia}$, 
and the second integral is the velocity of the external 
flow $\mathbf{v}_{ext}$ produced by the vortex $\mathbf{Y}(s,t)$. 
Then, we can decompose the velocity of the 
vortex filament into a sum:
$$
\frac{\partial}{\partial t}\mathbf{X}(s,t) = 
\mathbf{v}_{lia}(s,t) + \mathbf{v}_{ext}(s,t).
$$

\begin{figure}
    \centering
    \begin{tikzpicture}
        \draw[black, dashed] (0, -2.2) -- (0, 2.2);
        \draw[black, <->] (0, -1.7) -- (1.5, -1.7) node[pos=0.5, anchor=north] {$b$};
        \draw[black, ->] (0, 0) -- (1, 0) node[pos = 1, anchor=west] {$\mathbf{e}_1$};
        \draw[black, ->] (0, 0.) -- (0, 1) node[pos = 1, anchor=west] {$\mathbf{e}_3$};
        \draw[black, ->] (0, 0) -- (-0.7, -0.7) node[pos = 1, anchor=west] {$\mathbf{e}_2$};
        \draw[line width = 0.5mm, blue] (1.5, -2) -- (1.5, 2) node[pos = 1, anchor=west] {$\mathbf{X}(s,t)$};
        \draw[dashed, line width = 0.5mm, red] (-1.5, -2) -- (-1.5, 2) node[pos = 1, anchor=east] {$\mathbf{Y}(s,t)$};;
        \draw[line width = 0.5mm, blue, ->] (1.9, 1) arc(0:180:0.4);
        \draw[dashed, line width = 0.5mm, red, <-] (-1.1, 1) arc(0:180:0.4);
    \end{tikzpicture}    
    \caption{Initial configuration of vortices}\label{fig:vortex_init}
\end{figure}
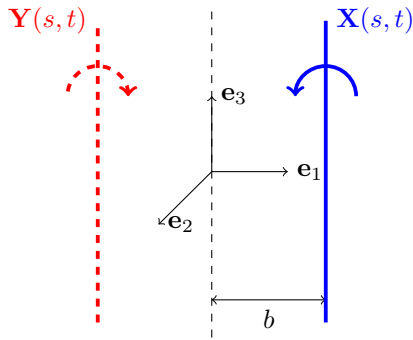
        
Let us call the components of vector 
$\mathbf{X}(s,t) = \begin{pmatrix} x_1(s,t) & x_2(s,t) & x_3(s,t) \end{pmatrix}^T$.
The initial configuration of the vortices is depicted in figure~\ref{fig:vortex_init}. Now we make the following assumptions: 
\begin{enumerate}
\item The vortices are symmetric respect to the plane $x_1 = 0$ 
hence we can reduce our problem to only one unknown 
curve $\mathbf{X}(s,t)$ obtaining the second one by:
\begin{equation}
\mathbf{Y}(s,t) = \begin{pmatrix}
-x_1(s,t) & x_2(s,t) & x_3(s,t)
\end{pmatrix}^T.\label{eq:def_of_y}
\end{equation}

\item For the velocity of the external flow we will 
use Rosenhead regularization~\cite{rosenhead1930}:
\begin{equation}
\mathbf{v}_{ext}(s,t) = \frac{\Gamma}{4\pi}\int_{-\infty}^{\infty} 
\frac{(\mathbf{X}(s,t)\!-\!\mathbf{Y}(q,t))\wedge 
\frac{\partial\mathbf{Y}}{\partial q}(q,t)}
{{\left(4r_c^2 + |\mathbf{X}(s,t)\!-\!\mathbf{Y}(q,t)|^2\right)}^
{\frac{3}{2}}} dq,\label{eq:rosenhead_regularization}
\end{equation}
where $r_c$ is a parameter related to the vortex core, $s$ is the parameter of the curve, and 
$q$ is the parameter of the curve used inside the integral. This regularization
prevents the singularity in the interaction term and can be understood as a 
viscosity effect during the merging of the cores of the vortices. For the 
self-induction part we do not use the regularization because there is no
core merging.

\item The second vortex can be linearized near any point $s$
that is in the interval $q\in (s - \alpha\sqrt{x_1^2(s,t) + r_c^2}, s + \alpha\sqrt{x_1^2(s,t) + r_c^2})$ 
for some parameter $\alpha$, furthermore the tails can be neglected in the Biot-Savart law. The length of 
this interval is almost proportional to the distance between vortices. When the second vortex is far, 
a large part of it makes a relevant contribution. On the other hand, when it is close the situation is similar to
LIA because the second integral in~\eqref{eq:biot-savart} is close to singular, thus we 
can consider only a small leading piece of it.

\item Either $x_1(s,t)$ or $\frac{\partial x_1(s,t)}{\partial s}$ are small so the product
$x_1(s,t) \frac{\partial x_1(s,t)}{\partial s}$ is neglectable. It means that a piece of vortex can be 
oriented in $\mathbf{e}_1$ direction only when it is close to the second vortex.

\end{enumerate}

\paragraph*{Self-induction}
Let us start from the first term in~\eqref{eq:biot-savart} which corresponds 
to the self-induction. Here we follow the standard derivation of the vortex 
filament equation for binormal flow~\cite{shaffman1992} using LIA. Fix a time moment and select a point 
$\mathbf{X}(s)$. We choose a perturbation 
$\rho \boldsymbol{\xi}(s) = \rho \left(\mathbf{N}(s) \cos(\theta) + 
\mathbf{B}(s)\sin(\theta)\right)$ with small $\rho$, normal $\mathbf{N}(s)$, 
binormal $\mathbf{B}(s)$, and  some angle $\theta$. 
The self-induced velocity of the vortex is found as the limit of Biot-Savart integral 
when $\rho$ goes to zero. Since the expression in the integral becomes singular we can estimate it 
using only a neighborhood of the point $s$ with a cut-off $L_{lia}$ where we can 
expand $\mathbf{X}(q)$ in the Taylor series up to the third order of $|q-s|$:
\begin{multline}
\mathbf{v}_{lia}(s) \approx -\frac{\Gamma}{4\pi}\int_{s-L_{lia}}^{s + L_{lia}} 
\frac{(\mathbf{X}(s)\!+\!\rho \boldsymbol{\xi}(s)\!-\!\mathbf{X}(q))\wedge \mathbf{X}_q(q)} 
{|\mathbf{X}(s)\!+\!\delta \mathbf{X}(s)\!-\!\mathbf{X}(q)|^3} dq\\ 
\approx -\frac{\Gamma}{4\pi}\int_{s-L_{lia}}^{s + L_{lia}}\Biggl( 
\frac{\rho \boldsymbol{\xi}(s)\wedge\bigl(\mathbf{X}_s(s)\!+\!(q\!-\!s)\mathbf{X}_{ss}(s)\bigr)} 
{{\left(\rho^2\!+\!{(q\!-\!s)}^2|\mathbf{X}_s(s)|^2\right)}^\frac{3}{2}} \\- 
\frac{{(q\!-\!s)}^2\mathbf{X}_s(s) \wedge \mathbf{X}_{ss}(s)} 
{2 {\left(\rho^2\!+\!{(q\!-\!s)}^2|\mathbf{X}_s(s)|^2\right)}^\frac{3}{2}} \Biggr) dq, \label{eq:si_expansion} 
\end{multline}
subindeces here designate corresponding derivatives. We have also used that $| \boldsymbol{\xi}(s)| = 1$ 
and that it is orthogonal to $\mathbf{X}_s(s)$. The first term in~\eqref{eq:si_expansion} represents 
rotation around the vortex central line without its alteration whereas the second one 
is the movement of the central line that gives the self-induced velocity:
\begin{equation}
\mathbf{v}_{lia} \approx \frac{\Gamma}{4\pi}\frac{\mathbf{X}_s\wedge \mathbf{X}_{ss}}{2 |\mathbf{X}_s|^3} 
\left(-\frac{2|\mathbf{X}_s|}{\sqrt{|\mathbf{X}_s|^2 + 4\frac{\rho^2}{L_{lia}^2}}} + 
\ln\left(\frac{{\sqrt{|\mathbf{X}_s|^2 + 4\frac{\rho^2}{L_{lia}^2}} + |\mathbf{X}_s|}} 
{{\sqrt{|\mathbf{X}_s|^2 + 4\frac{\rho^2}{L_{lia}^2}}-|\mathbf{X}_s|}}\right)\right).\nonumber 
\end{equation} 
In the limit $\rho/L_{lia} \to 0$ we obtain the local induction approximation (or the binormal flow):
\begin{equation}
\mathbf{v}_{lia}(s) = \frac{\Gamma}{4\pi}\left(-1 + \ln\left(\frac{L_{lia}}{\rho}|\mathbf{X}_s(s)|\right)\right) 
\frac{\mathbf{X}_s(s)\wedge \mathbf{X}_{ss}(s)}{|\mathbf{X}_s(s)|^3}.\label{eq:vfe} 
\end{equation}
Assume that the cut-off is inverse to the modulus of tangential vector, 
that is $L_{lia} = \tilde{L}_{lia} / |\mathbf{X}_s(s)|$, and introduce the first 
parameter of our model:
\begin{equation}
\varepsilon = \frac{2}{-1 + \ln\left(\tilde{L}_{lia}/\rho\right)}.\label{eq:epsilon}
\end{equation} 
We will see later that the modulus of the tangential vector is growing when the time is close 
to the reconnection moment. However, due to the regularization in the interaction term we can 
bound it from above with a power of $r_c$. The parameter $\varepsilon$ represents the strength
of the interaction between vortices.

\paragraph*{External flow}
The second integral in~\eqref{eq:biot-savart} after applying symmetry assumption, 
Rosenhead regularization~\eqref{eq:rosenhead_regularization}, and fixing the time reads
\begin{equation}
\mathbf{v}_{ext}(s)\!=\!\frac{\Gamma}{4\pi}\!\int_{-\infty}^{\infty}\!
\frac{(2 x_1(s)\mathbf{e}_1 + \mathbf{Y}(s) - \mathbf{Y}(q))\wedge \mathbf{Y}_q(q)} 
{{\left(4 r_c^2\!+\!4 x_1^2(s)\!+\!4 x_1(s) (y_1(s)\!-\!y_1(q))\!+ 
\!|\mathbf{Y}(s) \!-\!\mathbf{Y}(q)|^2\right)}^{\frac{3}{2}}}dq. \label{eq:external_flow} 
\end{equation}
where component $y_1(s)$ of vector $\mathbf{Y}(s)$ is given by~\eqref{eq:def_of_y}. 
We can apply the assumption 3 obtaining: 
\begin{equation}
\mathbf{v}_{ext}(s)\!\approx\!\frac{\Gamma}{32\pi {(x_1^2(s)\!+\!r_c^2)}^{3/2}}\!
\int_{s\!-\!\alpha\sqrt{x_1^2(s)\!+\!r_c^2}}^{s\!+\!\alpha\sqrt{x_1^2(s)\!+\!r_c^2}} 
\frac{2 x_1(s) \mathbf{e}_1 \wedge \mathbf{Y}_s(s) dq} 
{{\left(1 \!-\!\frac{(q\!-\!s) x_1(s) x_{1,s}(s) \!-\!{(q\!-\!s)}^2 |\mathbf{Y}_s(s)|^2} 
{x_1^2(s)\!+\!r_c^2}\right)}^{3/2}}. \nonumber 
\end{equation}
Neglecting the term $x_1(s) x_{1,s}(s)$ due to the assumption 4 we obtain that the external 
velocity is
\begin{equation}
\mathbf{v}_{ext}(s)\!=\!\frac{\Gamma x_1(s) \mathbf{e}_1\wedge\mathbf{Y}_s(s)} 
{16\pi (x_1^2(s)\!+\!r_c^2) |\mathbf{Y}_s(s)|} 
\frac{2\alpha |\mathbf{Y}_s(s)|}{\sqrt{1\!+\!\alpha^2 |\mathbf{Y}_s(s)|^2}}\!\approx\! 
\frac{\Gamma x_1(s) \mathbf{e}_1\wedge\mathbf{Y}_s(s)} 
{8\pi (x_1^2(s)\!+\!r_c^2) |\mathbf{Y}_s(s)|},\label{eq:ext_velocity} 
\end{equation}
if we choose $\alpha$ large enough. 

\paragraph*{Main equations}
Summing up~\eqref{eq:vfe} with~\eqref{eq:ext_velocity} and
rescaling the time with $[t] = 8\pi\varepsilon / \Gamma$ we obtain the main equations 
describing the evolution of a pair of symmetric vortices due to self-induction and interaction:
\begin{equation}
\mathbf{X}_t = \frac{\mathbf{X}_s \wedge \mathbf{X}_{ss}}{|\mathbf{X}_s|^3} - 
\frac{\varepsilon x_1}{x_1^2 + r_c^2} \frac{\mathbf{X}_s\wedge \mathbf{e}_1 }{|\mathbf{X}_s|}.\label{eq:main_equation}
\end{equation}
In~\eqref{eq:main_equation} we have also used that $|\mathbf{X}_s| = |\mathbf{Y}_s|$ 
and $\mathbf{X}_s \wedge  \mathbf{e}_1 = \mathbf{Y}_s  \wedge  \mathbf{e}_1$.
The equations should be equipped with an initial condition $\mathbf{X}_0(s) = \mathbf{X}(s,0)$
which is usually a small perturbation of a line $b\mathbf{e}_1 + s \mathbf{e}_3$ shifted
from the origin by a positive value $b$ in $\mathbf{e}_1$ direction and oriented in $\mathbf{e}_3$,
and the boundary conditions which we will suppose periodic on an interval $s\in (0, S)$.

The parameter $\varepsilon$ represents the strength of the vortex interaction when compared with the
self-induction. The bigger it is, the faster the reconnection happens. A more detailed relation between $\varepsilon$ and the 
velocity of the vortices is shown in 
section~\ref{sec:analysis}. In that section we also show how $\varepsilon$ influences on the vortex stretching. 
The parameter $r_c$ is necessary in order to avoid the singularity in the interaction term when $x_1 = 0$. However,
it has a physical meaning of viscosity. According to~\cite{shaffman1992} the radius of the vortex core $r_c \sim \sqrt{\nu t}$,
where $\nu$ is viscosity and $t$ is time. Since we are mainly interested in non-viscous reconnection the parameter
$r_c$ should be as small as possible. Even though the presence of $r_c$ does not allow to see a sharp corner we still 
can see the effect and complexity of wave interaction at later times.


\section{Some properties of the derived equations}\label{sec:analysis}

\paragraph*{Relation to previous models}

The Klein-Majda system of equations \cite{klein1995} for a pair of nearly parallel 
counter-rotating vortices in the symmetric case can be obtained from the 
equations~\eqref{eq:main_equation}. Indeed, in the considered case, taking into 
account that the vortices are nearly parallel to $\mathbf{e}_3$ and 
including regularization, the Klein-Majda system reads:
\begin{eqnarray}
&&\frac{\partial \mathbf{X}}{\partial t} = \frac{\Gamma}{4\pi} \mathbf{e}_3\wedge 
\left(\sigma\frac{\partial^2 \mathbf{X}}{\partial s^2} - 
\frac{\mathbf{X} - \mathbf{Y}}{|\mathbf{X} - \mathbf{Y}|^2 + r_c^2}\right),\nonumber\\ 
&&\frac{\partial \mathbf{Y}}{\partial t} = -\frac{\Gamma}{4\pi} \mathbf{e}_3 \wedge 
\left(\sigma\frac{\partial^2 \mathbf{Y}}{\partial s^2} - 
\frac{\mathbf{Y} - \mathbf{X}}{|\mathbf{Y} - \mathbf{X}|^2 + r_c^2}\right).\label{eq:klein_majda} 
\end{eqnarray}
Parameter $\sigma$ here depends on the structure of the vortex core, $\mathbf{X}$ and $\mathbf{Y}$ 
here are 2-dimensional vectors, the third component $x_3(s,t) = y_3(s,t) = s$ is known 
and ignored in the system. We can reduce the number of equations using symmetry:
\begin{equation}
x_{1,\tilde{t}} = -x_{2,ss},\quad x_{2,\tilde{t}} = x_{1,ss} - 
\varepsilon \frac{x_1}{x_1^2 + r_c^2},\label{eq:km_reduced}
\end{equation}
where $\tilde{t} = \frac{\Gamma}{4\pi} \sigma t$ is the rescaled time and $\varepsilon = 1 / \sigma$. 
Now we will consider the equation~\eqref{eq:main_equation} supposing that the vortex 
central line is given in the shape
\begin{equation}
\mathbf{X}(s,t) = s\mathbf{e}_3 + a \mathbf{X}^{(2)}\left(\xi s, \tau t\right),\label{eq:long_waves}
\end{equation}
where $\mathbf{X}^{(2)} \cdot \mathbf{e}_3 = 0$ and $\xi a \ll 1$. It means that similarly 
to~\eqref{eq:km_reduced} we have only 2 unknowns, and that the vortices may deviate 
from the straight line in only a long wave shape, comparing with the distance between vortices. 
In the Klein-Majda paper~\cite{klein1991a} these waves are called short-waves. However, if we compare them with 
the distance between vortices they are long. Plugging~\eqref{eq:long_waves} into~\eqref{eq:main_equation} 
and cancelling the amplitude $a$ we obtain:
\begin{equation}
\tau \mathbf{X}^{(2)}_t = \xi^2 \mathbf{e}_3 \wedge \mathbf{X}_{ss}^{(2)} - 
\varepsilon \frac{x_1^{(2)}}{r_c^2 + x_1^2} \mathbf{e}_2 + O(a\xi).\nonumber
\end{equation}
We suppose that the left-hand side and the first two terms of the right-hand side are of the same order 
whereas the rest is smaller, thus multiplying by $a$ we get:
\begin{equation}
\mathbf{X}_t = \mathbf{e}_3 \wedge \mathbf{X}_{ss} - 
\varepsilon \frac{x_1}{r_c^2 + x_1^2} \mathbf{e}_2,\nonumber
\end{equation}
that is equivalent to~\eqref{eq:km_reduced}.

\paragraph*{Crow instability}
Linear stability analysis of~\eqref{eq:main_equation} predicts a long-wave instability described by Crow 
in~\cite{crow1970}. Suppose that initially the vortices are parallel to $\mathbf{e}_3$, as depicted in figure~\ref{fig:vortex_init}, and add a perturbation:
\begin{equation}
    \mathbf{X}(s,t) = b\mathbf{e}_1 + v t \mathbf{e}_2 + s \mathbf{e}_3 + \delta e^{\mu t}\begin{pmatrix}
        \alpha \cos{\omega s} \\ \beta \cos{\omega s} \\ \gamma \sin{\omega s}
    \end{pmatrix} + O(\delta^2),\label{eq:perturbed_line}
\end{equation}
where $b$ is a half of the initial distance between vortices, $\delta \ll 1$. Without the perturbation 
the pair of vortices will move in $\mathbf{e}_2$ direction with velocity
\begin{equation}
    v = -\varepsilon \frac{b}{b^2 + r_c^2},\label{eq:velocity}
\end{equation}
that is proportional to $\varepsilon$ and almost inverse to the distance between vortices. This result 
coincide with many previous researches~\cite{crow1970}, with experiments, and with numerical simulation 
using the Navier-Stokes equations. Moreover, here we get another physical meaning of the parameter $\varepsilon$:
the bigger it is the faster the pair of vortices moves in the $\mathbf{e}_2$ direction.

Now let us find the frequencies of perturbations for which this straight line solution is not stable. 
Plugging~\eqref{eq:perturbed_line} into~\eqref{eq:main_equation} and keeping only linear terms respect 
to $\delta$ we obtain
\begin{equation}
    \mu \begin{pmatrix}
        \alpha \cos{\omega s} \\ \beta \cos{\omega s} \\ \gamma \sin{\omega s}
    \end{pmatrix} = -\omega^2 \begin{pmatrix}
        -\beta \cos{\omega s} \\ \alpha \cos{\omega s} \\ 0
    \end{pmatrix} - \varepsilon \alpha \mathbf{e}_2 \frac{r_c^2 - b^2}{{\left(b^2 + r_c^2\right)}^2} \cos{\omega s} -  
    \frac{\varepsilon \beta \mathbf{e}_3  b}{b^2 + r_c^2} \sin{\omega s}.\nonumber
\end{equation}
We have an eigenvalue problem
\begin{equation}
    \mu \begin{pmatrix}
        \alpha \\ \beta \\ \gamma 
    \end{pmatrix} = \begin{pmatrix}
        0 & \omega^2 & 0 \\ 
        -\omega^2 - \varepsilon \frac{r_c^2 - b^2}{{\left(b^2 + r_c^2\right)}^2} & 0 & 0\\
        0 & -\frac{\varepsilon  b}{b^2 + r_c^2} & 0
    \end{pmatrix}\begin{pmatrix}
        \alpha \\ \beta \\ \gamma 
    \end{pmatrix},\label{eq:eigen_value_problem}
\end{equation}
and the perturbed solution~\eqref{eq:perturbed_line} is unstable if at least one eigenvalue 
of~\eqref{eq:eigen_value_problem} has positive real part. It happens for the following 
frequencies~$\omega$ and wavelengths~$\lambda~=~2\pi/\omega$:
\begin{equation}
    \omega < \frac{\sqrt{\varepsilon (b^2 - r_c^2)}}{b^2 + r_c^2},\quad
    \lambda > \frac{2 \pi b}{\sqrt{\varepsilon}}\left(\frac{1 + {(r_c / b)}^2}{\sqrt{1 - {(r_c / b)}^2}}\right).
    \label{eq:crow_waves}
\end{equation}
Since $r_c$ represents the radius of the vortex core and is always smaller than the initial distance 
between vortices the square root in the expression~\eqref{eq:crow_waves} is always real. These waves are long and
called Crow waves since they were firstly described in~\cite{crow1970}. 

It is important to note that since we consider the equations~\eqref{eq:main_equation} on the interval 
$s\in(0, S)$ with periodic boundary conditions we have to be sure that the Crow waves~\eqref{eq:crow_waves}
fit in this interval, that is $S \ge \lambda$. It is also interesting that the velocity~\eqref{eq:velocity}
depends on $\varepsilon / b$ whereas the wavelength~\eqref{eq:crow_waves} depends on $\sqrt{\varepsilon} / b$
therefore we can not reduce number of parameters and consider only the ratio. Further we will see that $\varepsilon$
affects not only on the speed of the reconnection but also on the angle the vortices make at that moment. 

\paragraph*{The modulus of the tangential vector}
The vortex filament equations preserve the modulus of the tangential vector $\mathbf{T} = \mathbf{X}_s$. However,
when we have the interaction term as in~\eqref{eq:main_equation} it is not true anymore. Nevertheless, we can derive
a closed expression for the modulus $|\mathbf{T}|$. In order to do it we take a derivative of~\eqref{eq:main_equation}
respect to $s$ and calculate the inner product with $\mathbf{T}$:
\begin{equation}
    |\mathbf{T}|\frac{\partial}{\partial t} |\mathbf{T}| = 
    -\varepsilon \frac{x_1}{x_1^2 + r_c^2}\frac{\left(\mathbf{T}_s \wedge \mathbf{e}_1\right)\cdot \mathbf{T}}{|\mathbf{T}|} = 
    -\varepsilon \frac{x_1}{x_1^2 + r_c^2}\frac{\left( \mathbf{T} \wedge \mathbf{T}_s\right) \cdot \mathbf{e}_1}{|\mathbf{T}|}. 
    \label{eq:T_t_dot_T}
\end{equation}
We can also calculate the inner product of $\mathbf{X}_t$ with  $\mathbf{e}_1$ simplifying the expression of the 
right-hand side of~\eqref{eq:T_t_dot_T}:
\begin{equation}
    x_{1,t} = \frac{\left( \mathbf{T} \wedge \mathbf{T}_s\right) \cdot \mathbf{e}_1}{|\mathbf{T}|^3}.\label{eq:X_t_dot_e1}
\end{equation}
Combining~\eqref{eq:T_t_dot_T} with~\eqref{eq:X_t_dot_e1} and integrating respect to time we find the expression 
for the modulus of the tangential vector:
\begin{equation}
    |\mathbf{T}(s,t)| = L_0(s) {\left(x_1^2(s,t) + r_c^2\right)}^{-\varepsilon / 2},\label{eq:arclength}
\end{equation}
where $L_0(s)$ is a function which does not depend on time and is given by the initial conditions. The modulus $\mathbf{T}$
is growing in the reconnection region (that is $x_1$ goes to $0$) and even tends to a singularity when $r_c$ tends to $0$. It
can be understood as a vortex stretching phenomenon.

\paragraph*{Self-similar solution}
When $r_c$ goes to zero the equations~\eqref{eq:main_equation} have self-similar solutions. Let us define
$\eta = s / \sqrt{t}$ and plug $\mathbf{X}(s,t) = \sqrt{t}\mathbf{G}(\eta)$ into~\eqref{eq:main_equation} 
assuming that $r_c = 0$:
\begin{equation}
    \frac{1}{2 \sqrt{t}}\mathbf{G}(\eta) - \frac{\eta}{2\sqrt{t}}\mathbf{G}'(\eta) = 
    \frac{\mathbf{G}'(\eta)\wedge\mathbf{G}''(\eta)}{\sqrt{t}|\mathbf{G}'(\eta)|^3} - 
    \frac{\varepsilon}{\sqrt{t} G_1(\eta)} \frac{\mathbf{G}'(\eta)\wedge\mathbf{e_1}}{|\mathbf{G}'(\eta)|}.\nonumber
\end{equation}
It is easy to see that after multiplying by $\sqrt{t}$ we get an ODE for $\mathbf{G}(\eta)$:
\begin{equation}
    \frac{1}{2}\mathbf{G}(\eta) - \frac{1}{2}\eta\mathbf{G}'(\eta) = 
    \frac{\mathbf{G}'(\eta)\wedge\mathbf{G}''(\eta)}{|\mathbf{G}'(\eta)|^3} - 
    \frac{\varepsilon}{G_1(\eta)} \frac{\mathbf{G}'(\eta)\wedge\mathbf{e_1}}{|\mathbf{G}'(\eta)|}.\label{eq:self-similar-start}
\end{equation}
In order to extract the highest derivative we can calculate the cross product of~\eqref{eq:self-similar-start} with
$\mathbf{G}'(\eta)$:
\begin{equation}
    \frac{1}{2}\mathbf{G}\wedge \mathbf{G}' = \frac{\mathbf{G}''}{|\mathbf{G}'|} - 
    \frac{\mathbf{G}''\cdot\mathbf{G}'}{|\mathbf{G}'|^3}\mathbf{G}' - 
    \frac{\varepsilon |\mathbf{G}'|}{G_1}\mathbf{e}_1 + 
    \frac{\varepsilon G_1'}{|\mathbf{G}'| G_1}\mathbf{G}'.\nonumber
\end{equation}
Observe that $\mathbf{G}''\cdot\mathbf{G}' = |\mathbf{G}'|\frac{d}{d\eta}|\mathbf{G}'|$, and it can be expressed by
lower derivatives similarly to the previous paragraph. Thus, the final equation is:
\begin{equation}
    \mathbf{G}'' = |\mathbf{G}'|\frac{1}{2}\mathbf{G}\wedge \mathbf{G}' + \frac{\varepsilon |\mathbf{G}'|}{G_1}\mathbf{e}_1 + 
    \varepsilon\frac{\mathbf{G}'}{|\mathbf{G}'|}\left(\frac{1}{\eta} - 2\frac{G_1'}{G_1}\right).\label{eq:self-similar}
\end{equation}
The equation~\eqref{eq:self-similar} should be equipped with two initial conditions: $\mathbf{G}(0)$ and $\mathbf{G}'(0)$. 
It is not clear which initial conditions we have to impose for the reconnection problem. 

The self-similar solution for the model~\eqref{eq:intro_klein_model} is studied in~\cite{banica2017}. It appears that if
a singularity is introduced at the beginning it will persist for the infinite time. Thus, the self-similar reconnection in 
the model~\eqref{eq:intro_klein_model} will never have a clear horseshoe in difference with the corner vortex studied in~\cite{vega2003}.
One of the reasons can be that in the model~\eqref{eq:intro_klein_model} the LIA term has a linear approximation. In the equation~\eqref{eq:self-similar}
we include this term in the complete nonlinear form.

\paragraph*{Behavior close to the reconnection point}
Using formula~\eqref{eq:arclength} we can bound $T_1 / \|\mathbf{T}\|$ before the reconnection moment from bellow. 
Assume that in the interval $s\in[s_0, s_1]$ the component $x_1(s)$ is growing monotonically (so $T_1(s) \ge 0$) for value
$m = x_1(s_0)$ to $M = x_1(s_1)$. This assumption is correct before reconnection, but the numerical simulation shows that
it does not hold true after it since the helical waves emerge. Subject to the proposed assumption we can write the following
estimation:
\begin{multline}
    \sup_{s\in (s_0, s_1)} \frac{|T_1|}{\|\mathbf{T}\|} \ge \frac{1}{s_1 - s_0} \int_{s_0}^{s_1}\frac{|T_1|}{\|\mathbf{T}\|}ds = 
    \frac{1}{s_1 - s_0} \int_{s_0}^{s_1}\frac{x_1' ds}{ L_0(s) {\left(x_1^2 + r_c^2\right)}^{-\varepsilon / 2}} \\\ge 
    \left.\frac{ r_c^\varepsilon x_1(s) }{(s_1 - s_0)\|L_0\|_{C([s_1,s_0])}} 
    {_2 F_1}\left(\frac{1}{2}, -\frac{\varepsilon}{2}, \frac{3}{2}, -\frac{x^2_1(s)}{r_c^2}\right)\right|_{s_0}^{s_1},\label{eq:general_estimation}
\end{multline}
where $\|L_0\|_{C([s_1,s_0])} = \sup_{s\in  (s_0, s_1)} |L_0(s)|$, and ${_2 F_1}(a,b,c,d)$ is the hypergeometric function, the
modulus of the tangential vector is given by~\eqref{eq:arclength}. Assume now for simplicity that $r_c = 0$, the reconnection happens
at $s_0 = 0$, and designate $s_1 = s$. Then, the estimate~\eqref{eq:general_estimation} reads
\begin{equation}
    \sup_{q\in (0, s)} \frac{|T_1|}{\|\mathbf{T}\|} \ge \frac{x_1^{1+\varepsilon}(s)}{s(1+\varepsilon)\|L_0\|_{C([0,s])}}.\label{eq:simplified_estimation}
\end{equation}
If initially the vortices were oriented into $x_3$ direction and separated by value $2b$ the norm $\|L_0\|_{C([0,s])} = b^\varepsilon$.
The furthest point between vortices corresponds to $s = \lambda / 2$, where $\lambda$ is the wavelength of Crow waves given by~\eqref{eq:crow_waves}. 
The value of $x_1(s)$ in this point is not smaller than $b$, so we can use it to estimate the ratio $\frac{|T_1|}{\|\mathbf{T}\|}$:
\begin{equation}
    \sup_{q\in (0, s)} \frac{|T_1|}{\|\mathbf{T}\|} \ge \frac{\sqrt{\varepsilon}}{(1+\varepsilon)\pi}.
    \label{eq:final_estimate}
\end{equation}
This bound is not optimal and is very far from it. Nevertheless, we can expect that the first component of the tangent vector will grow when we increase
$\varepsilon$ tending it to $1$. The numerical experiments in section~\ref{sec:results} show that it becomes almost parallel to $\mathbf{e_1}$, so the shape 
of the vortex after reconnection is very close to a horseshoe.

We can also use the formula~\eqref{eq:simplified_estimation} to estimate the maximal possible value of $x_1$. Indeed, the left hand side can not be bigger than $1$,
so assuming that $\|L_0\|_{C([0,s])} = b^\varepsilon$ and the maximum for $x_1(s)$ is achieved at $s = \lambda / 2$ we obtain;
\begin{equation}
    x_1(s) \le b {\left(\frac{(1 + \varepsilon) \pi}{\sqrt{\varepsilon}}\right)}^{\frac{1}{1 + \varepsilon}}.\nonumber
\end{equation}

\section{Numerical method}\label{sec:numerics}

In this section we describe the numerical method we use to solve the system~\eqref{eq:main_equation}.
The main problem for numerical stability of the method is related to the interaction term which
grows when the vortices are close to each other. We consider a simpler case of Klein-Majda equations~\eqref{eq:km_reduced} to derive
possible restrictions for the numerical method. They come from the relation between the time step $\tau$,
the spatial discretization step $h$, and the regularization parameter $r_c$. Even though we use a more advanced
Runge-Kutta-Felhberg method and the equations~\eqref{eq:main_equation} have higher nonlinearity the derived restrictions hold
true in a qualitative way.

\paragraph*{Necessary stability conditions}
We derive the necessary stability conditions for a simpler case of Klein-Majda equations~\eqref{eq:km_reduced}. 
Consider the following semi-implicit numerical scheme:
\begin{eqnarray}
&&x_n^{(k+1)} = x_n^{(k)} - \frac{\tau}{h^2}\left(y_{n+1}^{(k)} - 2y_{n}^{(k)} + y_{n-1}^{(k)}\right),\label{eq:stab_scheme_x}\\    
&&y_n^{(k+1)} = y_n^{(k)} + 
\frac{\tau}{h^2}\left(x_{n\!+\!1}^{(k\!+\!1)}\!-\!2x_{n}^{(k\!+\!1)}\!+\!x_{n\!-\!1}^{(k\!+\!1)}\right) - 
\varepsilon \tau \frac{x_{n}^{(k\!+\!1)}}{{x_{n}^{(k\!+\!1)}}^2\!+\!r_c^2};\label{eq:stab_scheme_y}
\end{eqnarray}
where $x_n^{(k)}$ and $y_n^{(k)}$ are approximation of first and second components respectively 
of the solution $\mathbf{X}(s_n,t_k)$, $h$ and $\tau$ are discretization steps for the parameter $s$ and time
respectively. Assume now that there is a high-frequency but small amplitude numerical error $\delta^{(k)} \cos(\omega n)$ 
in the second component and let us analyse how it will grow on the next time step. Plugging the perturbed values 
$x_n^{(k)}+ \delta^{(k)}_x \cos(\omega n)$ and $y_n^{(k)} + \delta^{(k)}_y \cos(\omega n)$ into~\eqref{eq:stab_scheme_x},\eqref{eq:stab_scheme_y} 
we get for the following expression for the linear approximation of the error:
\begin{eqnarray}
&&\delta^{(k\!+\!1)}_x \cos(\omega n)\!=
\!\left(\delta_x^{(k)} + \lambda \zeta \delta_y^{(k)}\right)\cos(\omega n),\label{eq:stab_basic_x}\\
&&\delta^{(k\!+\!1)}_y \cos(\omega n)\!=
\!\left(-(\lambda\zeta + \mu)\delta_x^{(k)} + (1 - \lambda^2 \zeta^2 - \lambda\zeta\mu) \delta_y^{(k)}\right)\cos(\omega n),\label{eq:stab_basic_y}
\end{eqnarray}
where 
\begin{equation}
    \lambda = \tau / h^2,\quad \mu = \epsilon \tau \frac{r_c^2 - {x_n^{(k)}}^2}{{\left(r_c^2 + {x_n^{(k)}}^2\right)}^2},\quad\zeta=2(1 - \cos\omega).\label{eq:lambda_mu}
\end{equation}
The necessary stability condition requires the eigenvalues of the error transformation matrix be not bigger than $1$ by modulus.
From~\eqref{eq:stab_basic_x},\eqref{eq:stab_basic_y} we obtain the equation of the eigenvalues $\nu$:
\begin{equation}
    \begin{vmatrix}1 - \nu & \lambda \zeta \\ -\lambda\zeta - \mu & 1 - \lambda^2 \zeta^2 - \lambda\zeta\mu - \nu\end{vmatrix} = 
    \nu^2 - (2 - \lambda^2 \zeta^2 - \lambda\zeta\mu)\nu + 1 = 0.\nonumber
\end{equation}
The product of the eigenvalues is always $1$, therefore the necessary stability condition is satisfied if and only if the roots
are complex, that is:
\begin{equation}
    0 \le \lambda^2 \zeta^2 + \lambda\zeta\mu \le 4,\label{eq:condition}
\end{equation}
implying
\begin{equation}
    \lambda \zeta < \frac{-\mu + \sqrt{\mu^2 + 16}}{2}.\label{eq:condition_resolved}
\end{equation}
This condition can be resolved providing a constraint for $\tau$:
\begin{equation}
\tau \le \frac{2 h^2}{\sqrt{\zeta^2 + \varepsilon h^2 \zeta \frac{r_c^2 - {x_n^{(k)}}^2}{{\left(r_c^2 + {x_n^{(k)}}^2\right)}^2}}},\label{eq:tau_x_condition}
\end{equation}
for any $\omega$ and $x_n^{(k)}$. The condition has the strongest form when $\zeta$ achieves it maximal value ($\zeta = 4$, see~\eqref{eq:lambda_mu}) 
and $x_n^{(k)} = 0$, so the condition for $\tau$ reads: 
\begin{equation}
    \tau \le \frac{h^2}{\sqrt{4 + \varepsilon  \frac{h^2}{r_c^2}}}.\label{eq:tau_condition}
\end{equation}
Note that for smaller $r_c$ we have to use smaller time step. Furthermore, the stability of the scheme is lost when $r_c$ tends to $0$ what
corresponds to the emergence of a singularity in the interaction term. The formula~\eqref{eq:tau_condition} implies that for say a two times smaller
regularization parameter $r_c$ we have to use a two times smaller time step $\tau$. However, this relation does not hold true for a higher order
scheme. Indeed, in that case we will have a higher derivative of the interaction term respect to $x$. That is the terms $h^3/r_c^3$ and further will
be presented in the constraint for $\tau$. The stability can be obtained by choosing $h$ proportional to $r_c$ so all terms in the Taylor expansion
of the interaction term will be bounded. However, this choice leads to a very fast growth of computations making it very hard to perform the simulation
for small $r_c$.

\paragraph*{Description of the numerical scheme}
The main challenge in the numerical solution of the equations~\eqref{eq:main_equation}
is that at the reconnection moment the behavior of the interaction term is close
to singular. This time period should be passed with very small time step which is not needed
when the vortices are far from each other. Therefore, we use an adaptive time step tecnique:
an embedded 5th Runge-Kutta method in time with 8th order finite difference discretization in 
the filament parameter $s$.  The 8th order scheme gives the best results of those we have tried. 
On one hand, the spectral method that has a higher order requires a higher order time scheme. 
On the other hand, a lower order spatial discretization does not provide sufficient accuracy. 
We have also studied the possibility of use of implicit methods,
such as~\cite{buttke} but these methods suffer the same requirement of the small time step 
at the reconnection moment. Besides, we are interested in the multifractal behavior of the 
trajectories of the vortex filament points therefore we need data with very high discretization
in time thus such advantage of implicit method as large time step can not be used. %

In the implementation of the method we follow the book of J. Butcher~\cite{butcher}. The idea
of the embedded Runge-Kutta method consists in realization of two Runge-Kutta schemes on the 
same points one of order $p$ and another one of order $p + 1$. The difference between the outputs
of these methods on each step is used for the error estimation which should have the decay $\tau^5$.
If the error is bigger then we decrease the time step until the accuracy test is not passed. The 
explicit $k$-steps Runge-Kutta scheme for an ODE $x'(t) = f(x, t)$ at step $n$ is given by
\begin{eqnarray}
&&q_i = f\left(t_n + \alpha_i \tau, x_n + \tau \sum_{l = 1}^{i - 1} \beta_{il} q_l\right),\ 1 \le i < k;\\
&&x_{n + 1} = x_n + \tau \sum_{i = 1}^{k - 1} c_i q_i;
\end{eqnarray}
where $\tau$ is the time step, $\alpha_i$, $\beta_{il}$, and $c_i$ are the coefficients of the scheme.
Usually the coefficients are given in the Butcher table:
    \begin{tabular}{c|c}
        $\boldsymbol{\alpha}$ & $\boldsymbol{\beta}$ \\\hline
        & $\mathbf{c}^T$    
    \end{tabular}.
In order to add the accuracy test we have to add another vector of coefficients $\mathbf{\hat{c}}$ for
the embedded method. We are using Runge-Kutta-Felhberg method with the coefficients obtained in~\cite{fehlberg1969}.
The adaptive time step allows to decrease the time step when it is necessary. In our case when the reconnection happens
the interaction term is very close to singular and therefore, a much smaller time step comparing with the rest of the simulation has to be used.

To make the solution more stable we use the idea of \cite{delahoz2014} and resolve the equations for $\mathbf{X}$ and $\mathbf{T}$
at the same time adding also the arclength correction according to~\eqref{eq:arclength}. The new equations read
\begin{eqnarray}
&&\mathbf{X}_t\!=\!\frac{\mathbf{T}\!\wedge\!\mathbf{T}_{s}}{|\mathbf{T}|^3}\!-\! 
\frac{\varepsilon x_1}{x_1^2 + r_c^2} \frac{\mathbf{T}\wedge \mathbf{e}_1 }{|\mathbf{T}|}.\label{eq:num_system_X}\\
&&\mathbf{T}_t\!=\!\frac{\mathbf{T}\!\wedge\!\mathbf{T}_{ss}}{|\mathbf{T}|^3}\!-\! 
3 \frac{\mathbf{T}\!\wedge\!\mathbf{T}_{s}}{|\mathbf{T}|^4} \frac{\partial |\mathbf{T}|}{\partial s} \nonumber\\
&&\qquad-\frac{\varepsilon}{\left(x_1^2\!+\!r_c^2\right)|\mathbf{T}|} \left(\!x_1 \mathbf{T}_s\!\wedge\!\mathbf{e}_1\!+\! 
\left(\!\frac{r_c^2\!-\!x_1^2}{x_1^2\!+\!r_c^2}\!-\!
\frac{x_1}{|\mathbf{T}|}\frac{\partial |\mathbf{T}|}{\partial s}\!\right) 
\mathbf{T}\!\wedge\!\mathbf{e}_1\!\right).\label{eq:num_system_T}
\end{eqnarray}
Expression~\eqref{eq:arclength} for the modulus of the tangential vector allows us to avoid the calculation of the derivative:
\begin{equation}
    \frac{\partial |\mathbf{T}(s,t)|}{\partial s} = \left(\frac{L'_0(s)}{L_0(s)}  - 
    \varepsilon  \frac{x_1(s,t) T_1(s,t)}{x_1^2(s,t) + r_c^2}\right)|\mathbf{T}(s,t)|.
\end{equation}
Here $L_0(s)$ and its derivative are given as initial condition. Following~\cite{delahoz2014} we also perform the 
correction of the tangential vector modulus after each interaction using formula~\eqref{eq:arclength}. The method
is explained in detail in the appendix~\ref{sec:algorithm}.

Even though there are two connected unknowns $x_1$ and $T_1$, the solution of the system 
\eqref{eq:num_system_X}-\eqref{eq:num_system_T} provides the correct result $T_1 = x_{1,s}$. 
One can wonder if it is possible to reduce the number of unknowns and what is the minimal
number of independent functions which describe the evolution of the vortex reconnection. This is studied in
the appendix~\ref{sec:frenet} using the Frenet frame. It turns out that the reconnection of vortices can
be described in terms of just two functions: $x_1(s,t)$ and its derivative with respect to time $x_{1,t}(s,t)$. However, the equations in
this case contain derivatives respect to $s$ up to the 4th order. Thus, it is much more complicated for the numerical 
solution.


\section{Results}\label{sec:results}

We use the method described in the section~\ref{sec:numerics} to solve the equations~\eqref{eq:main_equation}.
There are a few things we are mainly interested in: (i) the emergence of Crow waves and their length; (ii) the 
influence of $r_c$ to the solution; (iii) the formation of the horseshoe structure and the direction of the 
tangential vector in that region.

\begin{figure}
    \centering
    \begin{subfigure}[c]{0.32\textwidth}
        \centering
        \includegraphics[width=\textwidth, trim={16cm 4cm 12cm 2cm}, clip]{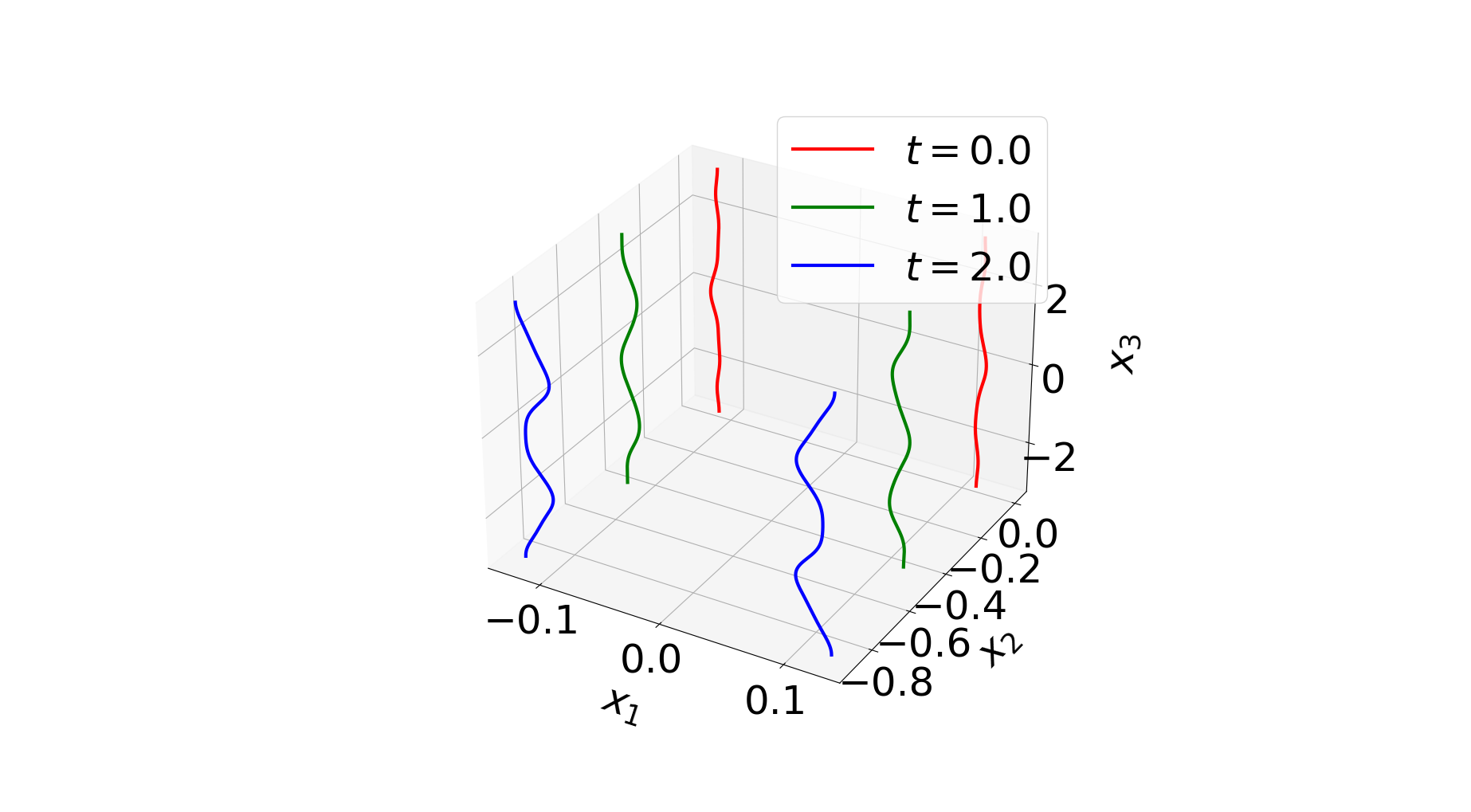}
        \caption{3D view}
    \end{subfigure}
    \begin{subfigure}[c]{0.32\textwidth}
        \centering
        \includegraphics[width=\textwidth, trim={16cm 4cm 12cm 2cm}, clip]{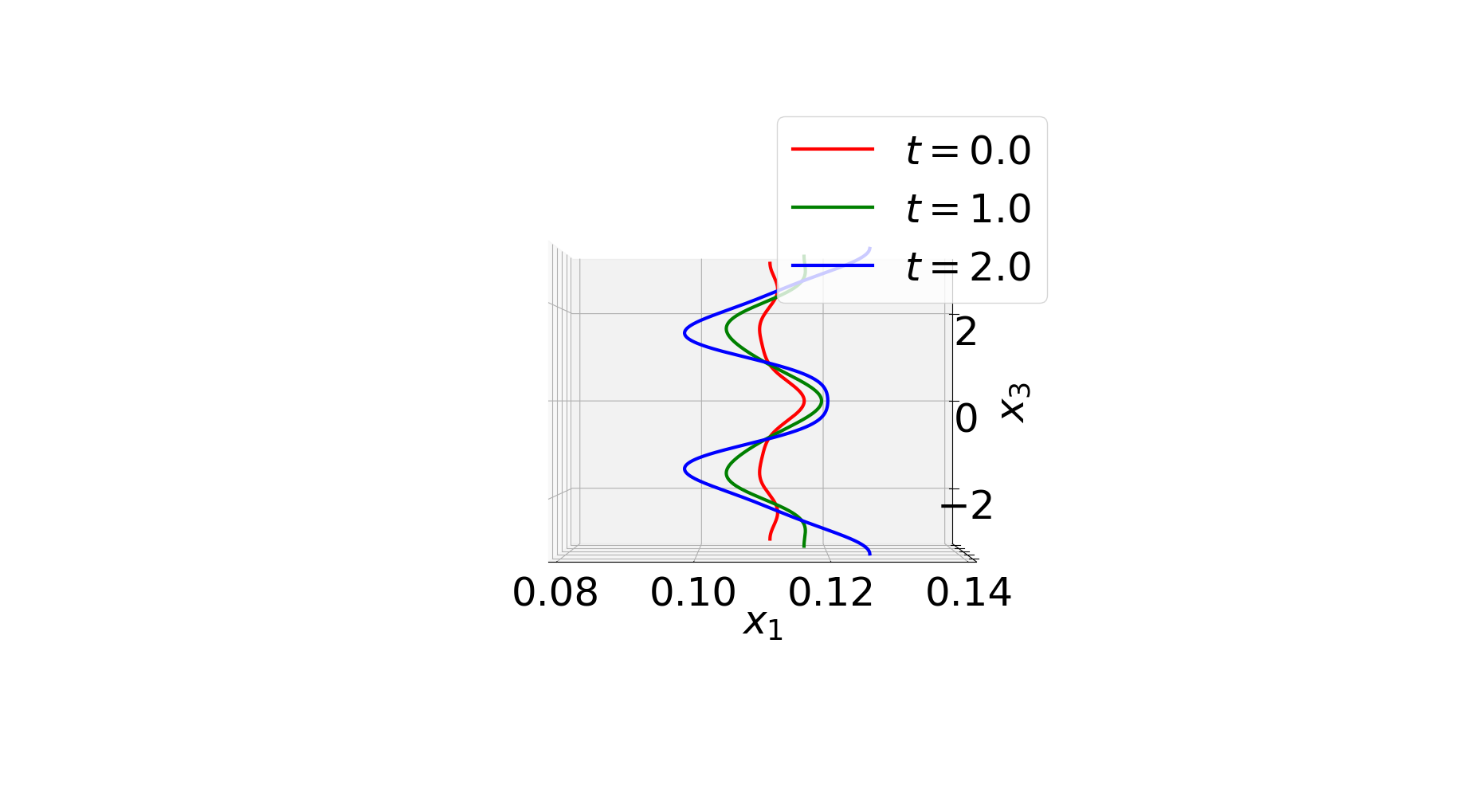}
        \caption{front view}
    \end{subfigure}
    \begin{subfigure}[c]{0.32\textwidth}
        \centering
        \includegraphics[width=\textwidth, trim={16cm 4cm 12cm 2cm}, clip]{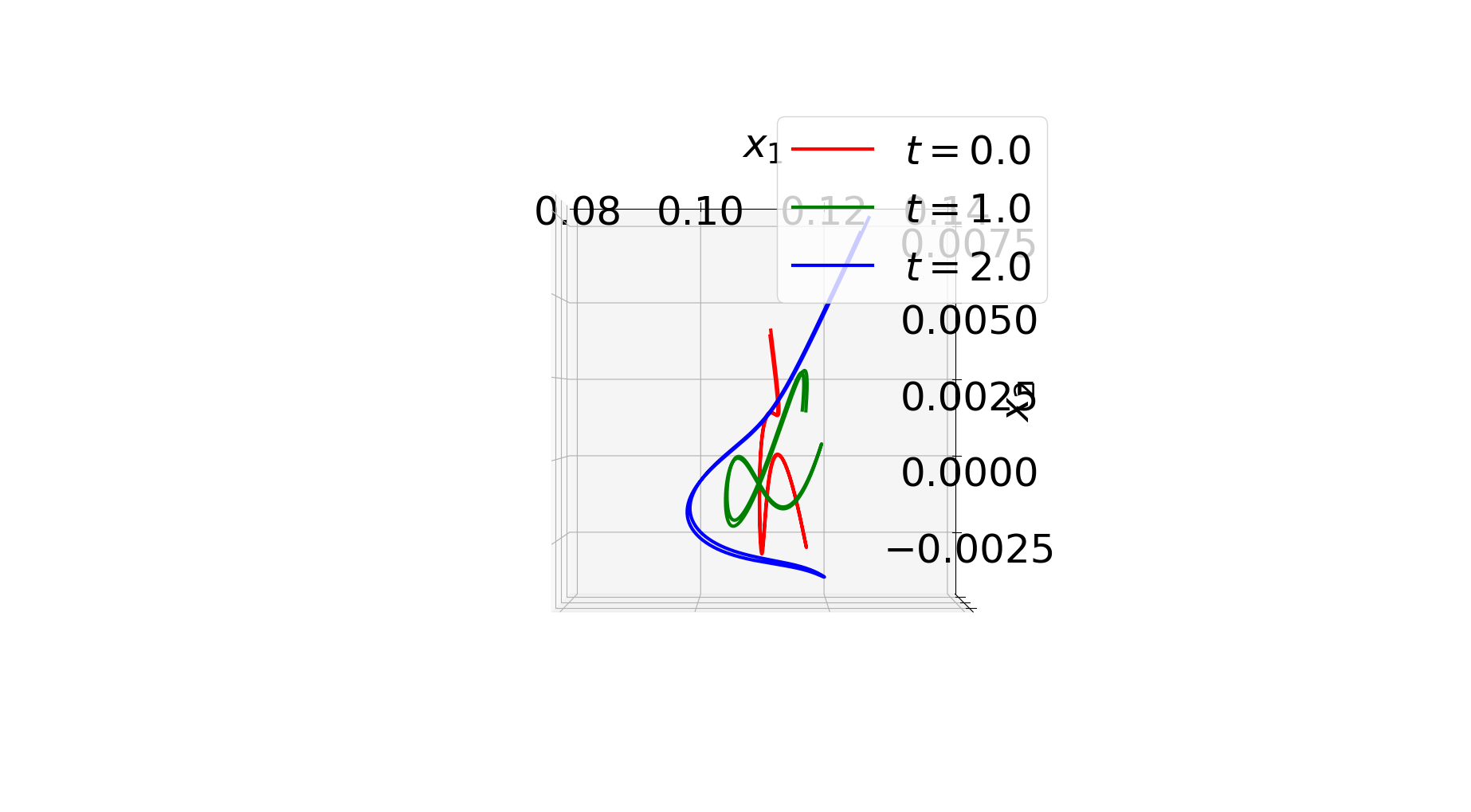}
        \caption{top view}
    \end{subfigure}    
    \caption{Formation of Crow waves for $\varepsilon = 0.05$ and $r_c = 0.0025$. $500$ nodes are used in the discretization. 
    In the front view picture we show only right vortex, in the top view the vortices are shifted by the mean value in $x_2$ direction.}             
    \label{fig:crow_waves}
\end{figure}

\paragraph*{Crow waves} In the first test we start from a random perturbation and check the formation of Crow waves. We start from a small
perturbation of a straight vortex separated by $b = 0.11$ and consider the evolution following the equations~\eqref{eq:main_equation}
selecting $\varepsilon = 0.05$ and $r_c = 0.025$. The results are depicted in figure~\ref{fig:crow_waves}. 
We can see that at time $t = 2$ we have almost sinusoidal waves. The wavelength is around $\pi$ that is very close to
one predicted by formula~\eqref{eq:crow_waves} for the given values of $\varepsilon$, $b$, and $r_c$.

\begin{figure}
    \centering
    \begin{subfigure}[c]{0.32\textwidth}
        \centering
        \includegraphics[width=\textwidth, trim={15cm 3cm 10cm 4cm}, clip]{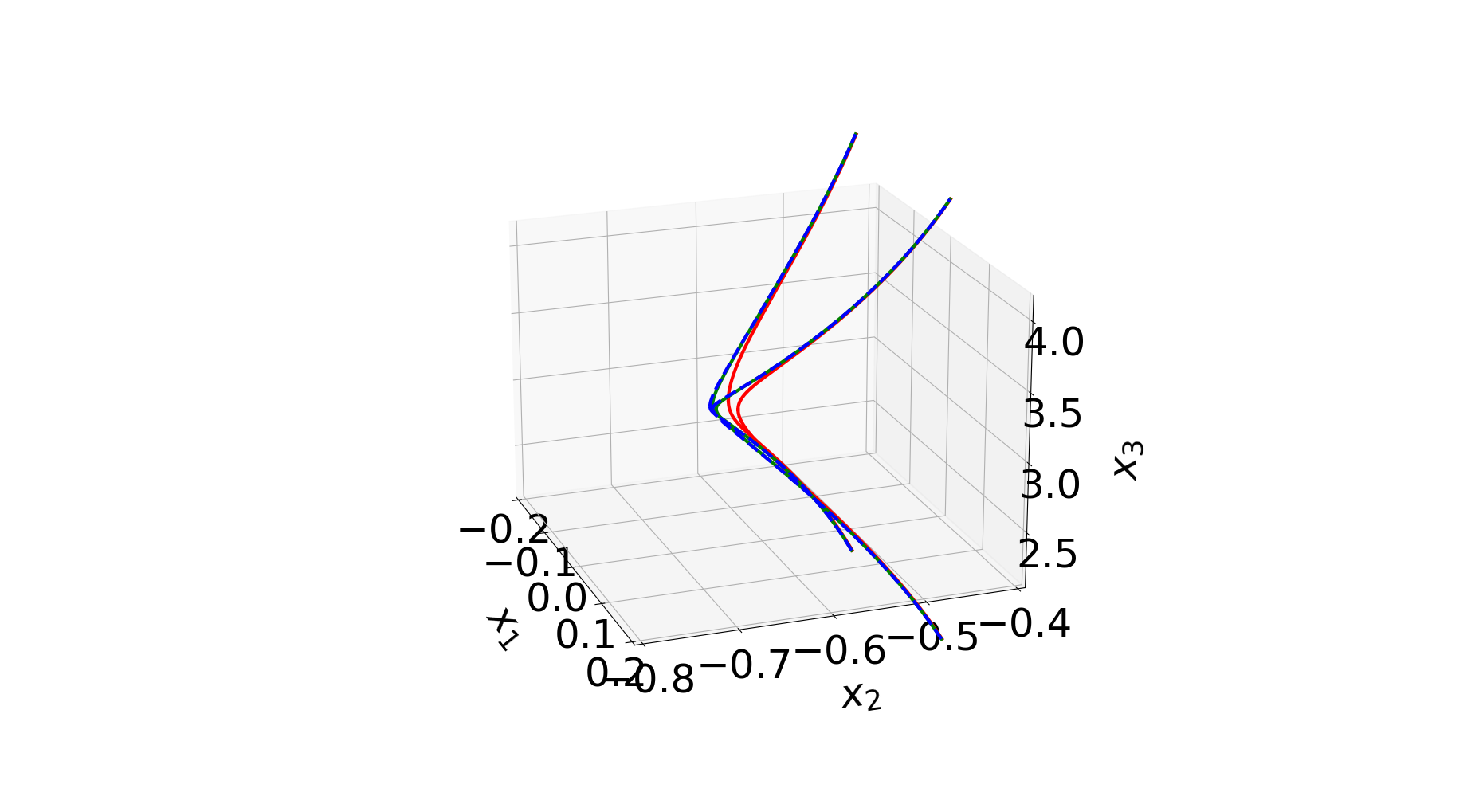}
        \caption{3D view}
    \end{subfigure}      
    \begin{subfigure}[c]{0.32\textwidth}
        \centering    
        \includegraphics[width=\textwidth, trim={15cm 3cm 10cm 4cmm}, clip]{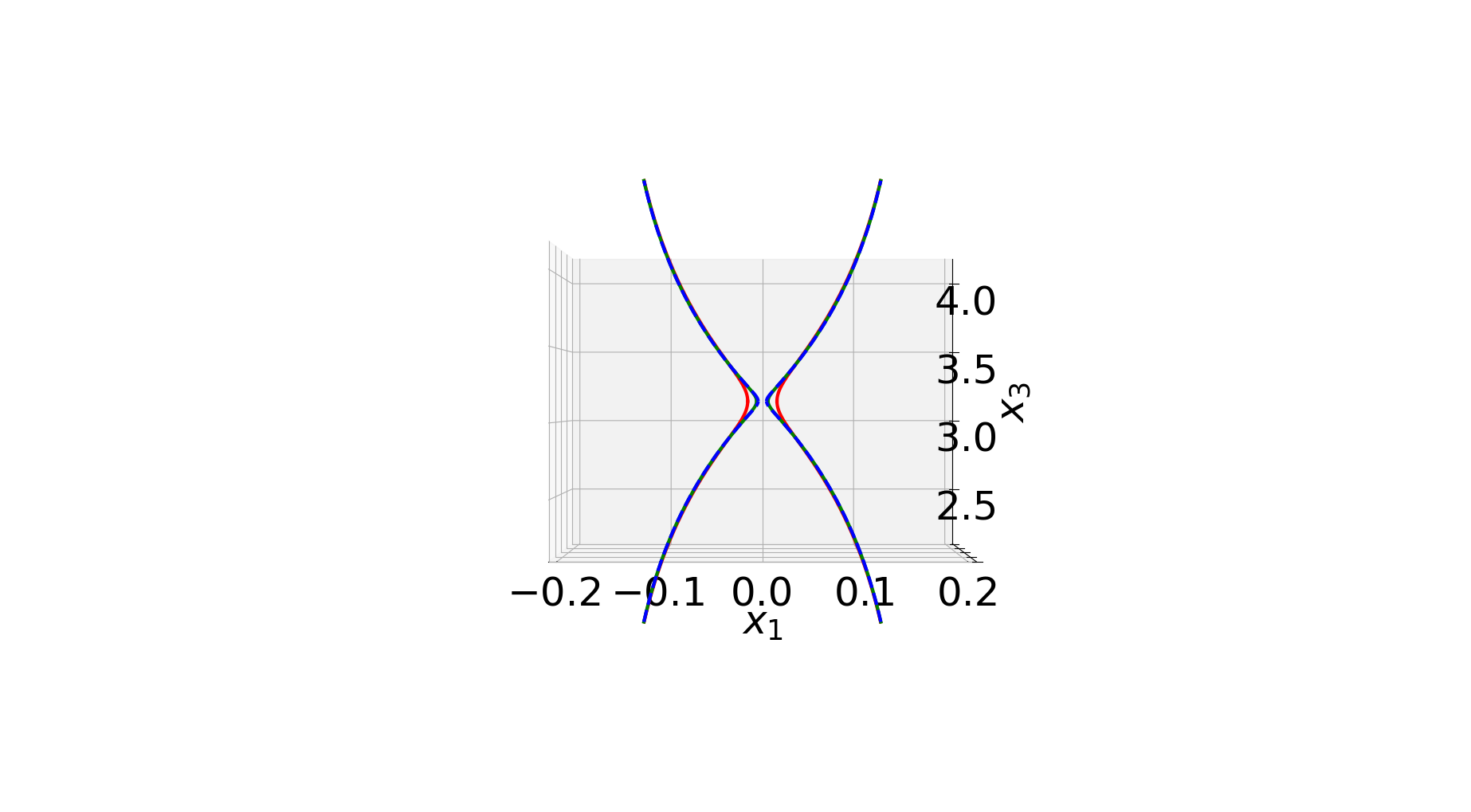}
        \caption{front view}
    \end{subfigure}          
    \begin{subfigure}[c]{0.32\textwidth}
        \centering    
        \includegraphics[width=\textwidth, trim={15cm 3cm 10cm 4cm}, clip]{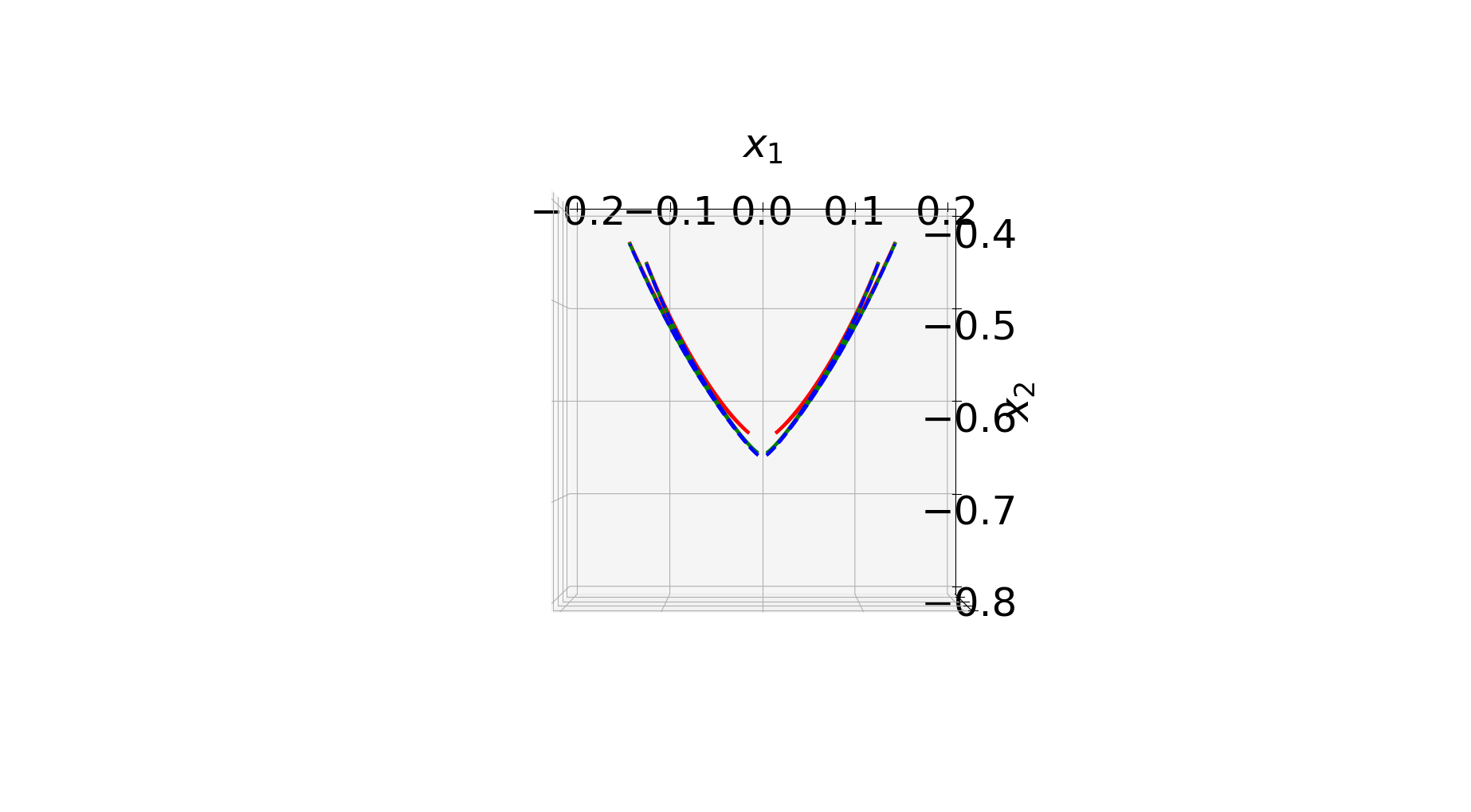}
        \caption{top view}
    \end{subfigure}          
    \caption{Configuration of the vortices before the reconnection at time $t = 1.01$ for different $r_c$ ($5 \cdot 10^{-3}$ -- red, 
    $1.25 \cdot 10^{-3}$ -- green, $3.125 \cdot 10^{-4}$ -- blue) and $\varepsilon = 0.05$}
    \label{fig:rc_before_rec}
\end{figure}

\begin{figure}
    \centering
    \begin{subfigure}[c]{0.32\textwidth}
        \centering    
        \includegraphics[width=\textwidth, trim={15cm 3cm 10cm 4cm}, clip]{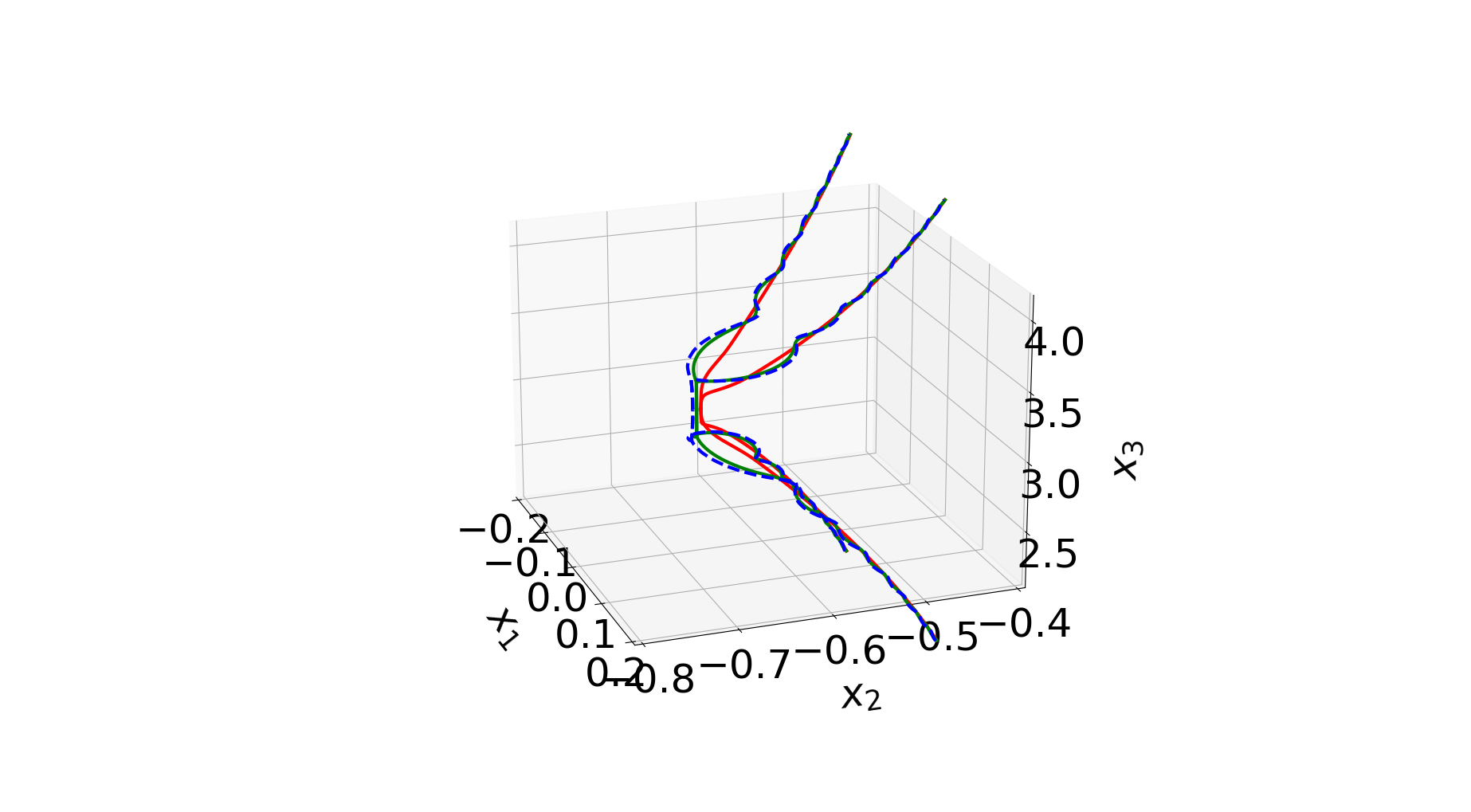}
        \caption{3D view}
    \end{subfigure}          
    \begin{subfigure}[c]{0.32\textwidth}
        \centering    
        \includegraphics[width=\textwidth, trim={15cm 3cm 10cm 4cm}, clip]{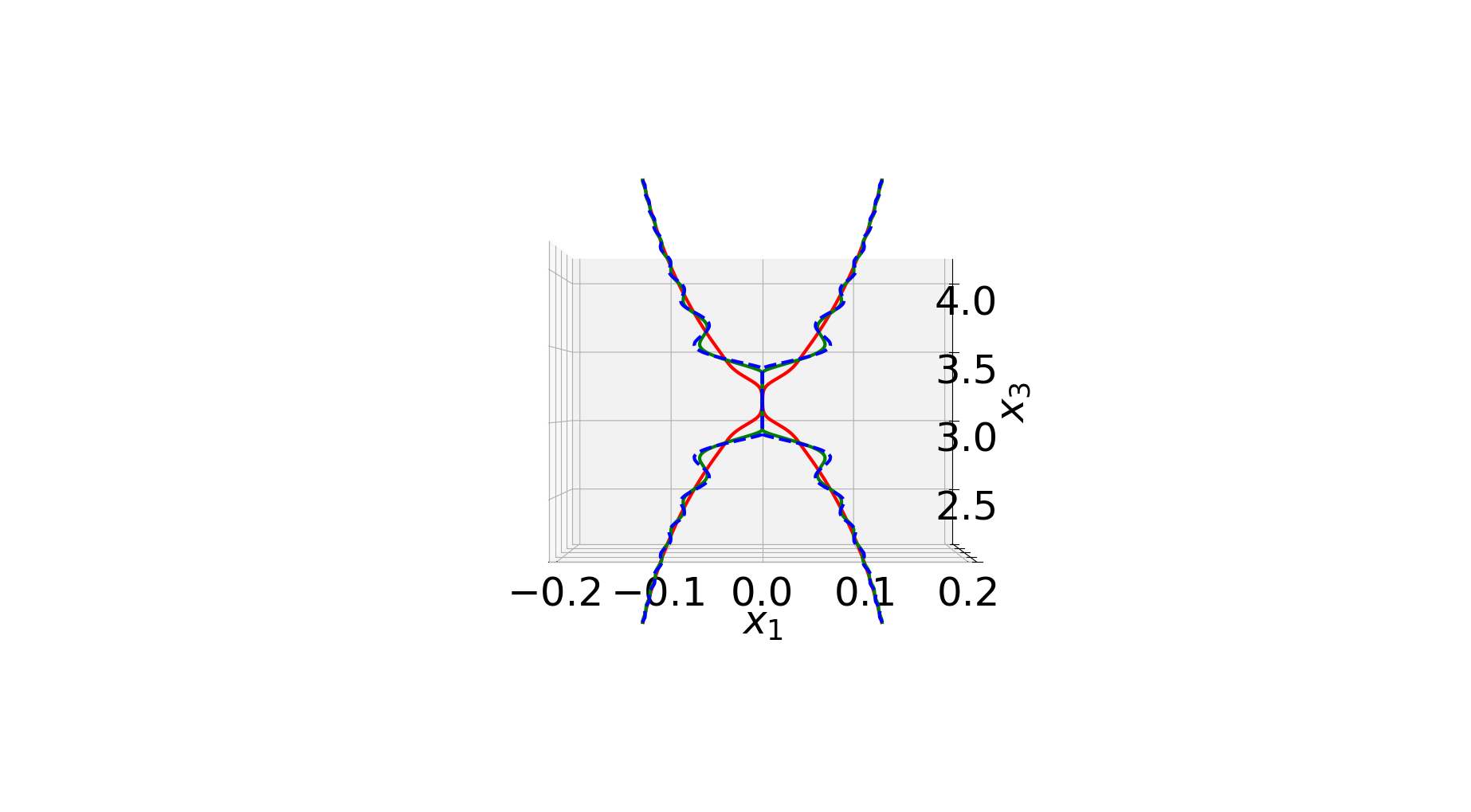}
        \caption{front view}
    \end{subfigure}      
    \begin{subfigure}[c]{0.32\textwidth}
        \centering    
        \includegraphics[width=\textwidth, trim={15cm 3cm 10cm 4cm}, clip]{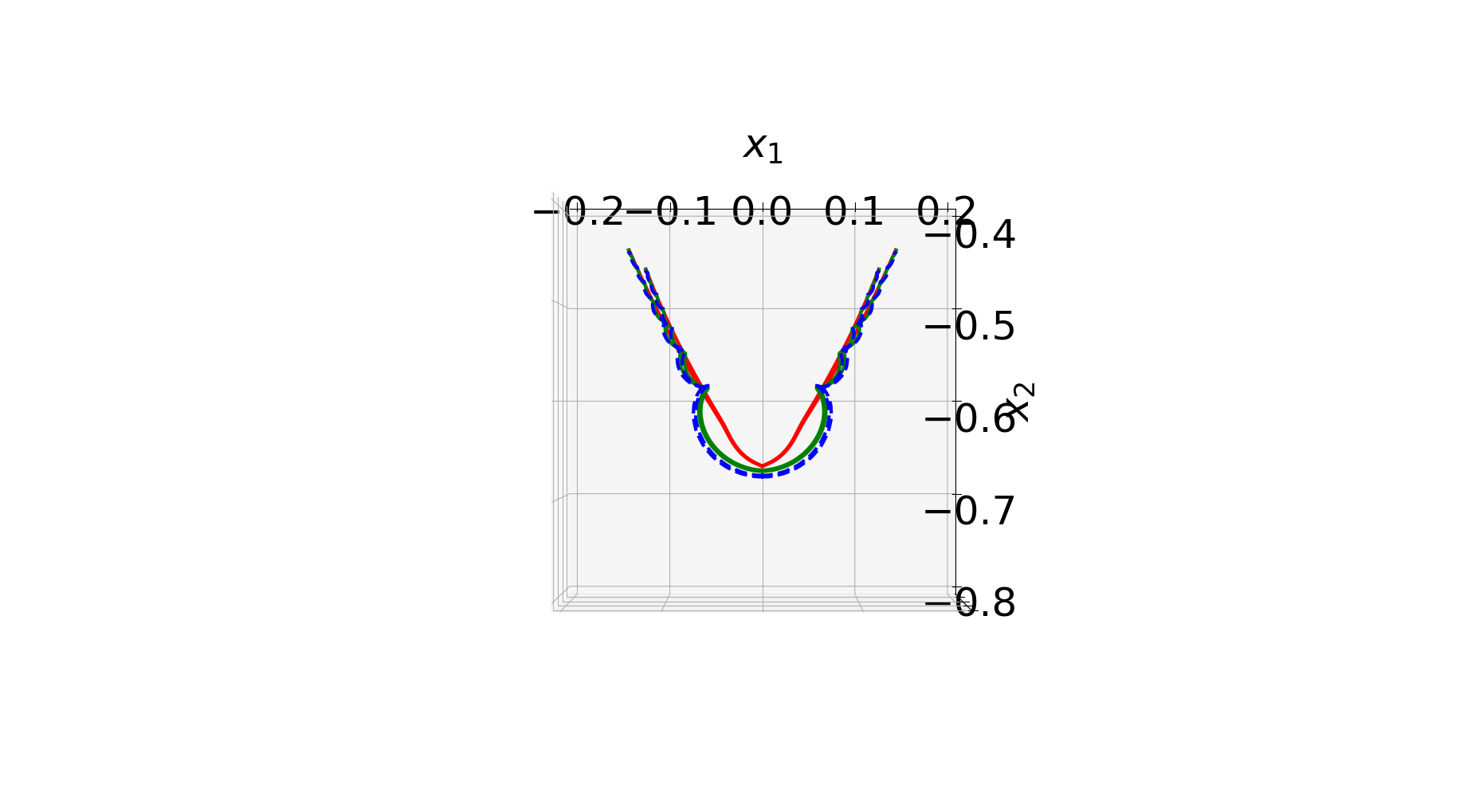}
        \caption{top view}
    \end{subfigure}      
    \caption{Configuration of the vortices at time $t = 1.25$ when the horseshoe structure emerged for all $r_c$ ($5 \cdot 10^{-3}$ -- red, 
    $1.25 \cdot 10^{-3}$ -- green, $3.125 \cdot 10^{-4}$ -- blue) and $\varepsilon = 0.05$}
    \label{fig:rc_at_rec}
\end{figure}

\begin{figure}
    \centering
    \begin{subfigure}[c]{0.32\textwidth}
        \centering  
        \includegraphics[width=\textwidth, trim={15cm 3cm 10cm 4cm}, clip]{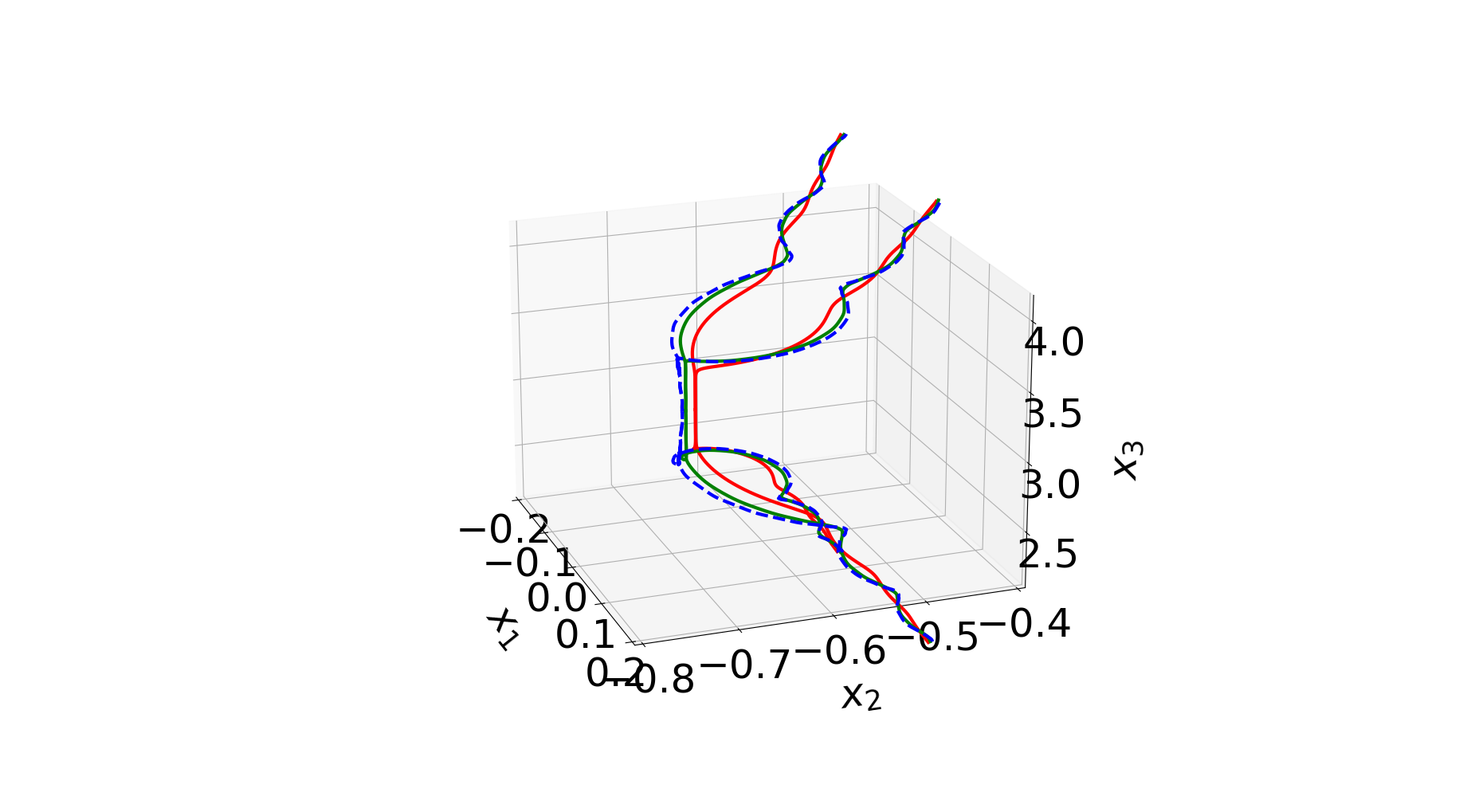}
        \caption{3D view}
    \end{subfigure}  
    \begin{subfigure}[c]{0.32\textwidth}
        \centering      
        \includegraphics[width=\textwidth, trim={15cm 3cm 10cm 4cm}, clip]{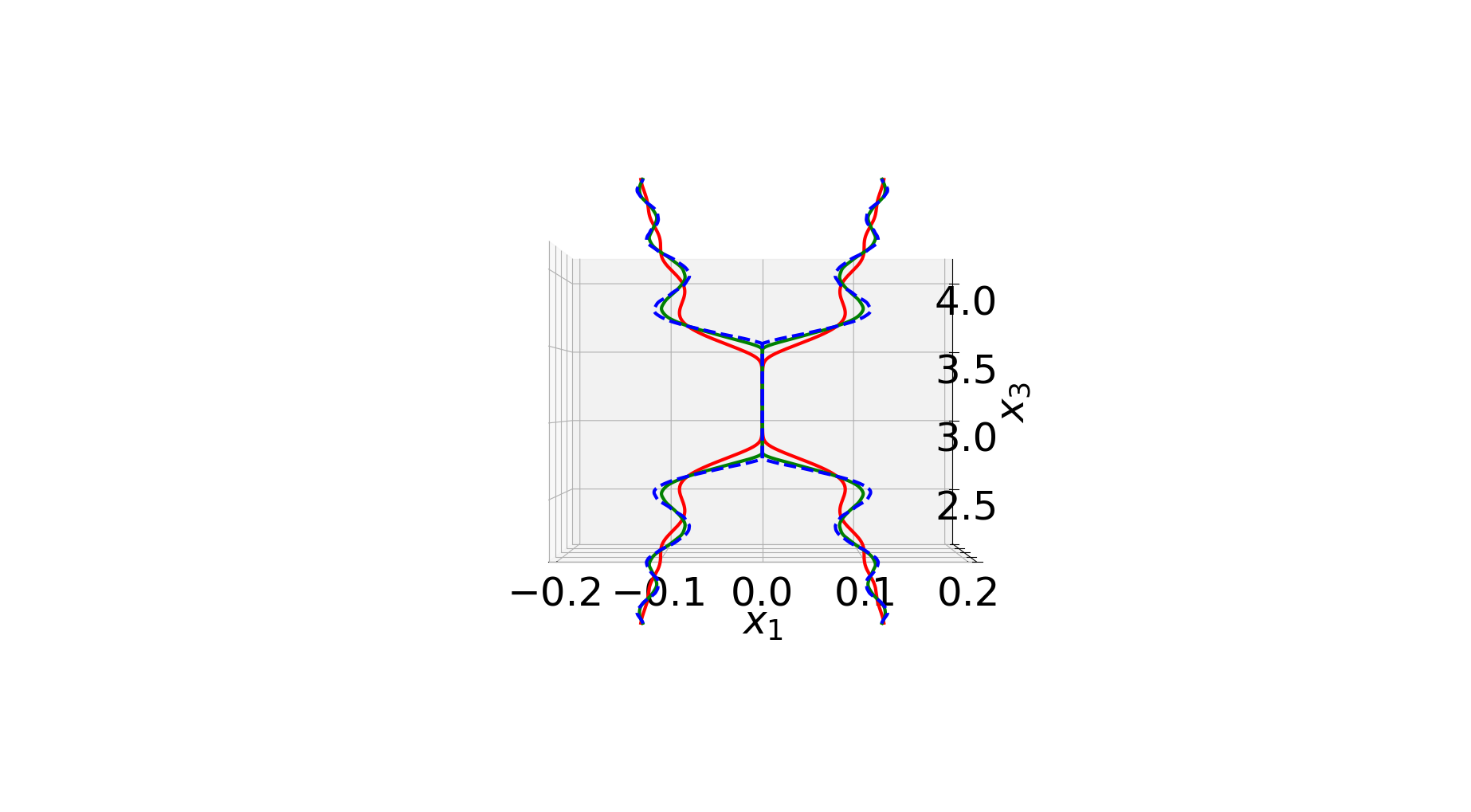}
        \caption{front view}
    \end{subfigure}  
    \begin{subfigure}[c]{0.32\textwidth}
        \centering      
        \includegraphics[width=\textwidth, trim={15cm 3cm 10cm 4cm}, clip]{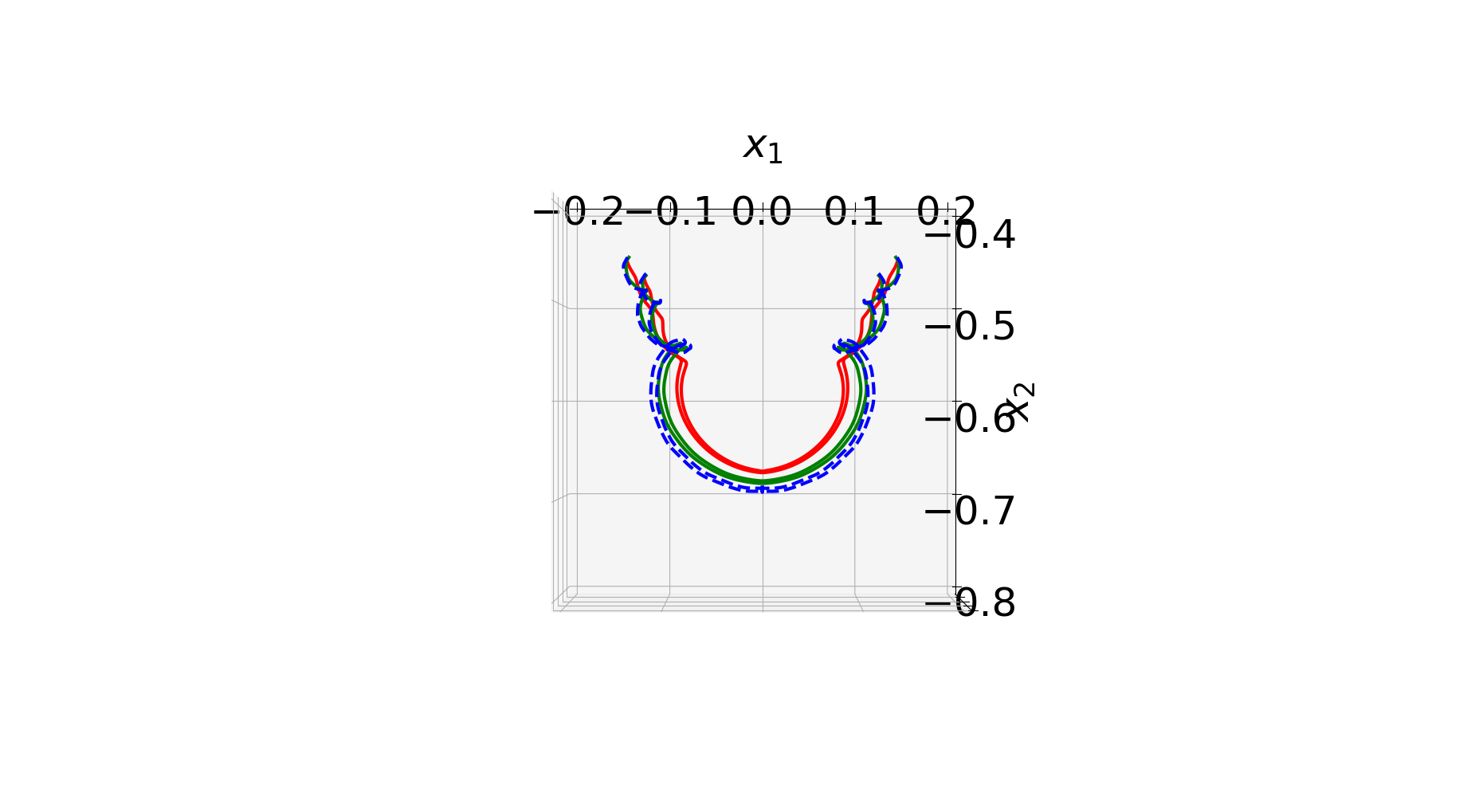}
        \caption{top view}
    \end{subfigure}  
    \caption{Configuration of the vortices after the reconnection at time $t = 1.5$ for different $r_c$ ($5 \cdot 10^{-3}$ -- red, 
    $1.25 \cdot 10^{-3}$ -- green, $3.125 \cdot 10^{-4}$ -- blue) and $\varepsilon = 0.05$}
    \label{fig:rc_after_rec}
\end{figure}

\paragraph*{Influence of the regularization parameter $r_c$} 
The behavior of the vortices far from the reconnection is not dependent of $r_c$. However, at times when the vortices are close to 
each other the regularization parameter starts to play a crucial role. Furthermore, the evolution of vortices after 
reconnection also depends on the parameter $r_c$. We may expect that when the regularization parameter tends to 
zero the shape of vortices will be less smooth. Therefore, the behavior of vortices after
reconnection will resemble the behavior of a corner vortex. These expectations are confirmed by the numerical solution,
see figures~\ref{fig:rc_before_rec},\ref{fig:rc_at_rec},\ref{fig:rc_after_rec}. 

We have consider a pair of symmetric vortices with initial conditions 
\begin{equation}
    \mathbf{X}(s,0) = \begin{pmatrix} b - \delta\cos(s) \\ -\delta\cos(s) \\ s - \pi \end{pmatrix}\label{eq:initial_condition}
\end{equation}
where $b = \sqrt{\varepsilon} / 2$, $\delta = b / 20$, in the interval $(0, 2\pi)$ discretized with $6000$ nodes. 
The boundary conditions are periodic.
The parameter $b$ is selected using formula~\eqref{eq:crow_waves} in such a way that there is exactly one Crow wave in the considered interval. 
If we decrease $b$ the reconnection may happen in multiple points thus complicating the analysis, whereas for larger values of $b$ the reconnection
does not happen due to the periodic boundary conditions. The computations until time $t = 1.5$ take around a day on a personal computer that is 
comparable with the performance of the solution of the Navier-Stokes equations from the section~\ref{sec:introduction} on the same computer. The advantage is that now we
can consider much thinner vortices and has $6000$ nodes along the vortex instead of $320$.

The vortices start to touch each other at time $t = 1.01$, figure~\ref{fig:rc_before_rec}. 
The influence of $r_c$ can be notice only close to the reconnection region and the shape of vortices is sharper for small $r_c$. 
It is not completely clear what we can call "the reconnection moment": the first touch or the moment when the horseshoe emerges. 
Both these moments are dependent on $r_c$. However, the second one has a more complicated dependence since the smallest size of the horseshoe is
dictated by $r_c$: the smaller $r_c$, the smaller the horseshoe structure will be. This effect is demonstrated in
figure~\ref{fig:rc_at_rec} where the configuration of vortices at time $t = 1.025$ is depicted. We can see that even though
for all $r_c$ we have a contact the horseshoe appears only for small values of $r_c$. Furthermore, only for small values
of the regularization parameter we can see the helical waves at this time moment. The fact that we can not define the reconnection
moment does not allow us to change the topology of vortices that leads to the formation of the bridge between the horseshoes.
The bridge is growing and represents a source of numerical difficulties and possible instabilities at later times. The bigger 
bridge at time $t = 1.05$ can be seen in figure~\ref{fig:rc_after_rec}. There we can also see the horseshoe and the helical waves
for large values of $r_c$. It is interesting to note that the difference between solutions for $r_c = 1.25\cdot 10^{-3}$ and 
$r_c = 3.125\cdot 10^{-4}$ is almost neglectable. We can expect that there is a convergence when $r_c$ tends to zero.

The analysis of the influence of $r_c$ to the solution shows an important phenomenon: the reconnection moment can not be determined 
so the "real corner" between vortices can not be seen. We can see only the consequences, such as the horseshoe starting from some minimal 
size related to $r_c$. This effect also appears in the experiments~\cite{fonda2019}. Since the regularization parameter have the physical meaning 
of the radius of vortex core we can expect that for
thick vortices we will never see the singularity and the shape after the reconnection will be closer to a vortex ring without any helical waves,
since the minimal size of the horseshoe is big and close to the length of Crow waves. When the vortices get thinner the shape after
reconnection becomes more complicated. However, we still do not see the singularity. One of the challenges related to this phenomenon is the
incapacity to perform the reconnection, that is to say the change topology, because we do not know when we have to do it. On the one hand if we 
reconnect the vortices when they touch each other, the corner they create is not the one that generates the horseshoe structure. On the other hand, if we wait
until the horseshoe structure emerges we always find an artifact on its tip related to the bridge.

\paragraph*{Influece of the parameter $\varepsilon$}

The results obtained in section~\ref{sec:analysis} predicts that the horseshoe will be closer to circular when the value of
$\varepsilon$ increases (formula~\eqref{eq:final_estimate}). We use the same initial condition~\eqref{eq:initial_condition} 
as before but now the initial distance $b = \sqrt{\varepsilon} / 2$, and the perturbation $\delta = b / 20$ both depend on the parameter $\varepsilon$.
It is necessary for two things: firstly, we have to use such initial distance, so there is at least one Crow wave in the interval, secondly
the change of perturbation amplitude allows to achieve the reconnection almost at the same time. The regularization parameter is set to 
$r_c = 0.05$ and the interval $(0, 2\pi)$ is discretized with $6000$ nodes.

\begin{figure}
    \centering
    \begin{subfigure}[c]{0.32\textwidth}
        \centering      
        \includegraphics[width=\textwidth, trim={15cm 3cm 10cm 4cm}, clip]{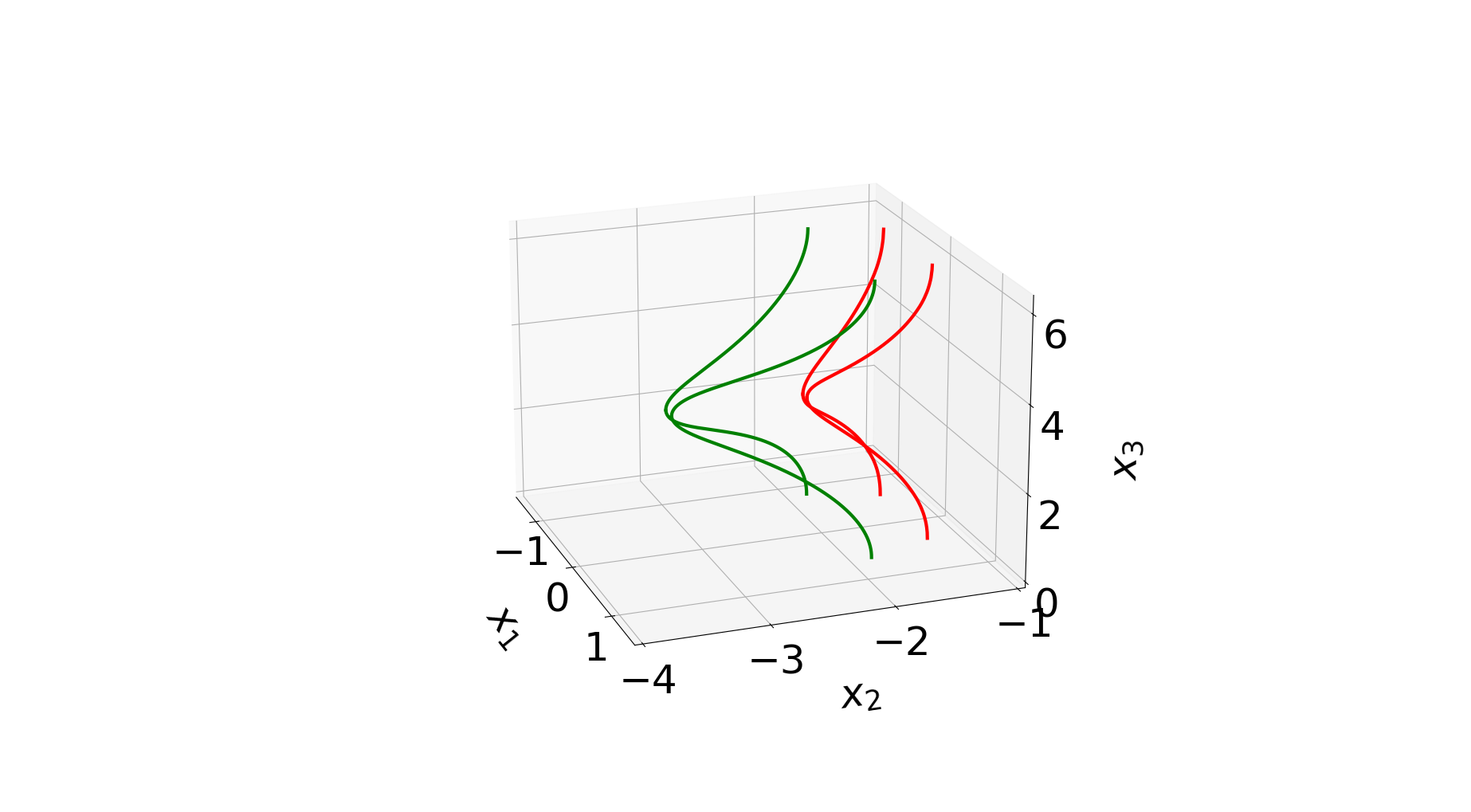}
        \caption{3D view}
    \end{subfigure}
    \begin{subfigure}[c]{0.32\textwidth}
        \centering   
        \includegraphics[width=\textwidth, trim={15cm 3cm 10cm 4cm}, clip]{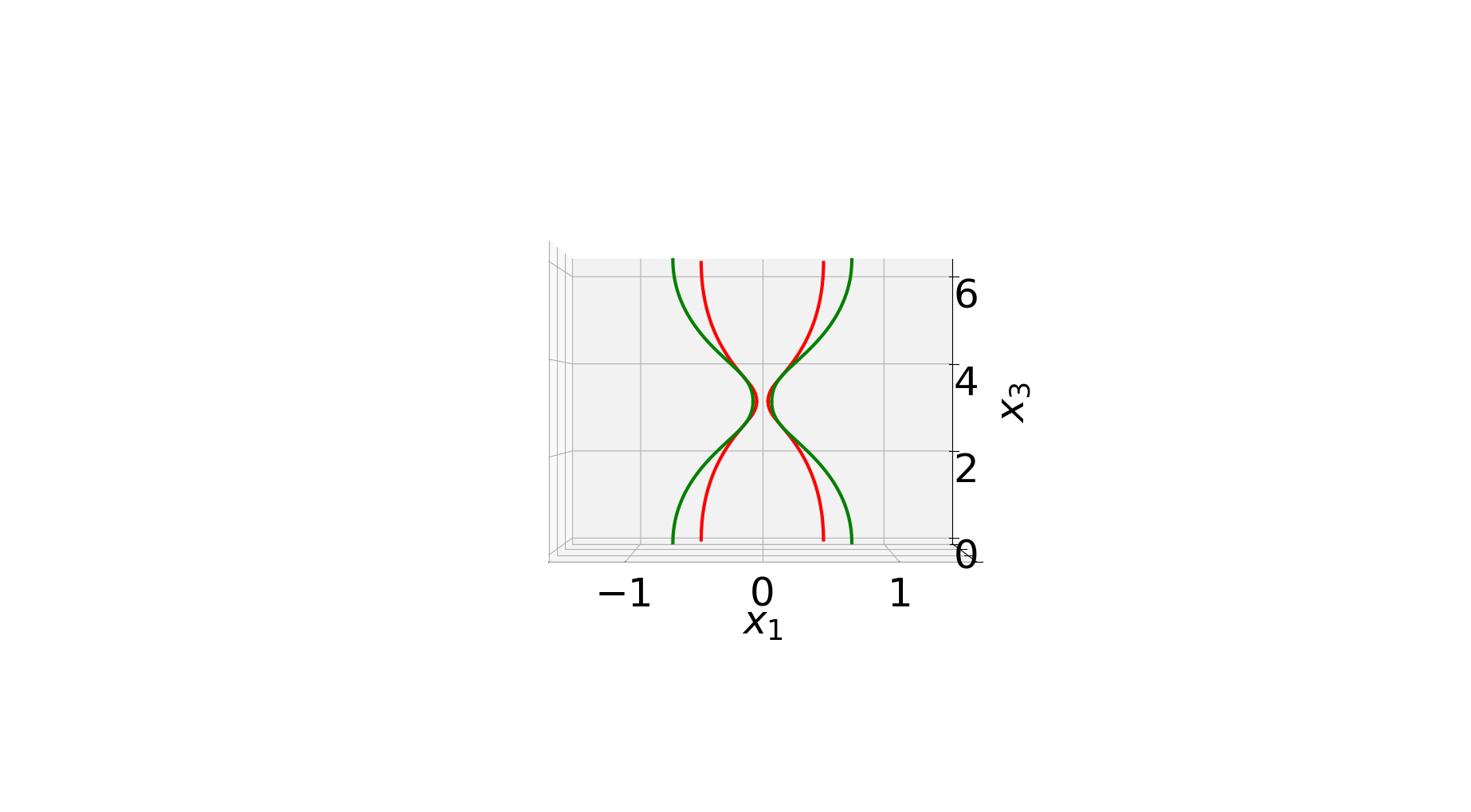}
        \caption{front view}
    \end{subfigure}
    \begin{subfigure}[c]{0.32\textwidth}
        \centering     
        \includegraphics[width=\textwidth, trim={15cm 3cm 10cm 4cm}, clip]{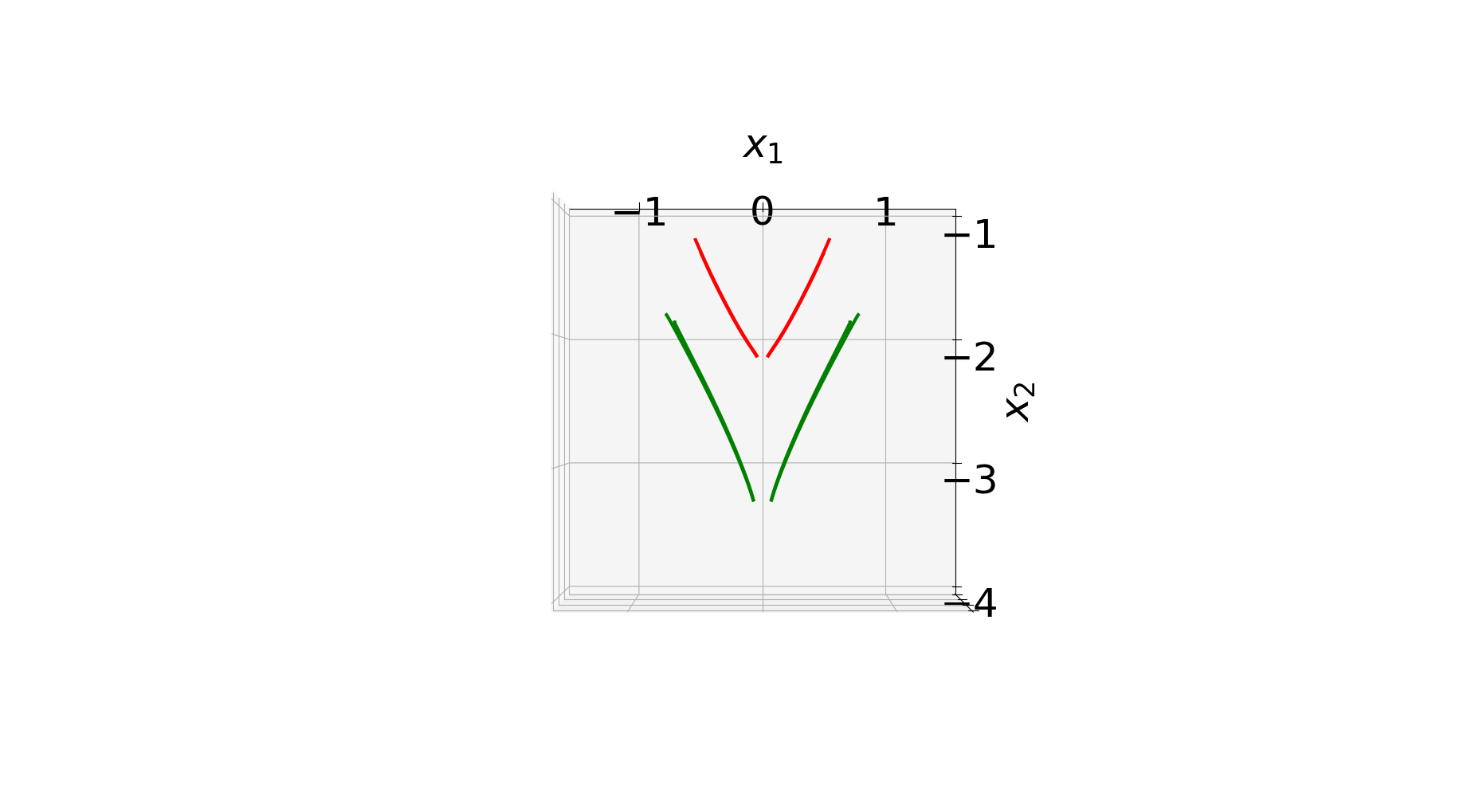}
        \caption{top view}
    \end{subfigure}
    \caption{Configuration of the vortices before the reconnection at time $t = 1.1$ for different $\varepsilon$ 
    ($\varepsilon = 0.5$ -- red, $\varepsilon = 1$ -- green) and $r_c~=~0.05$}
    \label{fig:eps_before_rec}
\end{figure}

\begin{figure}
    \centering
    \begin{subfigure}[c]{0.32\textwidth}
        \centering  
        \includegraphics[width=\textwidth, trim={15cm 3cm 10cm 4cm}, clip]{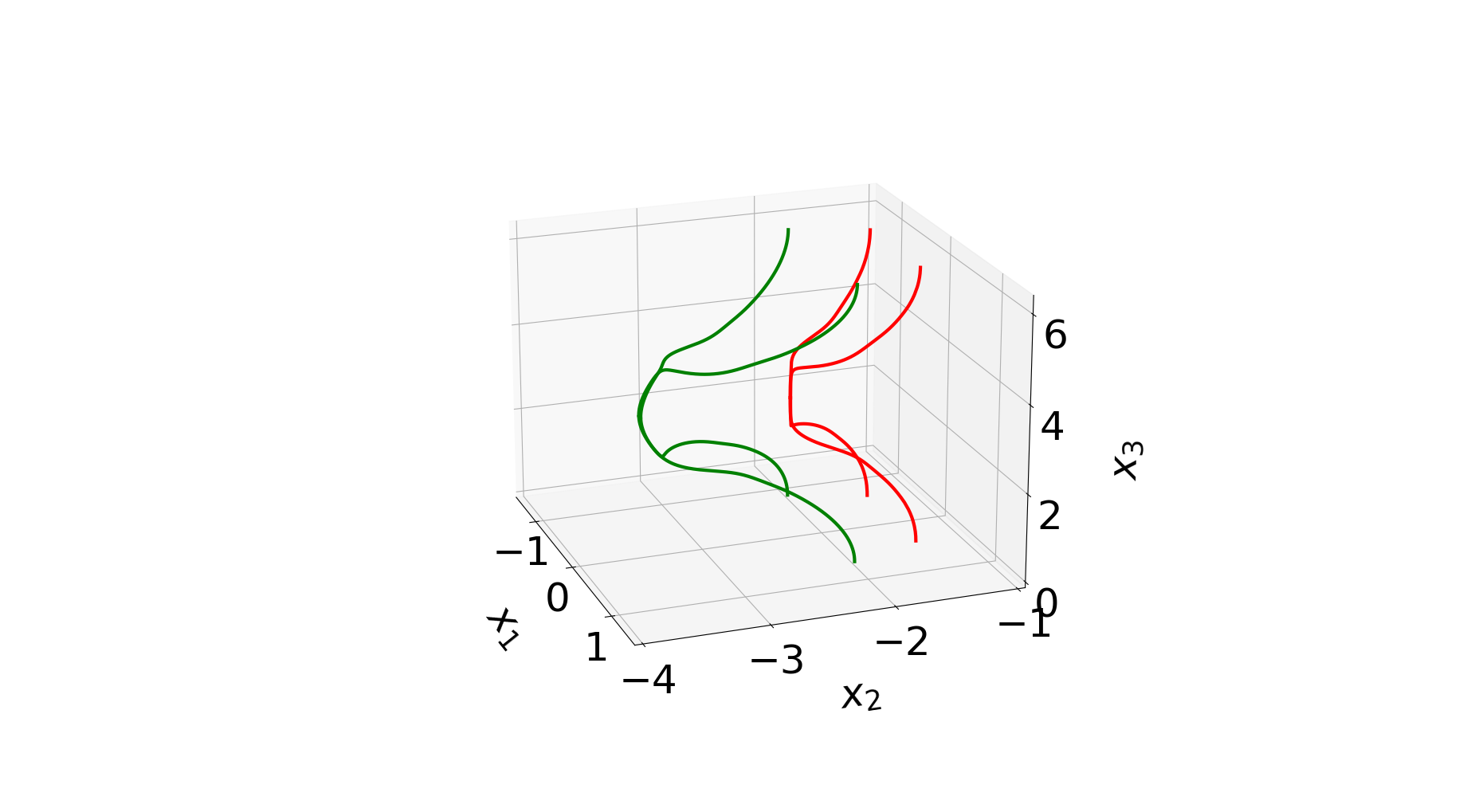}
        \caption{3D view}
    \end{subfigure}
    \begin{subfigure}[c]{0.32\textwidth}
        \centering  
        \includegraphics[width=\textwidth, trim={15cm 3cm 10cm 4cm}, clip]{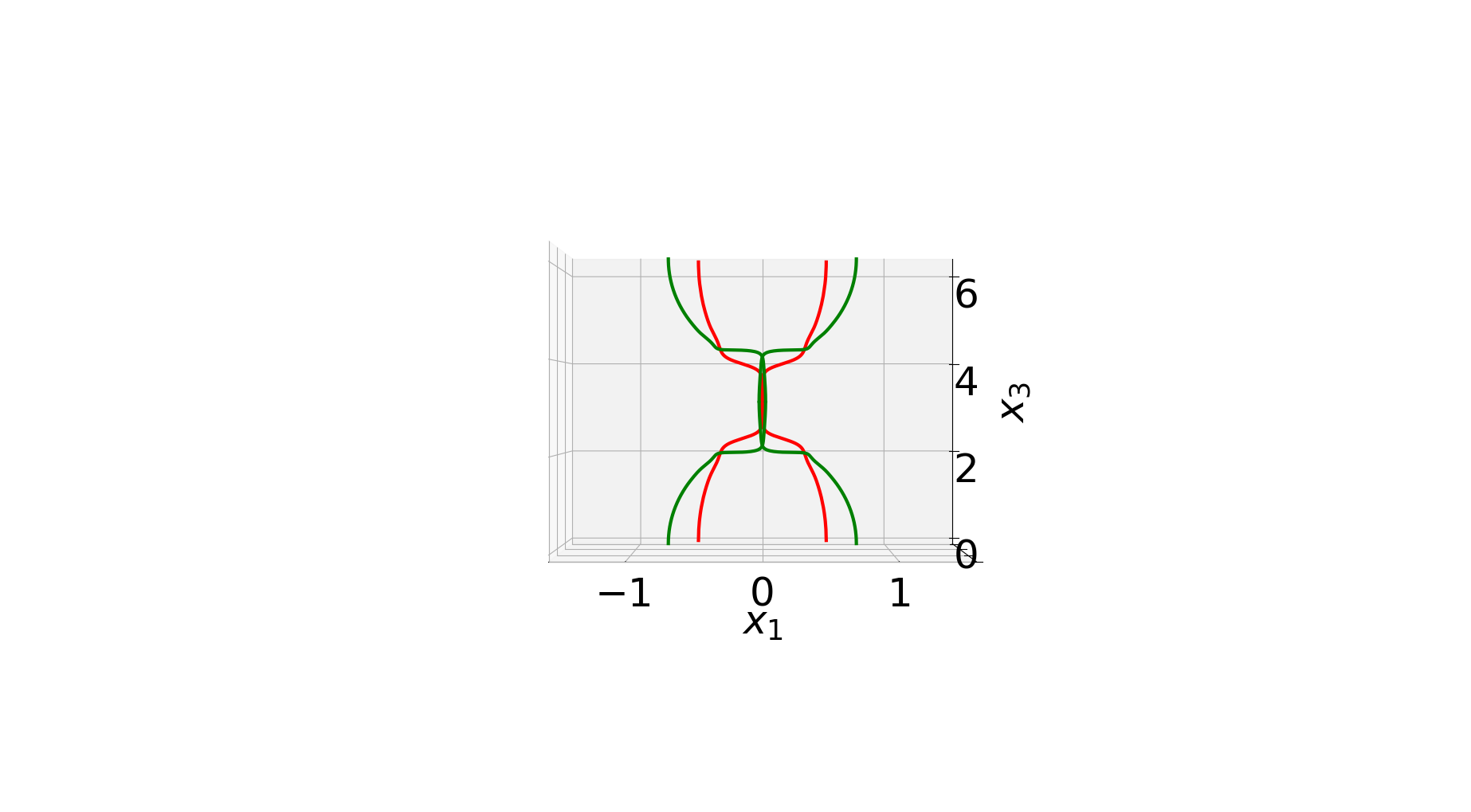}
        \caption{front view}
    \end{subfigure}
    \begin{subfigure}[c]{0.32\textwidth}
        \centering  
        \includegraphics[width=\textwidth, trim={15cm 3cm 10cm 4cmm}, clip]{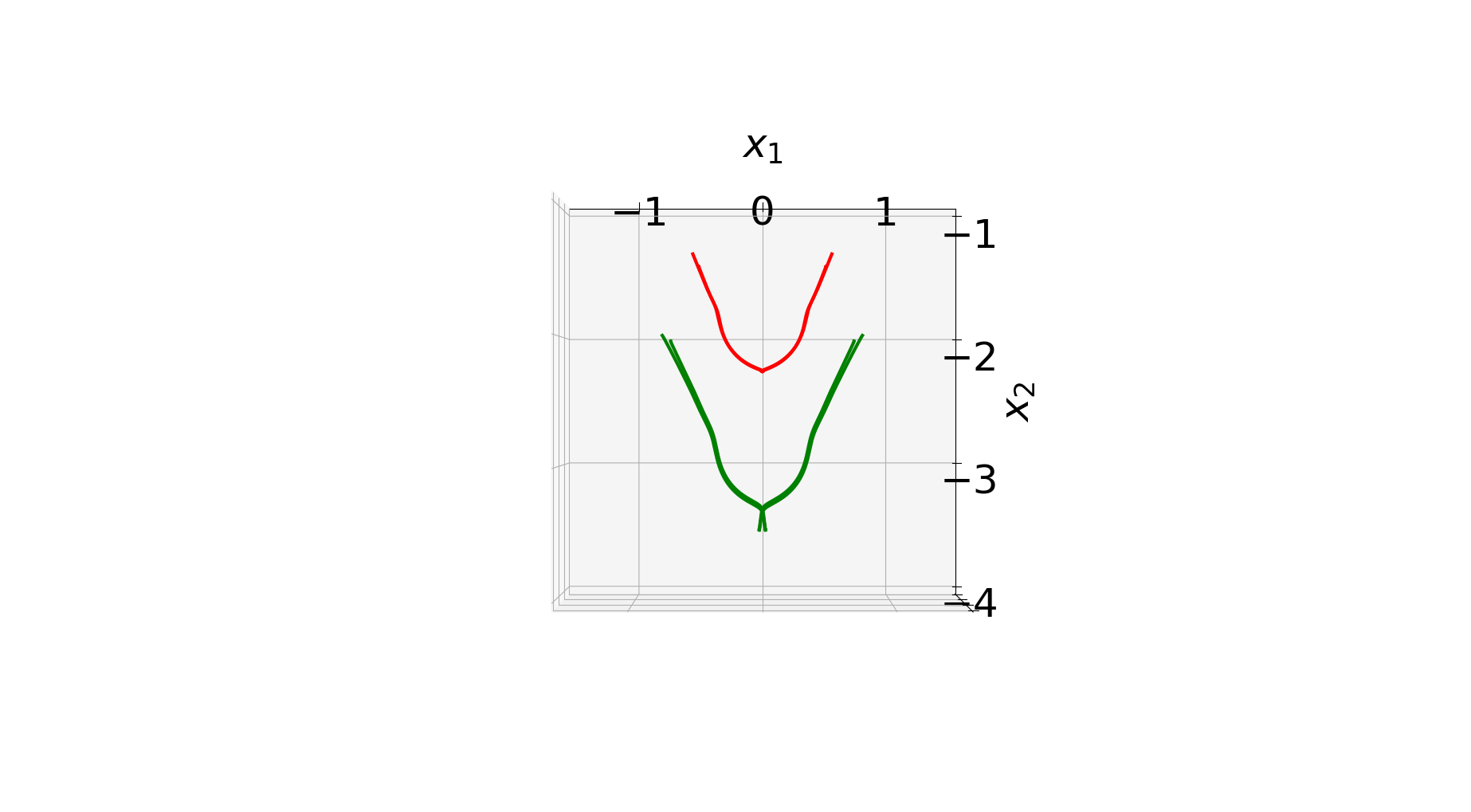}
        \caption{top view}
    \end{subfigure}
    \caption{Configuration of the vortices at time $t = 1.2$ when the horseshoe structure emerges for different $\varepsilon$ 
    ($\varepsilon = 0.5$ -- red, $\varepsilon = 1$ -- green) and $r_c = 0.05$}
    \label{fig:eps_just_after_rec}
\end{figure}

\begin{figure}
    \centering
    \begin{subfigure}[c]{0.32\textwidth}
        \centering  
        \includegraphics[width=\textwidth, trim={15cm 3cm 10cm 4cm}, clip]{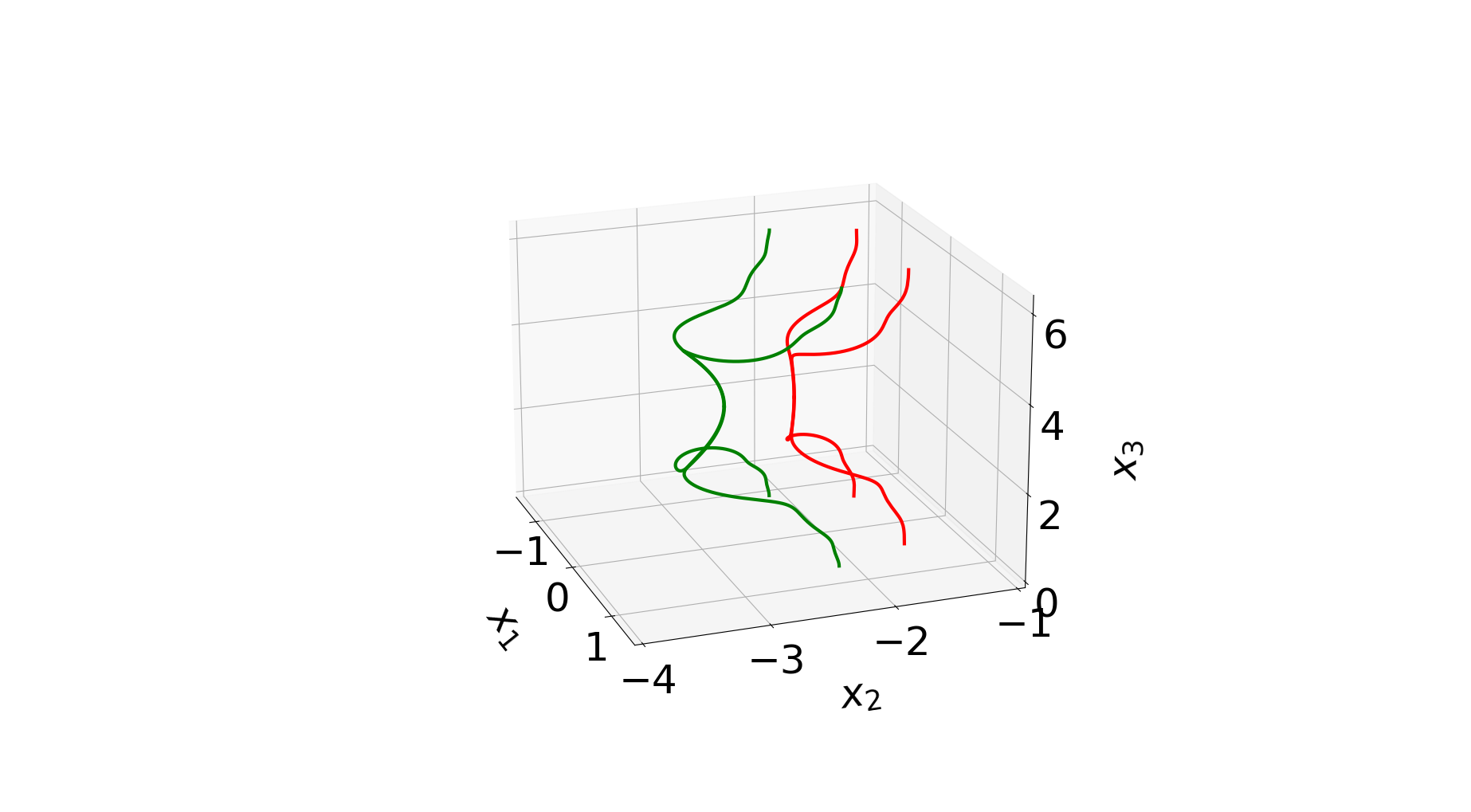}
        \caption{3D view}
    \end{subfigure}
    \begin{subfigure}[c]{0.32\textwidth}
        \centering      
        \includegraphics[width=\textwidth, trim={15cm 3cm 10cm 4cm}, clip]{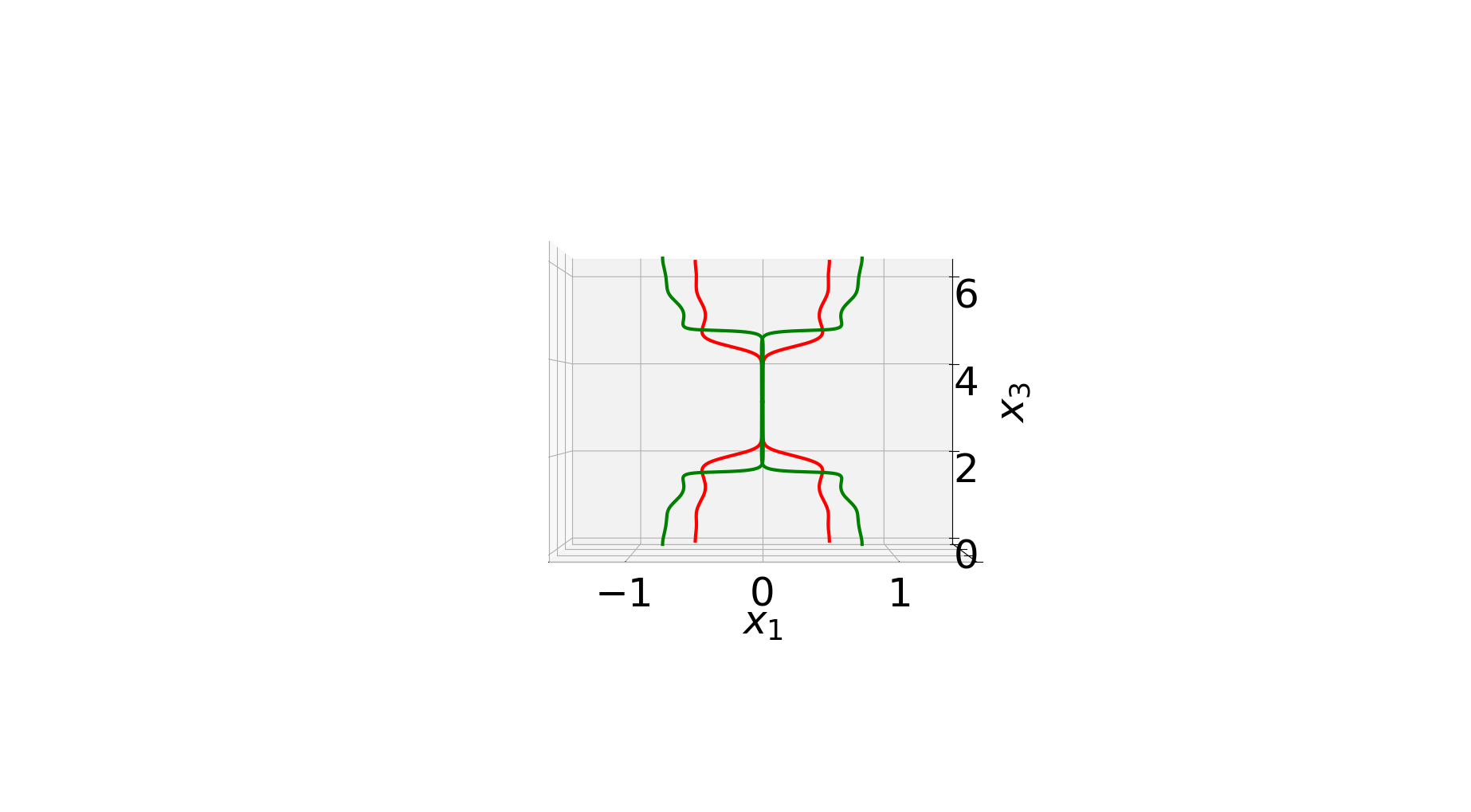}
        \caption{front view}
    \end{subfigure}
    \begin{subfigure}[c]{0.32\textwidth}
        \centering          
        \includegraphics[width=\textwidth, trim={15cm 3cm 10cm 4cm}, clip]{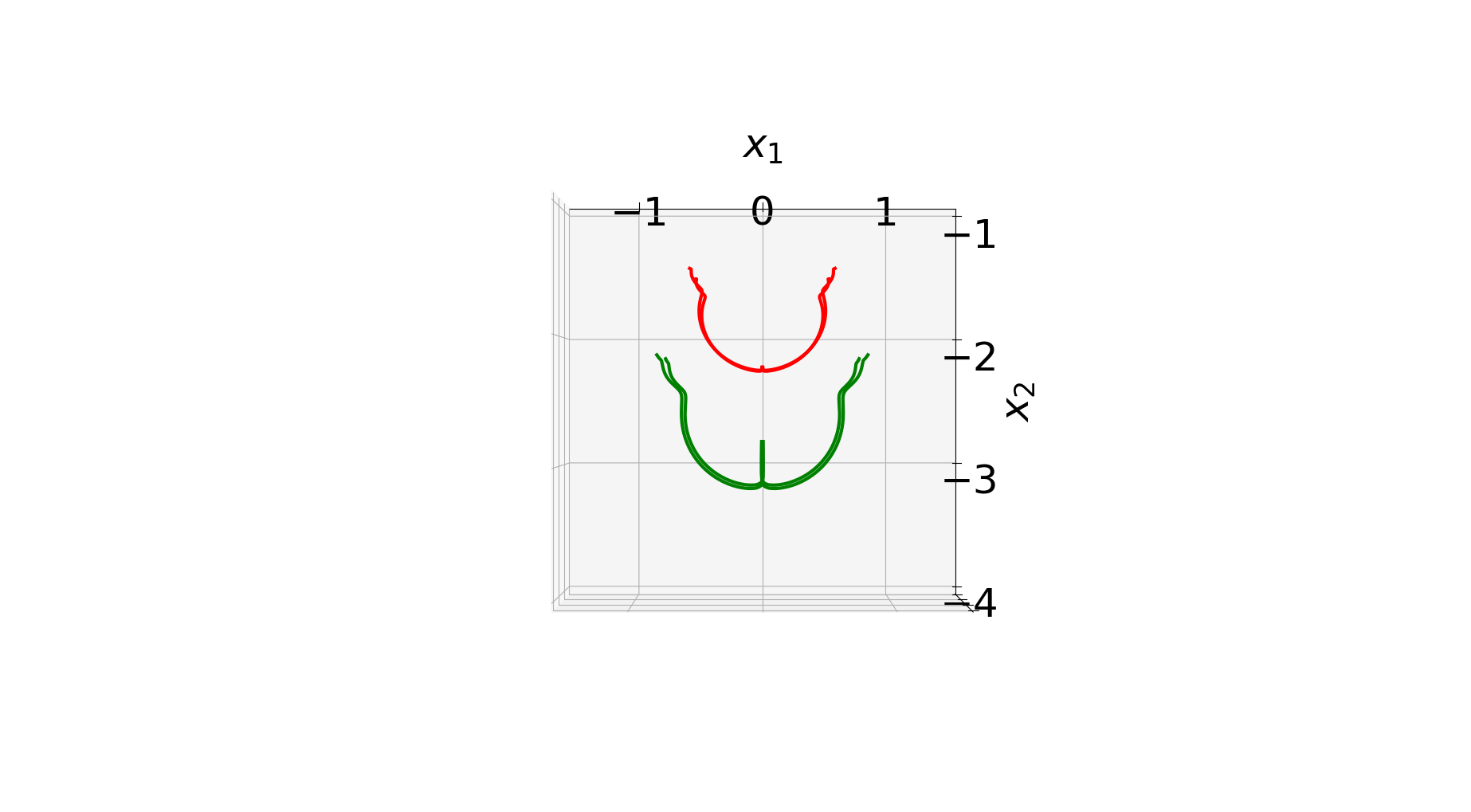}
        \caption{top view}
    \end{subfigure}
    \caption{Configuration of the vortices after the reconnection at time $t = 1.3$ for different $\varepsilon$ 
    ($\varepsilon = 0.5$ -- red, $\varepsilon = 1$ -- green) and $r_c = 0.05$}
    \label{fig:eps_after_rec}
\end{figure}

The configuration of the vortices at different time moments are depicted in figures~\ref{fig:eps_before_rec},\ref{fig:eps_just_after_rec}, and~\ref{fig:eps_after_rec}.
In the first figure the reconnection has not happened yet. However, we can see that the vortices with bigger values of $\varepsilon$ move faster
in $x_2$ direction and also the amplitude of Crow wave is bigger due to the bigger initial distance $b$. In the figure~\ref{fig:eps_just_after_rec}
the moment when the horseshoe structure emerges is depicted, but there are still no helical waves. For the smaller value of $\varepsilon$
the horseshoe structure is not planar and has a sharper tip. The configuration with helical waves and the horseshoe structure is shown in figure~\ref{fig:eps_after_rec}. For 
the large value of $\varepsilon$ the horseshoe is almost planar and the vortices look very similar to the shape we can see in evolution of a corner vortex in figure~\ref{fig:corner}.

\begin{figure}
    \centering
    \begin{subfigure}[c]{0.47\textwidth}
        \centering
        \includegraphics[width=\textwidth, trim={1cm 0cm 3cm 3cm}, clip]{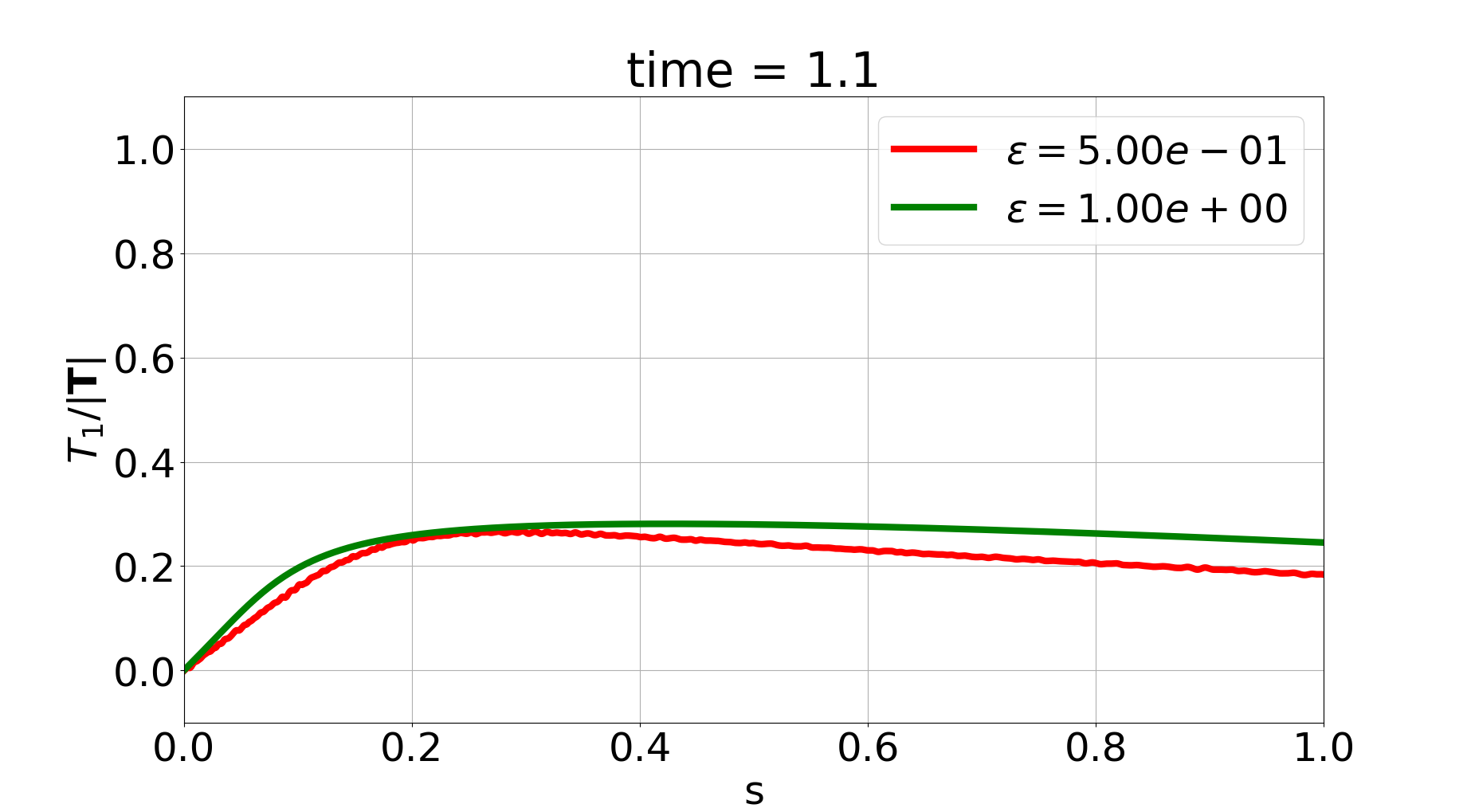}
        \caption{$t = 1.1$}\label{fig:T1T_ratio_t11}
    \end{subfigure}    
    \begin{subfigure}[c]{0.47\textwidth}
        \centering
        \includegraphics[width=\textwidth, trim={1cm 0cm 3cm 3cm}, clip]{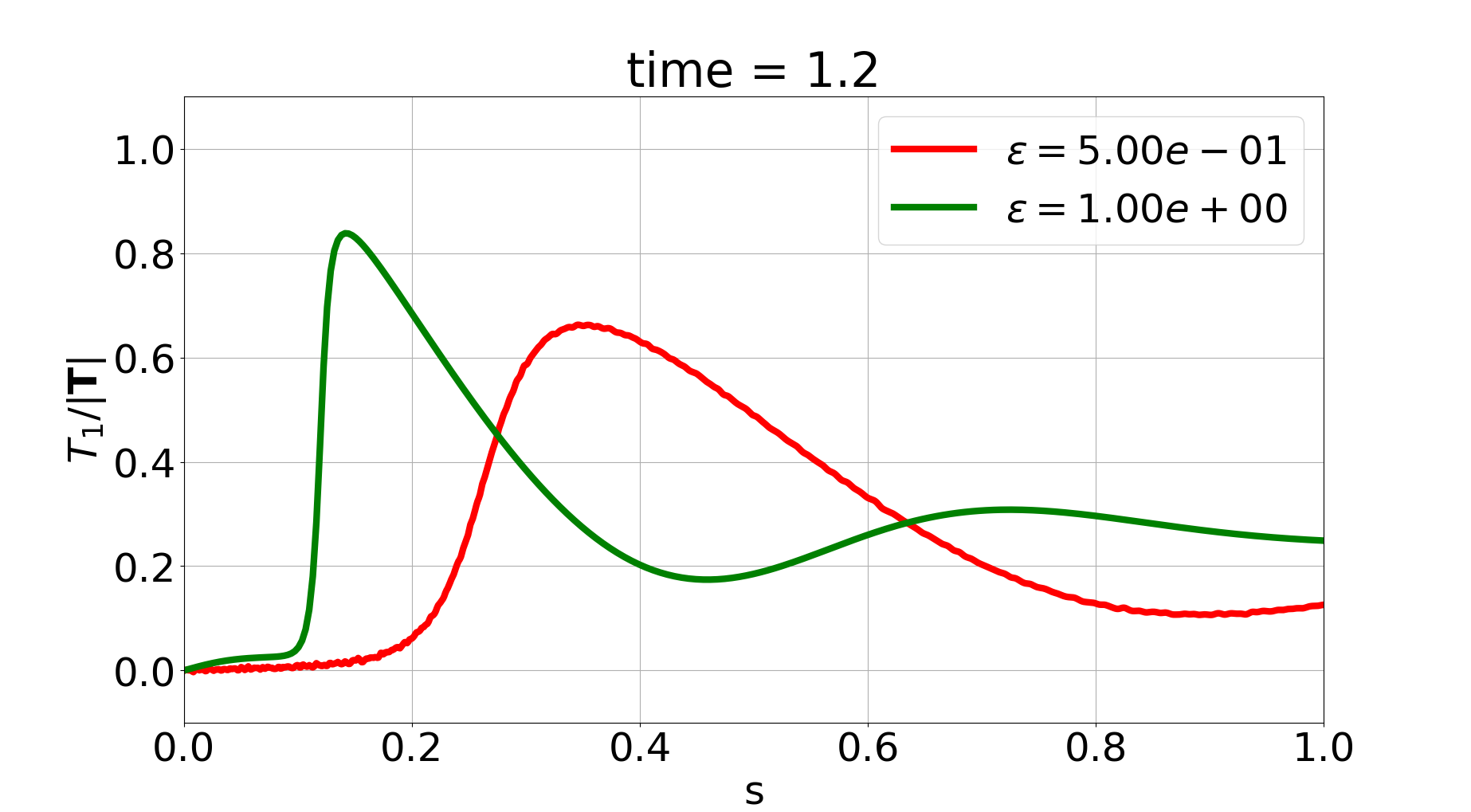}
        \caption{$t = 1.2$}\label{fig:T1T_ratio_t12}
    \end{subfigure}
    \begin{subfigure}[c]{0.47\textwidth}
        \centering
        \includegraphics[width=\textwidth, trim={1cm 0cm 3cm 3cm}, clip]{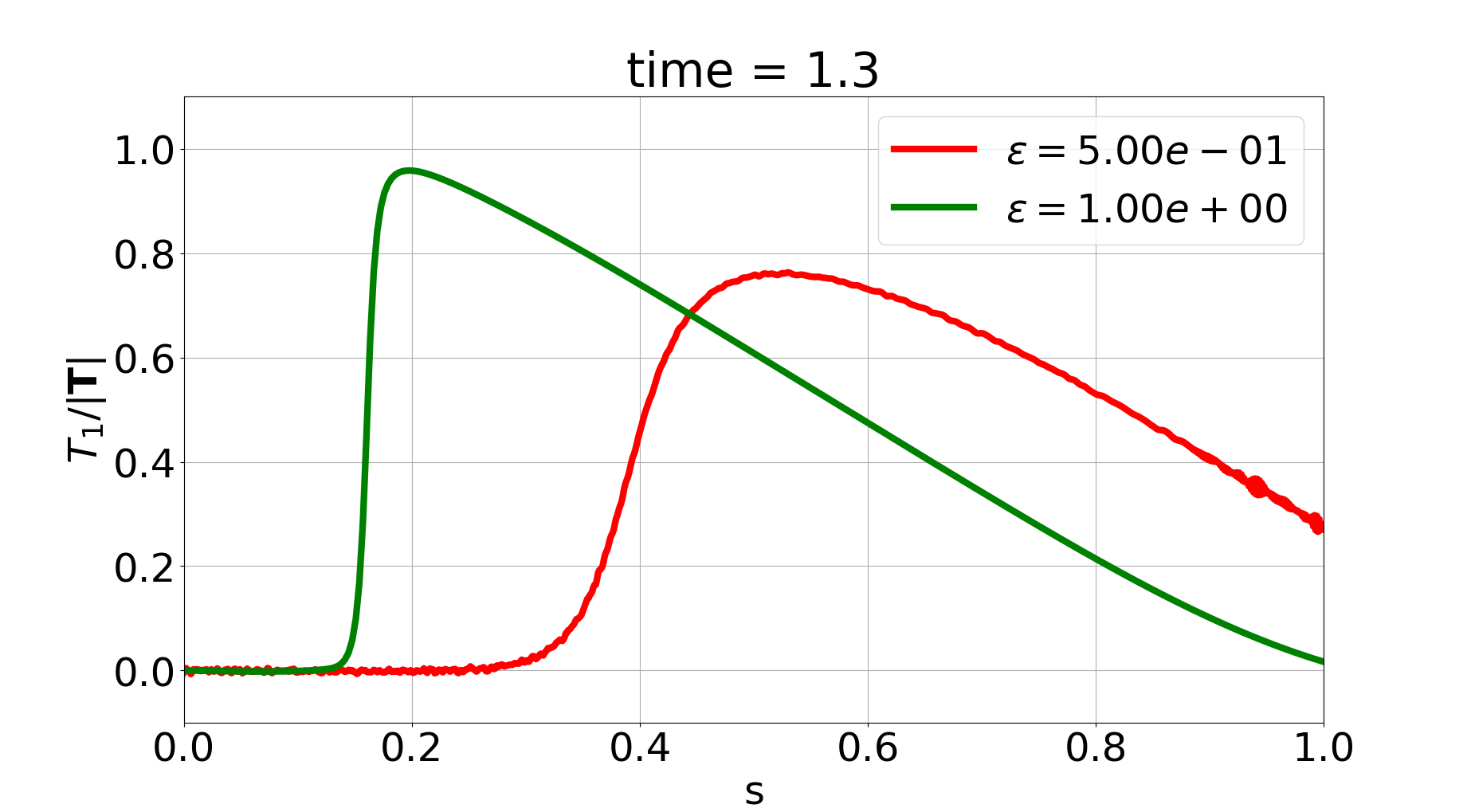}
        \caption{$t = 1.3$}\label{fig:T1T_ratio_t13}
    \end{subfigure}
    \caption{The ration $T_1 / \|\mathbf{T}\|$ for different time moments, different values of $\varepsilon$ and $r_c = 0.05$}
    \label{fig:T1T_ratio}
\end{figure}

Another important question one may ask is the behavior of the vortices around the reconnection. In particular, is there a cusp? Or can we expect a smooth
horseshoe? In order to do it we can analyze the components of the tangential vector $\mathbf{T}$. The bigger value of the first component $T_1$ corresponds to a smoother horseshoe. 
Moreover, when $\mathbf{T}$ is parallel to $\mathbf{e}_1$ the cusp disappears. A very rough analytical result on the ratio $T_1 / |\mathbf{T}|$ is presented in 
the section~\ref{sec:analysis}. Using numerical simulations we can see that the value of $T_1 / |\mathbf{T}|$ may be very close to $1$ for some 
$\varepsilon$. In figure~\ref{fig:T1T_ratio_t11} it is possible to see that before the reconnection the projection of normalized vector $\mathbf{T}$ to the direction
$\mathbf{e}_1$ is small because it is mostly oriented in the $\mathbf{e}_3$ direction. However, in figure~\ref{fig:T1T_ratio_t12} when the reconnection already
happens we can see that there is a bridge, where $T_1$ is close to $0$, and a jump almost of size $1$ at tip of the horseshoe. Moreover, the bigger is $\varepsilon$
the bigger is the maximum of $T_1 / |\mathbf{T}|$. This coincides with the analytical result from the section~\ref{sec:analysis}. In figure~\ref{fig:T1T_ratio_t13} 
we can see that the bridge is growing and the horseshoe persists.

\section{Comparison with the eye-shaped vortex}~\label{sec:eye-shaped}

We can compare the evolution of the reconnected vortices with an isolated vortex which has an eye-shape and
deforms obeying the local induction approximation. These not planar vortex can be considered as an approximation
of the shape after reconnection (see figure~\ref{fig:ns_horseshoe}) or as a "two-corner curvilinear polygon". We will see that
the evolution of the eye-shaped vortex has similarities with both the reconnection process and the evolution of polygonal 
vortex. It means that even though we can not see the corner at the reconnection time, the further evolution
of the vortices has a similar structure to the corner vortex (figure \ref{fig:corner}) under LIA.



\begin{figure}
    \centering
    \includegraphics[width=0.7\textwidth, trim={1cm 0cm 3cm 1cm}, clip]{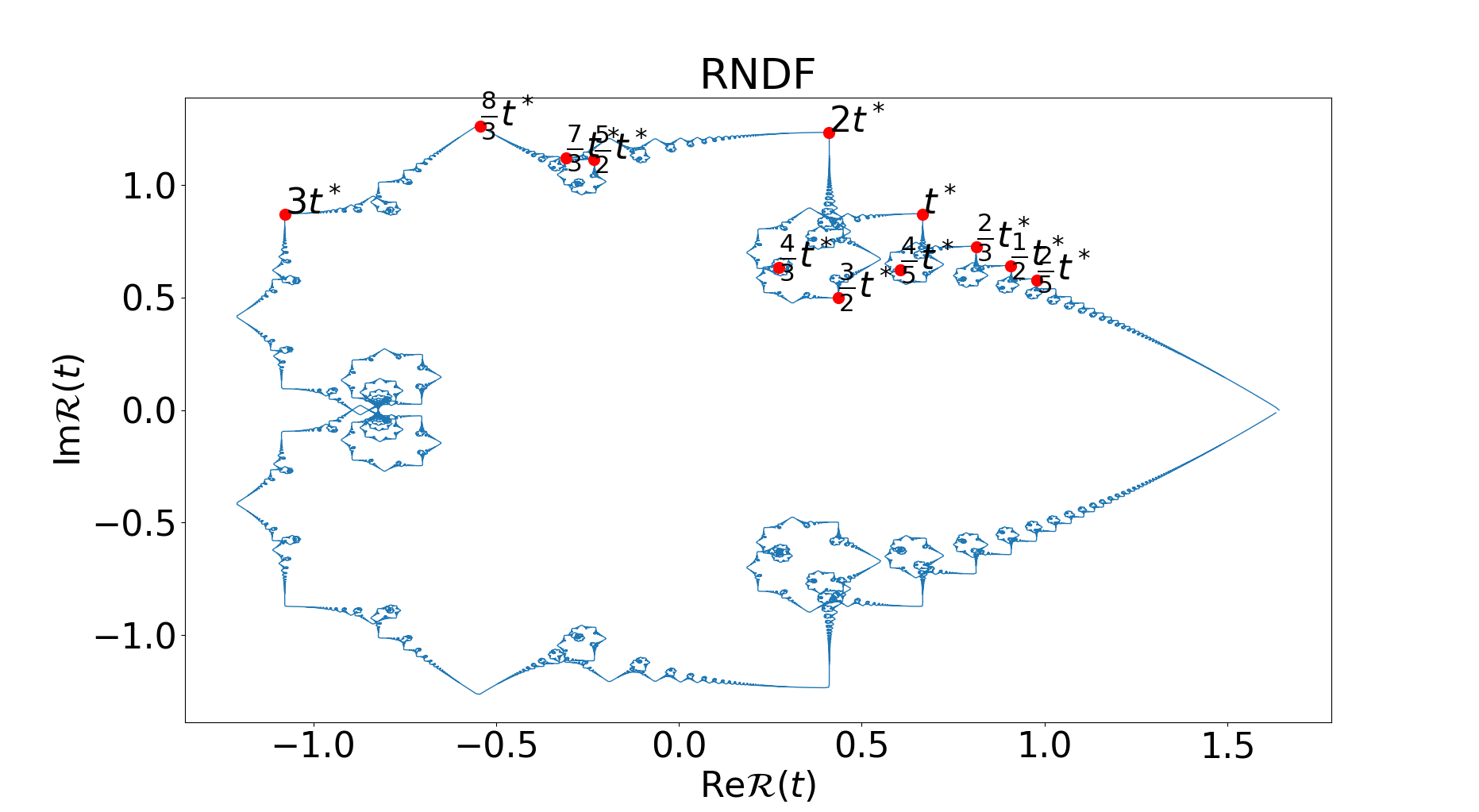}
    \includegraphics[width=0.7\textwidth, trim={1cm 0cm 3cm 1cm}, clip]{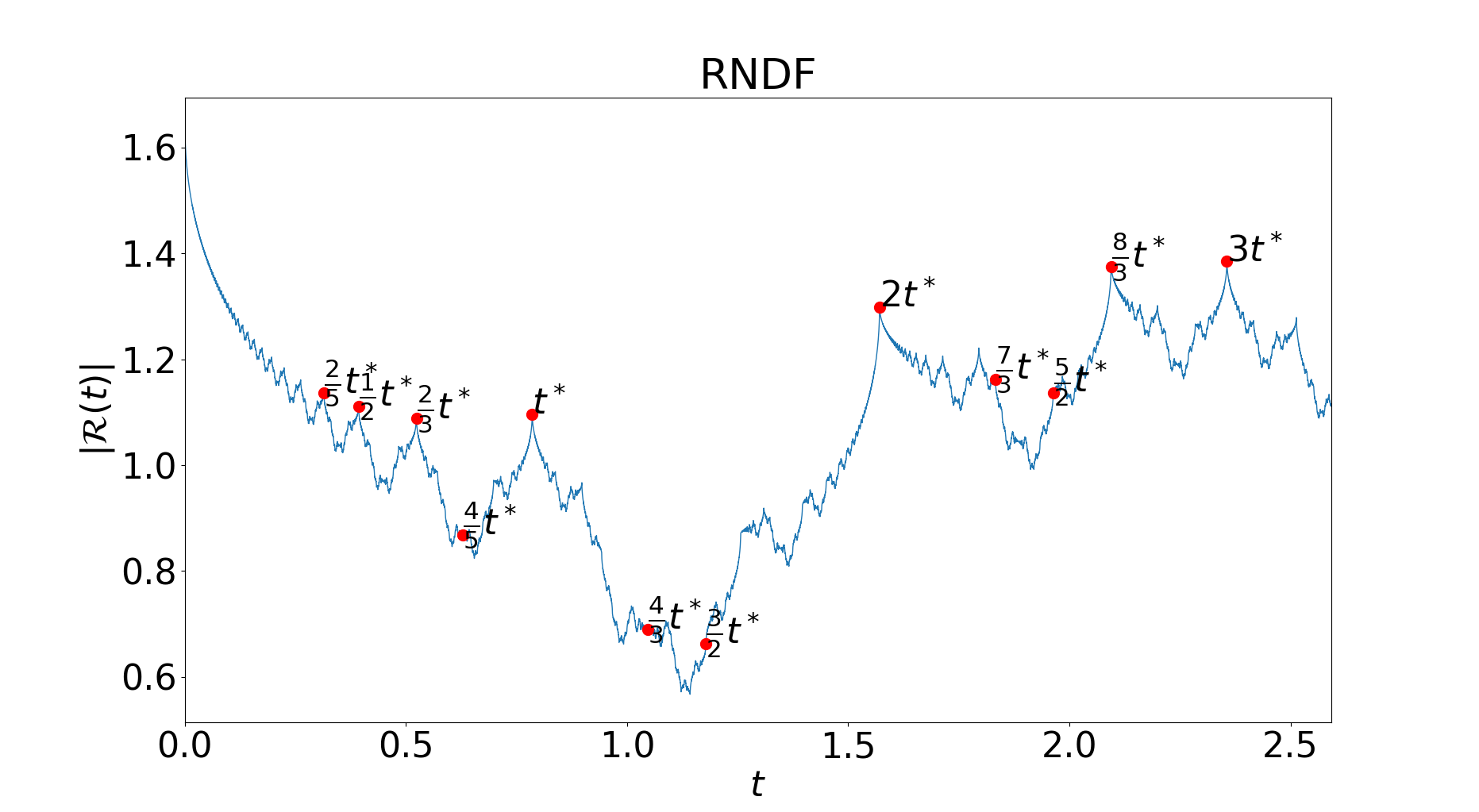}
    \caption{Riemann's non-differential function~\eqref{eq:rndf} on the complex plane and its absolute value}
    \label{fig:rndf}
\end{figure}

It is known~\cite{delahoz2018,jerrard2015} that a polygonal vortex with $M$ corners has a periodic behavior with period $T = 2\pi / M^2$ whereas 
the trajectory of the corner tends to a modification of Riemann's non-differential function (RNDF):
\begin{equation}
\mathcal{R}(t) = \sum_{k=1}^\infty \frac{e^{i t k^2}}{k^2}\label{eq:rndf}
\end{equation}
when $M$ tends to infinity, $i^2 = -1$. RNDF is a periodic multifractal that has a peculiar behavior at points corresponding
to rational multiples of the period~\cite{banica2022}. It can be seen in figure~\ref{fig:rndf}. We pick a point $t^* = \pi / 4$
and multiply it by different rational numbers. It is easy to see that most of them corresponds to local minima and maxima of
the absolute value of~\eqref{eq:rndf} and to corners of the graph in the complex plane. However, points that correspond
to the rational multiples with odd denominator fall into a cusp singularity of the absolute value of RNDF seen as a spirale on the complex plane. It was
shown in~\cite{delahoz2018} that the polygonal vortices have a similar behavior. Furthermore, at each rational multiple $p / q$ of the period
the shape of the corresponding vortex is also polygonal though not necessary planar. Moreover, the angle between two adjacent sides tends to $\pi$
when $q$ tends to infinity. This dependence on $q$ is more visible in the $H_\delta$ process studied in~\cite{kumar2022}.

The similar effect can be seen for the eye-shaped vortex which can be considered as a curvilinear polygon with two corners. The eye-shaped vortex is
more similar to the configuration we have after the reconnection than a polygonal one therefore the comparison with it is interesting to us.
In appendix~\ref{sec:eye-shaped-details} it is shown that the evolution of this vortex is quasi-periodic. Moreover, the Fourier analysis of the
trajectory that starts in a corner shows that the dominating frequencies are still the squares like in the polygonal case.

\paragraph*{The fluid impulse of the reconnected vortices}
We define the fluid impulse around the corner as
\begin{equation}
    \mathcal{F}_l(t) = \frac{1}{2}\int_{-l/2}^{l/2} \left(\mathbf{X}(q,t) - \mathbf{X}_0(t)\right)\wedge\mathbf{X}_s(q,t)d q,\label{eq:fluid_impulse}
\end{equation}
where the corner is located at $s = 0$ and the interval $l$ is $20\%$ of the whole perimeter of the vortex. The fluid impulse
is important for us because in the case of the reconnection we can not extract the singularity point hence we can not find its trajectory,
whereas the formula~\eqref{eq:fluid_impulse} is always applicable. It depends however on the position of the origin $\mathbf{X}_0$ which
in the case of the eye-shaped vortex we define as $\mathbf{X}_0(t) = \mathbf{X}(0,t)$. We will be mainly interested in the oscillations
and multi-fractal behavior for which the choice of origin does not have any influence. 

\begin{figure}
    \centering
    \includegraphics[width=0.7\textwidth, trim={0mm 0cm 15mm 0mm}, clip]{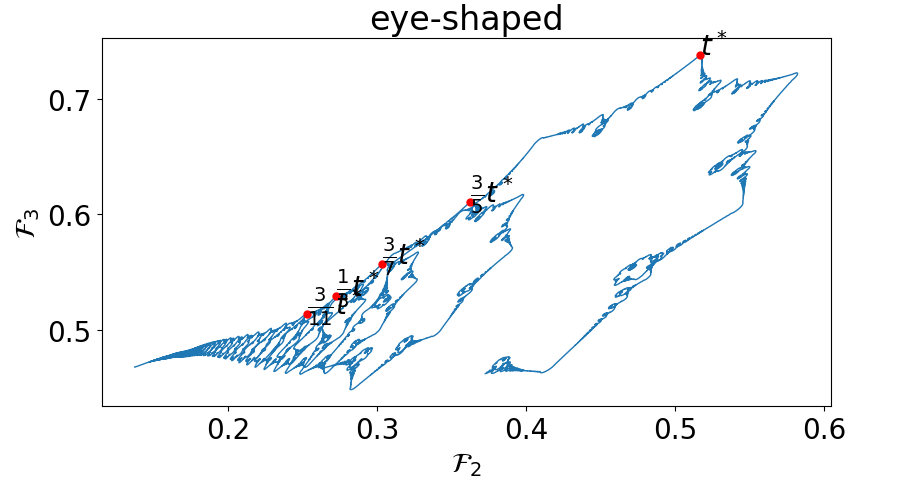}
    \includegraphics[width=0.7\textwidth, trim={0mm 0cm 15mm 0mm}, clip]{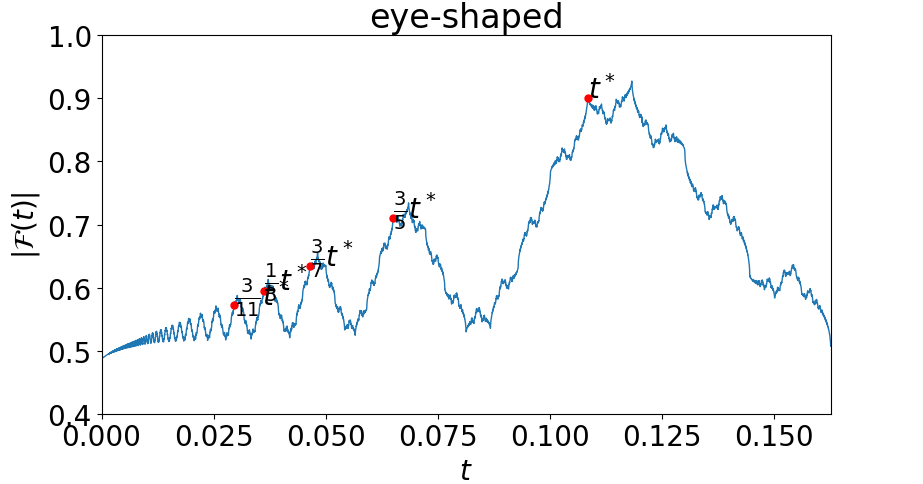}
    \caption{Fluid impulse of an eye-shaped vortex for $b = 0.4$, $\theta = \pi / 6$ and its modulus. Some
    maxima a located in times $3/11 t^*$, $1/3 t^*$, $3/7 t^*$, $3/5t^*$, $t^*$, and the bigger is the denominator,
    the smaller is the peak at that (similar to what happens in~\cite{kumar2022}).}
    \label{fig:eye_fi_rationals}
\end{figure}

Despite this similarity to the RNDF, the real structure of the fluid impulse of the eye-shaped vortex is much more complicated, 
as can be seen in figure~\ref{fig:eye_fi_rationals}. The curve is still planar since $\mathcal{F}_1(t) = 0$ for 
all time, the rational multiples of a local maxima $t^* = 0.10848$ also corresponds to the local maxima, minima,
and singular points but the scaling of self-similar structures are deformed and there is no real period. The dependence on 
the size of the maxima with respect to the size of the denominator is also presented.

\begin{figure}
    \centering
    \includegraphics[width=\textwidth, trim={1cm 0cm 3cm 1cm}, clip]{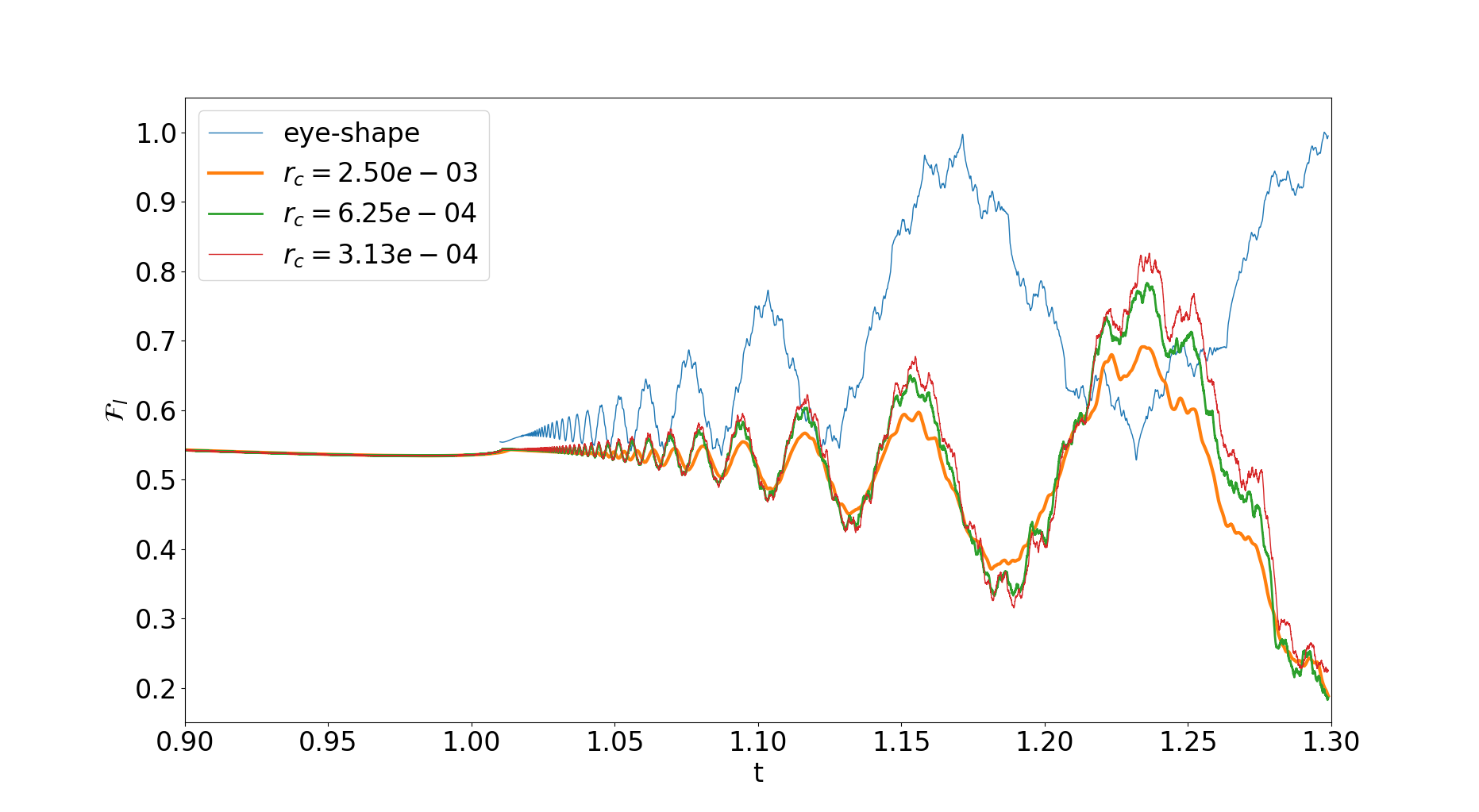}
    \caption{Fluid impulse~\eqref{eq:fluid_impulse} for the reconnection problem with different values of $r_c$,
             $\varepsilon = 0.05$, and $b = \sqrt{\varepsilon}$. The blue line corresponds to normalized and 
             shifted fluid impulse of the eye-shaped vortex with $\theta = \pi / 6$ and $b = 0.4$. At time $t\approx 1.01$
             the behavior suddenly changes from monotone growth to oscillation. We can use it as the definition of
             the reconnection time.}\label{fig:fi_rec}
\end{figure}

In the case of the reconnection we can not specify the reconnection point and follow its trajectory. However, we can calculate the fluid
impulse~\eqref{eq:fluid_impulse} around the reconnection region. The results for different $r_c$ and $\varepsilon = 0.05$ are depicted in
figure~\ref{fig:fi_rec}. Analysing the fluid impulse we can detect a sudden change from monotone to oscillating behavior at time $t\approx 1.01$.
We call this moment  the reconnection time. Note that from the configuration of the vortices we could not to define this time, thus analysis
of integral quantities such as the fluid impulse is beneficial for understanding of the reconnection phenomena. Moreover, one can note that the behavior after
the reconnection time is quite reminiscent to the one of the eye-shaped vortex though the period is different. The smaller is the 
regularization parameter $r_c$, the more details we can see in the fluid impulse. 

The extraction of squares however is not possible for the reconnection problem due to the noise generated by the bridge. 
Analysis of distribution of minima, maxima, and singular points using wavelet transform have also faced problems 
related to the noise. Maxima for the considered types of signals corresponds to singular points that can 
be studied using the multifractal analysis~\cite{turiel2006, muzy1993}. The main idea of this approach is to construct the singularity 
spectrum $D(h)$, that is the function for which each H\"older exponent $h$ yields the Hausdorff dimension of 
the set of points where the function has this exponent. For RNDF the singularity spectrum is known: 
$D(h) = 4 h - 2$ for $h\in[0.5, 0.75]$, $D(h) = 0$ for $h = 1.5$ and $D(h) = -\infty$ otherwise~\cite{jaffard1996}. 
We have tried to use the $p$-leaders method~\cite{wendt2007} to estimate the spectrum of RNDF and the fluid impulse.
Unfortunately these signals are quite delicate. Therefore, even in the case of RNDF the approximation of $D(h)$ is not very accurate 
(especially for exponents corresponding to $D(h) = 0$). For the fluid impulse of the polygonal or reconnected vortices the situation 
is much worse even though there is a tendency that $h$ corresponding to maximal $D(h)$ is decreasing getting closer to $0.75$.  

\begin{figure}
    \centering
    \includegraphics[width=0.7\textwidth, trim={1cm 0cm 3cm 1cm}, clip]{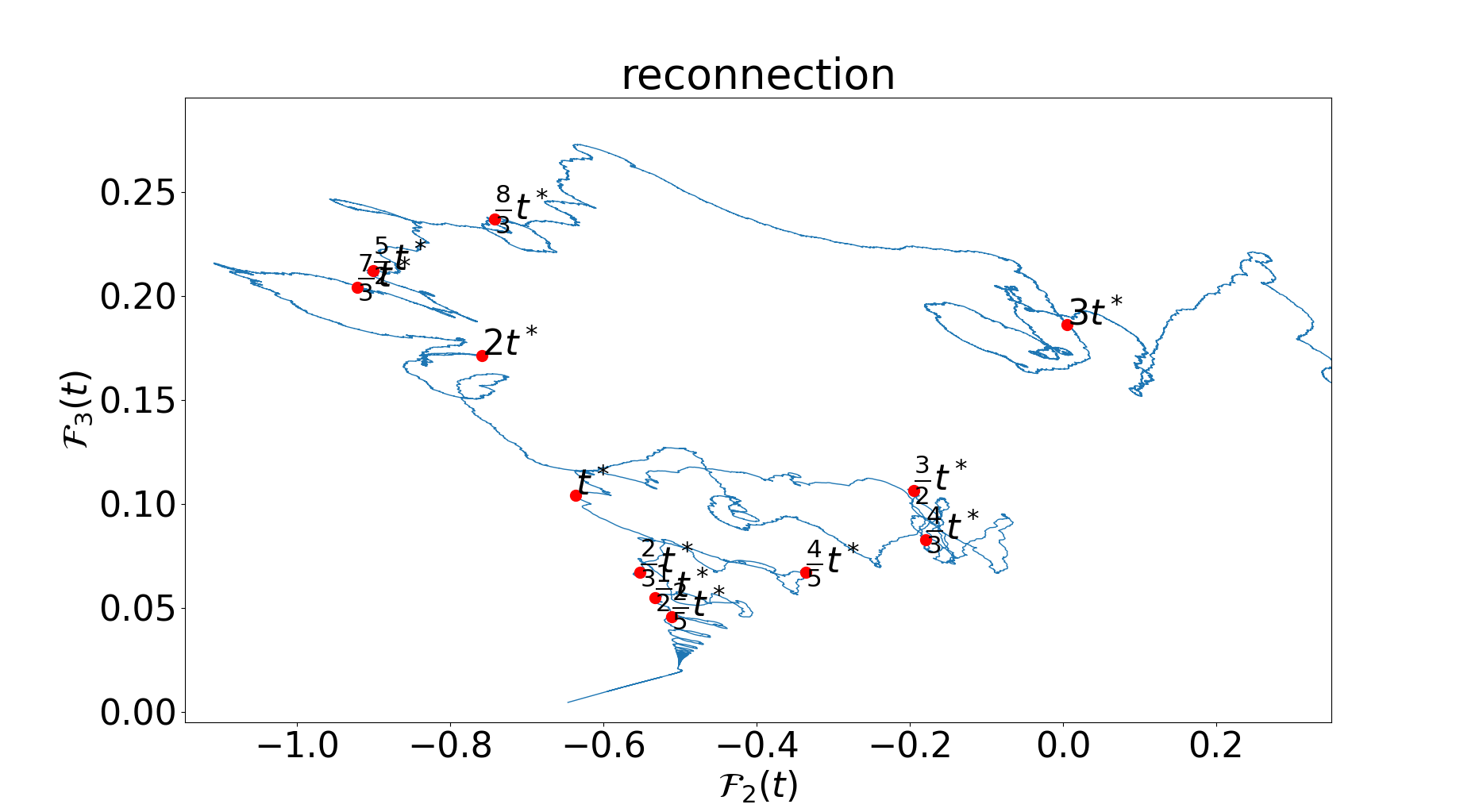}
    \includegraphics[width=0.7\textwidth, trim={1cm 0cm 3cm 1cm}, clip]{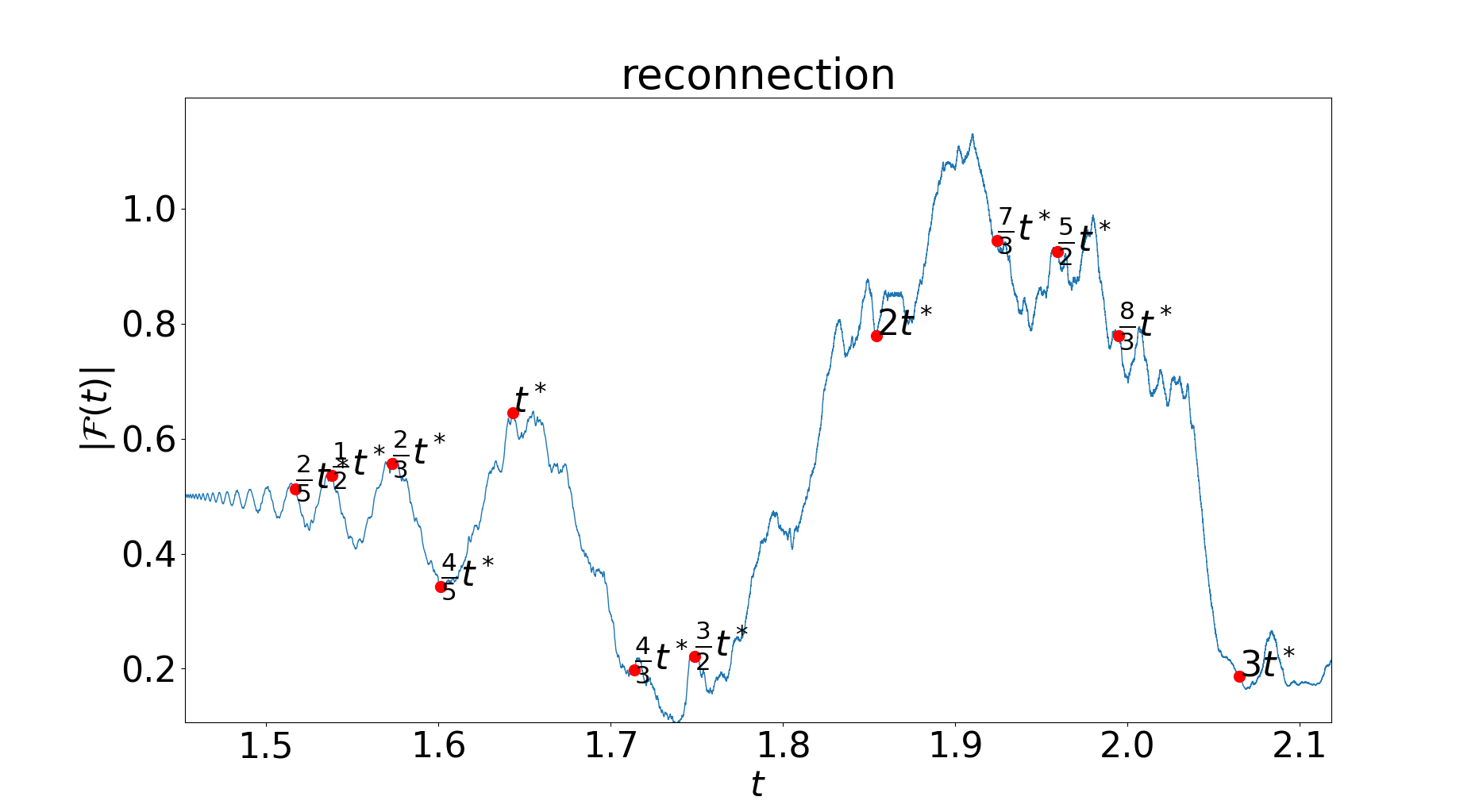}
    \caption{Fluid impulse  and its modulus for the reconnection problem when $\varepsilon = 0.03$, $b = 0.22$, $r_c = 6.25\cdot 10^{-4}$}
    \label{fig:rec_fi_rationals}
\end{figure}

Since the multifractal analysis failed, in this paper we perform only a qualitative analysis of similarity between the fluid impulse of 
the eye-shaped vortex and the one of the reconnection problem. We choose a point $t^*$ related to a local maxima of the reconnection 
fluid impulse and check if the rational multiples of this point also fall in maxima, minima, or singular points. 
The results are depicted in figure~\ref{fig:rec_fi_rationals} for $\varepsilon = 0.03$, $b = 0.22$, $r_c = 6.25\cdot 10^{-4}$. We can see
that the rational points mostly fall into local minima and maxima at least for a short time after reconnection. However later we 
can see that the self-similar structure of the fluid impulse vector is completely lost. In order to obtain a cleaner structure we have
to focus on a small time interval after the reconnection and choose a small regularization parameter $r_c$. This is quite challenging,
due to the stability condition~\eqref{eq:tau_condition} of the method. Thus, a new approach is required to find more similarities between
the reconnection process and the eye-shaped vortex.

\paragraph*{The vortex separation rate}

The scaling law which determines how the distance between vortices changes during the reconnection was studied in 
multiple works~\cite{krstulovic2017,fonda2019,hussain2020}. It is not completely clear if there are different laws before and after
the reconnection or if this law varies for quantum and classical fluids. Nevertheless, there are multiple evidences that the separation rate
$\delta(t)$ of the vortices after the reconnection should be of the scale $\sqrt{t - t_{rec}}$ where $t_{rec}$ is the reconnection time.
This rate can be observed in experiments~\cite{fonda2019} and also coincide with the analytical result for the corner vortex~\cite{vega2003}.

\begin{figure}
    \centering
    \begin{subfigure}[c]{0.7\textwidth}
        \centering
        \includegraphics[width=\textwidth, trim={1cm 0cm 3cm 1cm}, clip]{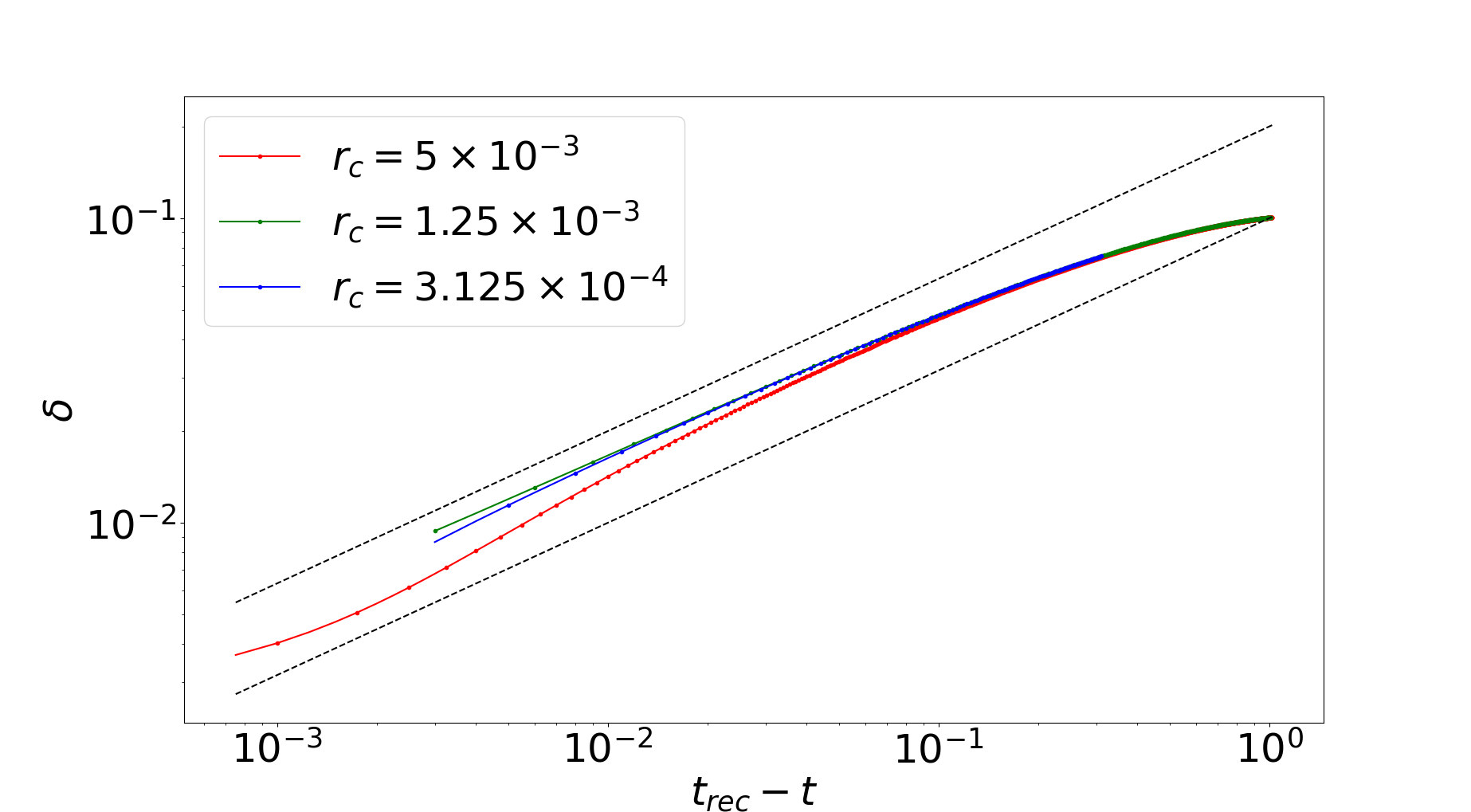}
        \caption{Before the reconnection}\label{fig:separation_rate_before}
    \end{subfigure}
    \begin{subfigure}[c]{0.7\textwidth}
        \centering
        \includegraphics[width=\textwidth, trim={1cm 0cm 3cm 1cm}, clip]{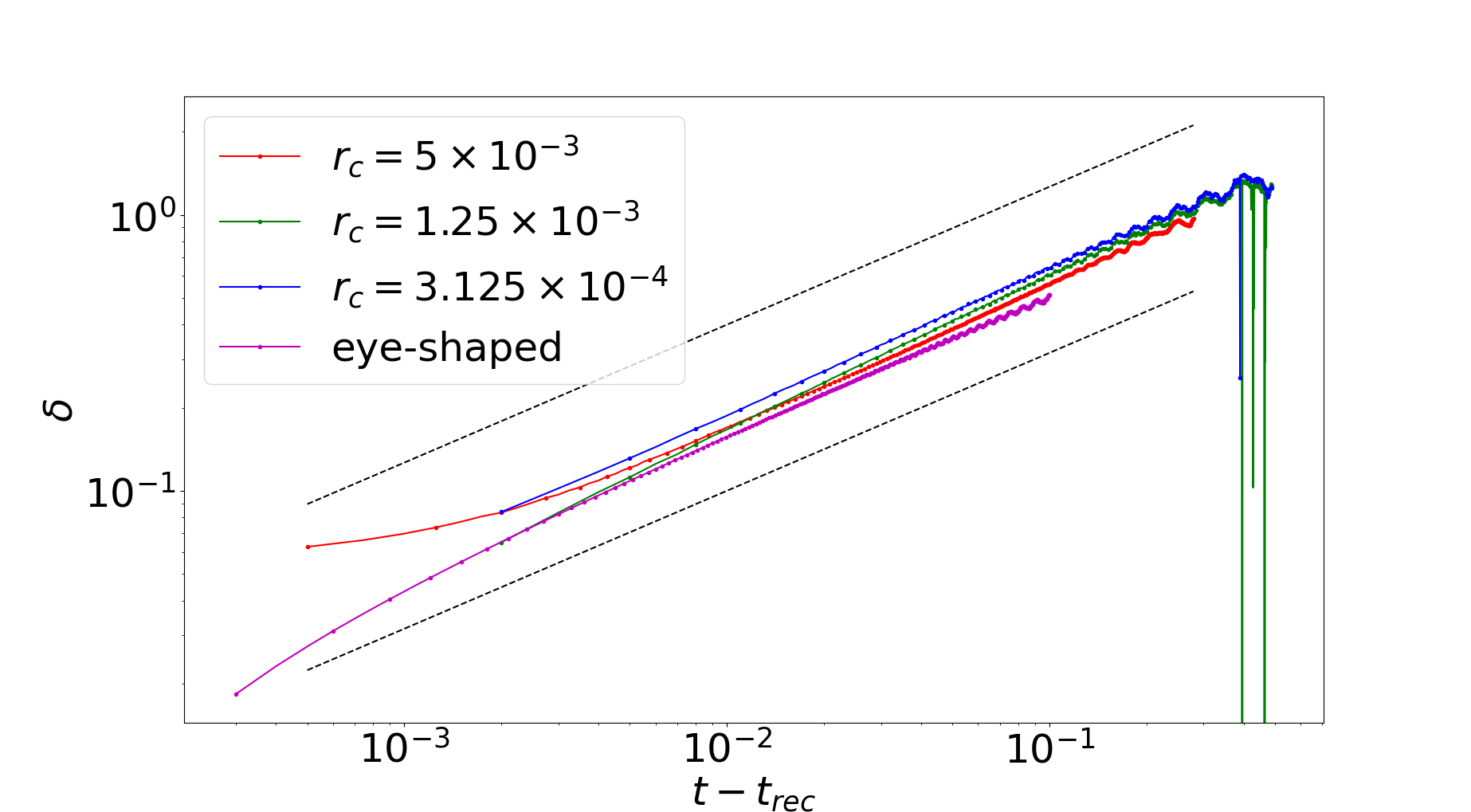}
        \caption{After the reconnection}\label{fig:separation_rate_after}
    \end{subfigure}        
    \caption{The minimal distance between vortices before and after the reconnection for $\varepsilon = 0.05$ and different values of $r_c$}
    \label{fig:separation_rate}  
    
\end{figure}

In figure~\ref{fig:separation_rate} the separation rates before and after the reconnection are depicted for $\varepsilon = 0.05$ different values of 
$r_c$. The black dashed lines corresponds to the scale $\sqrt{|t - t_{rec}|}$. We can see that before and after the reconnection the separation 
rate is very close to the square root law. The agreement is better for small values of $r_c$ that correspond to a case of smaller viscosity. 
We can also note that the $x_2$ component of the position of the eye-shaped vortex corner has the same square root timescale as the reconnected vortices (figure~\ref{fig:separation_rate_after}). This result can be considered as another evidence that the behavior of the vortices after the reconnection resembles the evolution
of the corner vortex even though we can not see the corner at the reconnection time due to the presence of the regularization parameter $r_c$.

\section{Conclusions}\label{sec:conclusions}

Even though it is not entirely clear what happens in the reconnection time the further
evolution of vortices contains coherent structures reminiscent to the ones generated by a polygonal vortex~\cite{delahoz2018, jerrard2015}. 
In order to extract and analyze these structures a new model~\eqref{eq:main_equation} describing the interaction of a pair of antiparallel 
vortices is developed.  
The introduction of the regularization parameter $r_c$ allows us to go beyond the reconnection time. Moreover, the model 
provides a closed expression~\eqref{eq:arclength} for the length of the tangent vector of the vortex central line.
This length is proportional to the circulation and, according to our result, it is increasing, when
the distance between vortices tends to $0$. It can be considered as the vortex stretching phenomenon.
The model also predicts the Crow waves firstly described in~\cite{crow1970} and the formation of coherent 
structures. These structures are a pair of horseshoes in the spirit of~\cite{schwarz1985, vega2003} connected by
a bridge which is artificial but can not be removed due to restrictions of the model. Nevertheless,
we have shown analytically in section~\ref{sec:analysis}, that the non smoothness of the horseshoe due 
to the bridge can be bounded. Finally, the model predicts the square root timescale for the separation rate after
the reconnection.

We performed a numerical simulation for different values of the vortex interaction parameter $\varepsilon$ and the 
regularization parameter $r_c$. There is a difference with the evolution of the corner vortex filament shown in 
figure~\ref{fig:corner} and the reconnection. In the first case the horseshoe emerges immediately and can 
be infinitely small, whereas for our model the smallest possible size of the horseshoe is finite and dictated by
the parameter $r_c$. This effect makes impossible to determine the reconnection time, so we can speak only about 
the reconnection interval which is getting smaller when $r_c$ tends to zero. However, the reconnection time can
be defined better if an integral quantity instead of the configuration of the vortices is considered. We focus on
the fluid impulse that is an integral of the cross product between the position and the vorticity calculated in the
reconnection region. It is possible to see a sudden change of the behavior of this quantity from monotone 
to oscillatory happening at the reconnection time (see figure~\ref{fig:fi_rec}). Note also that the smaller $r_c$
the more complicated the behavior after the reconnection is. The oscillations look quite reminiscent to
the ones of the fluid impulse of the polygonal vortex~\cite{delahoz2018} that it tends to Riemann's 
non-differential function (RNDF, formula~\eqref{eq:rndf}) when the number of polygon sides tends to infinity. It can be considered as
an evidence that the antiparallel vortices indeed form a corner at the reconnection. However, a further
research is necessary for estimation of the noise produced by the bridge, studying the influence of the interaction,
and calculation of the corner angle. 

A possible way to find more similarities between the fluid impulse of reconnecting vortices and RNDF is to
improve the numerical method. In particular the condition~\eqref{eq:tau_condition} is very restrictive because
we have to reduce both space and time discretization to be able to solve problems with small $r_c$. Probably this
restriction can be surpassed by choosing right orders in space and time discretizations. Another way
is to apply filtering to the already obtained fluid impulse or try more advanced methods of the analysis.
We have tried to apply methods of multifractal analysis such as the p-leader method~\cite{wendt2007}. This approach
is based on studying the distribution of the singular points of the function and the calculation of the singularity spectrum
which can be considered as a fingerprint of the signal. Similar singularity spectrum means similar multifractal
properties of the signals. Since, the RNDF is very delicate for such methods, and the fluid impulse of the reconnecting
vortices have a lot of noise, our analysis could not arrive to any conclusion. 

Advances in the study of the vortex reconnection may have a huge impact to the understanding of 
turbulence and its structure.


    \appendix

\section{Equations in the generalized Frenet frame}\label{sec:frenet}

The equations~\eqref{eq:main_equation} can also be rewritten in the generalized 
Frenet frame following~\cite{how1998}. In this part in we use $\mathbf{T}$ to designate the normalized
tangential vector, not just the derivative $\mathbf{X}_s$. Consider the orthonormal 
frame in $\mathbb{R}^3$:
\begin{equation}
\mathbf{T}=\frac{\mathbf{X}_s}{|\mathbf{X}_s|},\quad 
\mathbf{N}_1 \perp \mathbf{T},\quad 
\mathbf{N}_2 = \mathbf{T}\wedge\mathbf{N}_1.\label{eq:basis}
\end{equation}
For this frame we can write the generalized Frenet system:
\begin{eqnarray}
&&\frac{\partial}{\partial q} \mathbf{T} = \kappa_1 \mathbf{N}_1 - \kappa_2 \mathbf{N}_2,\label{eq:frenet_T}\\
&&\frac{\partial}{\partial q} \mathbf{N}_1 = -\kappa_1 \mathbf{T} + \omega \mathbf{N}_2,\label{eq:frenet_N1}\\
&&\frac{\partial}{\partial q} \mathbf{N}_2 = \kappa_2 \mathbf{T} - \omega \mathbf{N}_1,\label{eq:frenet_N2}
\end{eqnarray}
where $q$ is the arclength parameter, that is $\frac{\partial}{\partial q} = \frac{1}{|\mathbf{X}_s|}\frac{\partial}{\partial s}$. 
In terms of the system~\eqref{eq:frenet_T}-\eqref{eq:frenet_N2} the curvature and the torsion can be calculated as
\begin{eqnarray}
\kappa = \sqrt{\kappa_1^2 + \kappa_2^2},\label{eq:curvature}\\
\tau = \omega + \frac{\kappa_2 \kappa_{1,q} - \kappa_1 \kappa_{2,q}}{\kappa^2},\label{eq:torsion}.
\end{eqnarray}
On the orther hand the vectors $\mathbf{T}$, $\mathbf{N}_1$, and $\mathbf{N}_2$ change in time following the system
\begin{eqnarray}
&&\frac{\partial}{\partial t} \mathbf{T} = -\lambda_1 \mathbf{N}_2 + \lambda_2 \mathbf{N}_1,\label{eq:time_T}\\
&&\frac{\partial}{\partial t} \mathbf{N}_1 = -\lambda_2 \mathbf{T} + \lambda_3 \mathbf{N}_2,\label{eq:time_N1}\\
&&\frac{\partial}{\partial t} \mathbf{N}_2 = \lambda_1 \mathbf{T} - \lambda_3 \mathbf{N}_1.\label{eq:time_N2}
\end{eqnarray}
We can recover coefficients $\lambda_i,\ i\in\{1,2,3\}$ using the equation~\eqref{eq:main_equation}. 
Indeed, in the frame~\eqref{eq:basis} the final system~\eqref{eq:main_equation} reads
\begin{equation}
\mathbf{X}_t = u\mathbf{N}_1 + v\mathbf{N}_2,\label{eq:system_basis}
\end{equation}
where coefficients are given by
\begin{eqnarray}
&&u = \kappa_2 + \varepsilon \frac{x_1}{x_1^2 + r_c^2} n_2,\label{eq:u_coeff}\\
&&v = \kappa_1 - \varepsilon \frac{x_1}{x_1^2 + r_c^2} n_1.\label{eq:v_coeff}
\end{eqnarray}
Here $n_1$ and $n_2$ are first components of the vectors $\mathbf{N}_1$ and $\mathbf{N}_2$ respectively. 
Taking derivative of \eqref{eq:system_basis} respect to $s$ we obtain
\begin{equation}
\mathbf{X}_{st} = \mathbf{T} L (-u \kappa_1 + v \kappa_2) + 
\mathbf{N}_1 (u_s - L \omega v) + 
\mathbf{N}_2 (v_s + L \omega u),\label{eq:xst_eq}
\end{equation}
where $L = |\mathbf{X}_s|$. On the other hand we can write that $\mathbf{X}_s = L \mathbf{T}$ 
and take derivative with respect to $t$:
\begin{eqnarray}
\mathbf{X}_{st} = L_t \mathbf{T} + 
\lambda_2 L \mathbf{N}_1 - 
\lambda_1 L \mathbf{N}_2.\label{eq:xst_frenet}
\end{eqnarray}
Since for $\mathbf{X}(s,t)$ the order of differentiation does not matter we can find the equation for $L_t$:
\begin{equation}
L_t = L (-u \kappa_1 + v \kappa_2) ,\label{eq:L_t}
\end{equation}
and also expressions for coefficients $\lambda_1$ and $\lambda_2$:
\begin{eqnarray}
&&\lambda_1 = -\frac{v_s}{L} - \omega u,\label{eq:lambda_1}\\
&&\lambda_2 = \frac{u_s}{L} - \omega v.\label{eq:lambda_2}
\end{eqnarray}

We can continue this process and obtain the equations for $\kappa_{1,t}$ and $\kappa_{2,t}$:
\begin{eqnarray}
&&\kappa_{1,t} = \frac{1}{L}\frac{\partial}{\partial s}\left(\frac{u_s}{L}\right) - 
\frac{\omega_s v + 2 \omega v_s}{L} - \omega^2 u + \kappa_1 (u \kappa_1 - v \kappa_2) - 
\lambda_3 \kappa_2,\label{eq:kappa_1_t}\\
&&\kappa_{2,t} = -\frac{1}{L}\frac{\partial}{\partial s}\left(\frac{v_s}{L}\right) - 
\frac{\omega_s u + 2 \omega u_s}{L} + \omega^2 v + \kappa_2 (u \kappa_1 - v \kappa_2) + 
\lambda_3 \kappa_1.\label{eq:kappa_2_t}
\end{eqnarray}
What do we need to close the system? In the expression for $u$ and $v$ we use $x_1$, $n_1$, and $n_2$ 
whose equations can be obtained from~\eqref{eq:system_basis},\eqref{eq:time_N1}, and~\eqref{eq:time_N2} respectively:
\begin{eqnarray}
&&x_{1,t} = u n_1 + v n_2,\label{eq:x_t},\\
&&n_{1,t} = -\left(\frac{u_s}{L} - \omega v\right) \frac{x_{1,s}}{L} + \lambda_3 n_2,\label{eq:n_1_t}\\
&&n_{2,t} = -\left(\frac{v_s}{L} + \omega u\right) \frac{x_{1,s}}{L} - \lambda_3 n_1.\label{eq:n_2_t}
\end{eqnarray}
Here we got rid of $\lambda_1$ and $\lambda_2$ using~\eqref{eq:lambda_1},\eqref{eq:lambda_2}. 
Expressions for $\omega$ and $\lambda_3$ are still missing. There are two ways to find the first quantity. 
First one is to use equations~\eqref{eq:frenet_N1} and~\eqref{eq:frenet_N2}:
\begin{equation}
\omega = \frac{n_{1,s} n_2 + \kappa_1 x_{1,s} n_2 - n_{2,s} n_1 + \kappa_2 x_{1,s} n_1}{L(n_1^2 + n_2^2)}.\label{eq:omega_stat}
\end{equation}
This is a functional equation, and it does not include $\lambda_3$. However, in order the denominator to be nonzero we have 
to require $n_1^2 + n_2^2 \ne 0$ that is $\mathbf{T} \ne \mathbf{e}_1$. Another way is to take a 
derivative of~\eqref{eq:frenet_N1} respect to $t$ and a derivative of~\eqref{eq:time_N1} respect to $q$ making them equal, 
that gives a differential equation for $\omega_t$:
\begin{equation}
\omega_t = \frac{u_s \kappa_2 + v_s \kappa_1 + \lambda_{3,s}}{L} + \omega (-v \kappa_2 + u \kappa_1).\label{eq:omega_t}
\end{equation}

In order to find $\lambda_3$ we have to make assumptions about our frame. Suppose that $n_2 = 0$ and use the equation~\eqref{eq:n_2_t}:
\begin{equation}
    \lambda_3 = -\left(\frac{v_s}{L} + \omega u\right)\frac{x_{1,s}}{L n_1} = -\left(\frac{v_s}{L} + \omega u\right)\frac{x_{1,s}}{\sqrt{L^2 - x_{1,s}^2}}.\label{eq:lambda_3_interaction}
\end{equation}
Then the direct expressions for $\mathbf{N}_1$ and $\mathbf{N}_2$ are:
\begin{equation}
    \mathbf{N}_1 = \frac{T_1 \mathbf{T} - \mathbf{e}_1}{\sqrt{1 - T_1^2}},\quad 
    \mathbf{N}_2 = \frac{\mathbf{e}_1 \wedge \mathbf{T}}{\sqrt{1 - T_1^2}}. \nonumber
\end{equation}
To make the first vector not zero the tangential vector should not be oriented in the $\mathbf{e}_1$
direction. Numerical experiments show that if $r_c > 0$ this is true before the reconnection and a long time 
after it. The second vector corresponds to the direction of the interaction term in the
system~\eqref{eq:main_equation}. This choice of frame vectors has another advantage that
we can express $\kappa_1$ using $L$, $x_1$ and their derivatives. Indeed, calculating
inner product of~\eqref{eq:frenet_T} with $\mathbf{N}_1$ and taking into account that 
$\mathbf{T}_q \cdot \mathbf{T} = 0$ we obtain:
\begin{equation}
    \kappa_1 = \frac{1}{\sqrt{L^2 - x_{1,s}^2}}\frac{\partial}{\partial s}\left(\frac{x_{1,s}}{L}\right),\quad \omega=\frac{\kappa_2 x_{1,s}}{\sqrt{L^2-x_{1,s}^2}}.\label{eq:kappa_1_omega}
\end{equation}
Using this we can reduce the system of equations to only 3 unknowns: 
\begin{eqnarray}
&&x_{1,t} = -\kappa_2 \sqrt{L^2 - x_{1,s}^2},\label{eq:frenet_x1_finite}\\    
&&L_t = L \varepsilon \frac{x_1}{x_1^2 + r_c^2} \kappa_2 \sqrt{L^2 - x_{1,s}^2},\label{eq:frenet_L_finite}\\    
&&\kappa_{2,t} = -\frac{1}{L}\frac{\partial}{\partial s}\left(\frac{v_s}{L}\right) - 
\frac{\omega_s u + 2 \omega u_s}{L} + \omega^2 v + \kappa_2 (u \kappa_1 - v \kappa_2) + 
\lambda_3 \kappa_1.\label{eq:frenet_kappa2_finite}
\end{eqnarray}
We will show later that the equation~\eqref{eq:frenet_L_finite} can be resolved analytically. Besides, similarly to
$\kappa_1$ we can find $\kappa_2$
\begin{equation}
    \kappa_2 = \frac{L x_{1,t}}{\sqrt{L^2 - x_{1,s}^2}},\nonumber
\end{equation}
and see that it is proportional to the velocity $x_{1,t}$. Thus, the behavior of the system of antiparallel
vortices is governed by 2 quantities: the distance between vortices $x_1$ and the velocity of their approximation
represented by $\kappa_2$.

Even though the choice of the interaction frame allows us to reduce the number of unknowns it is not very useful in
numerical simulation since the 4th derivative of $x_1$ respect to $s$ is required: equation~\eqref{eq:frenet_kappa2_finite}
includes the second derivative of $v$, that according to~\eqref{eq:v_coeff} depends on $\kappa_1$, that is proportional to
the second derivative of $x_1$ due to~\eqref{eq:kappa_1_omega}. In practice, it is better to
use the initial formulation~\eqref{eq:main_equation}.


\section{Algorithm for the numerical solution}\label{sec:algorithm}

Assume that at time $t = t_0$ the initial values $\mathbf{X}$, $\mathbf{T}$ are given. Suppose also that we know function
$L_0(s)$ from the equation~\eqref{eq:arclength} and its derivative $L_0'(s)$. We use the following algorithm to obtain the numerical solution up to the time $t_{\text{end}}$
with space discretization step $h$ and time step $\Delta t$: 

\begin{algorithmic}[1]

\While{$t < t_{\text{end}}$}
    \State{accuracy\_test\_passed $\gets$ False}        
    \While{accuracy\_test\_passed is False}
    \State{$\Delta \mathbf{X} \gets \mathbf{0}$; $\Delta \mathbf{T} \gets \mathbf{0}$}
    \State{$\mathbf{X}_{\text{error}} \gets \mathbf{0}$; $\mathbf{T}_{\text{error}} \gets \mathbf{0}$}
    \For{$k \gets 1$ to $6$}
      \State{$\boldsymbol{\xi} \gets \mathbf{X}$; $\boldsymbol{\tau} \gets \mathbf{T}$}
      \For{$q \gets 1$ to $k-1$}
        \State{$\boldsymbol{\xi} \gets \boldsymbol{\xi} + \Delta t\ \alpha_{k-1,q}\ \mathbf{X}^{(q)}$;
               $\boldsymbol{\tau} \gets \boldsymbol{\tau} + \Delta t\ \alpha_{k-1,q}\ \mathbf{T}^{(q)}$}
      \EndFor
      \State{$\boldsymbol{\tau}_{s,h} \gets \text{FirstDerivative}(\boldsymbol{\tau}, h)$}
      \State{$\boldsymbol{\tau}_{ss,h} \gets \text{SecondDerivative}(\boldsymbol{\tau}, h)$}
      \State{$a \gets L_0 (\xi_1^2 + r_c^2)^{-\varepsilon/2}$ // modulus of the tangential vector~\eqref{eq:arclength}}
      \State{$b \gets \frac{L'_0}{L_0} - \varepsilon \frac{\xi_1}{\xi_1^2 + r_c^2}\tau_1$ // correction $|\boldsymbol{\tau}|_s / |\boldsymbol{\tau}|$}
      \State{$\mathbf{X}^{(k)} \gets \frac{\boldsymbol{\tau} \wedge \boldsymbol{\tau}_{s,h}}{a^3} - 
                                  \frac{\varepsilon \xi_1}{a (\xi_1^2 + r_c^2)} \boldsymbol{\tau}\wedge \mathbf{e}_1$}
      \State{$\mathbf{T}^{(k)} \gets \frac{\boldsymbol{\tau} \wedge \left(\boldsymbol{\tau}_{ss,h}\!-\!3b \boldsymbol{\tau}_{s,h}\right)}{{a^3}} - 
                                 \frac{\varepsilon  \tau_1 (r_c^2 - \xi_1^2)}{a (\xi_1^2 + r_c^2)^2} \boldsymbol{\tau}\wedge \mathbf{e}_1 - 
                                 \frac{\varepsilon \xi_1}{a (\xi_1^2 + r_c^2)} \left(\boldsymbol{\tau}_{s,h}\!-\!b \boldsymbol{\tau}\right)\wedge \mathbf{e}_1$}      
      \State{$\Delta \mathbf{X} \gets \Delta \mathbf{X} + \Delta t\ c_k\ \mathbf{X}^{(k)}$;
             $\Delta \mathbf{T} \gets \Delta \mathbf{T} + \Delta t\ c_k\ \mathbf{T}^{(k)}$}

      \State{$\mathbf{X}_{\text{error}} \gets \mathbf{X}_{\text{error}} + \Delta t\ \hat{c}_k\ \mathbf{X}^{(k)}$;
             $\mathbf{T}_{\text{error}} \gets \mathbf{T}_{\text{error}} + \Delta t\ \hat{c}_k\ \mathbf{T}^{(k)}$}
    \EndFor{// $k \gets 1$ to $6$}
    \State{$\text{error} \gets h \sqrt{\|\mathbf{X}_{\text{error}}\|^2 + \|\mathbf{T}_{\text{error}}\|^2}$}
    \If{$\text{error}  < \text{threshold}$}
        \State{accuracy\_test\_passed $\gets$ True}        
        \State{$\mathbf{X}\gets \mathbf{X} + \Delta \mathbf{X}$}
        \State{$\mathbf{T}\gets \mathbf{T} + \Delta \mathbf{T}$}
        \State{$t\gets t + \Delta t$}
    \EndIf
    \State{$\Delta t_{\text{new}} \gets 0.9 \Delta t {\left(\frac{\text{treshold}}{\text{error}}\right)}^{0.2}$}
    \State $\Delta t \gets 2^{\lfloor \log_2\left(\Delta t_{\text{new}}/\Delta t\right) \rfloor} \Delta t$            
    \EndWhile{// accuracy test not passed}    
\EndWhile{// $t < t_{\text{end}}$}

\end{algorithmic}
The coefficients $a_{kq}$, $c_k$, and $\hat{c}_k$ are given in Butcher table~\ref{tab:butcher_table}.

\begin{table}
    \centering
    \begin{tabular}{c|c c c c c}
        $2/9$ & $2/9$ & & & & \\
        $1/3$ & $1/12$ & $1/4$ & & & \\
        $3/4$ & $69/128$ & $-243/128$ & $135/64$ & & \\
        $1$ & $-17/12$ & $27/4$ & $-27/5$ &  $16/15$ & \\
        $5/6$ & $65/432$ & $-5/16$ & $13/16$ &  $4/27$ & $5/144$ \\\hline
        & $1/9$ & $0$ & $9 / 20$ & $16/45$ & $1/12$  \\\hline
        & $47/450$ & $0$ & $12 / 25$ & $32/225$ & $1/30$  
    \end{tabular}
    \caption{Butcher table for Runge-Kutta-Fehlberg method~\cite{fehlberg1969}}\label{tab:butcher_table}
\end{table}


\section{The eye-shaped vortex}~\label{sec:eye-shaped-details}

The initial configuration of the eye-shaped vortex is given by
\begin{equation}
    \mathbf{X}(s,0) = \begin{pmatrix}
        b \sin{s} \\ s - \pi / 2 \\ -b \sqrt{\frac{1 + \cos{\theta}}{1 - \cos{\theta}} - \frac{1}{b^2}}\cos{s}       
    \end{pmatrix},\ s\in(0,\pi],\label{eq:eye_vortex}
\end{equation}
where $b$ is the thickness of the eye, $\theta$ is the angle of the corner, and the part $s\in(\pi,2\pi]$ is
obtained by reflection respect to the plane $x = 0$. Note, that the component $x_3$ is real only if the 
expression below square root is positive therefore for large angles $\theta$ we also have to use a large separation $b$.

\begin{figure}
    \centering
    \includegraphics[width=0.49\textwidth, trim={16cm 4cm 12cm 2cm}, clip]{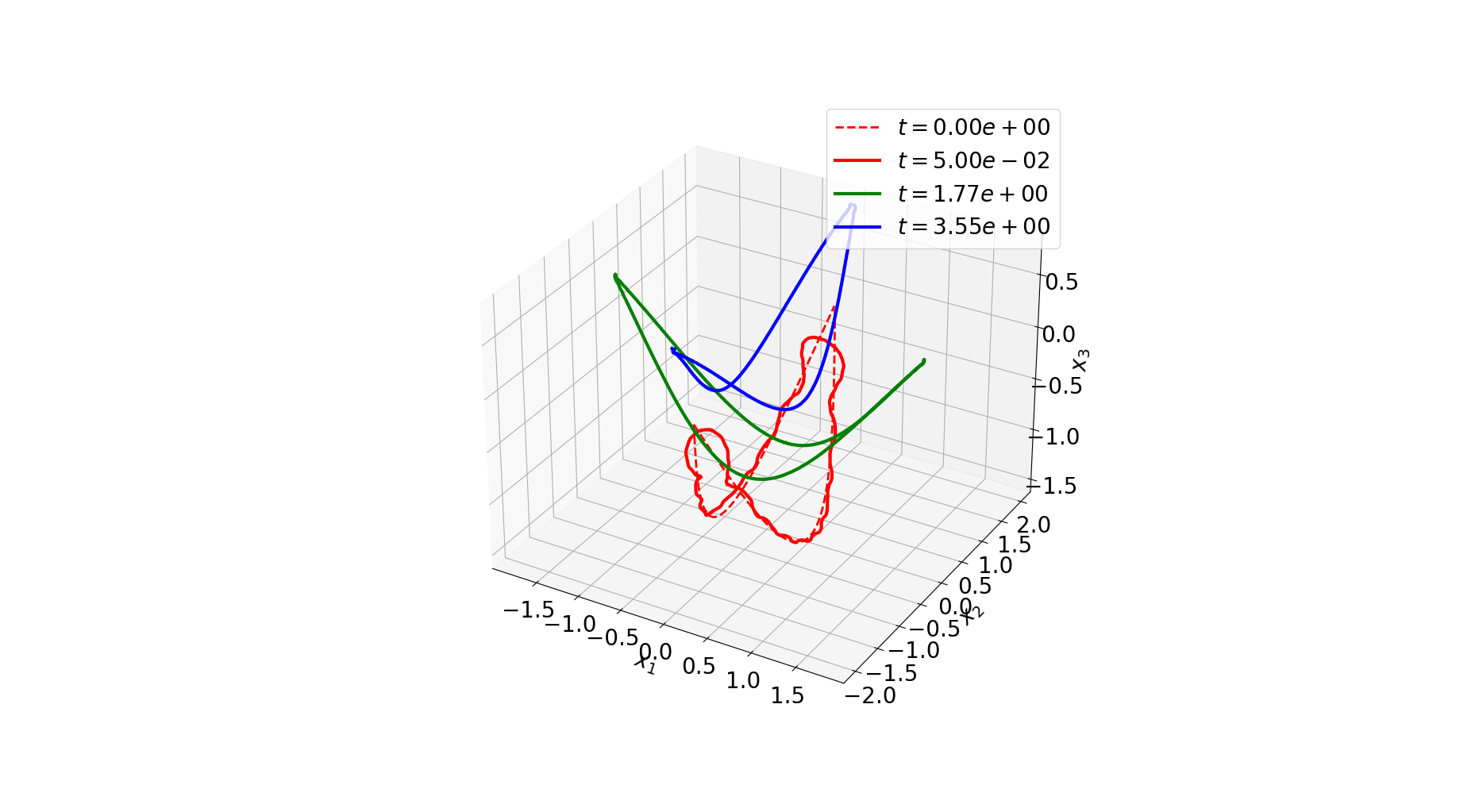}
    \includegraphics[width=0.49\textwidth, trim={16cm 4cm 12cm 2cm}, clip]{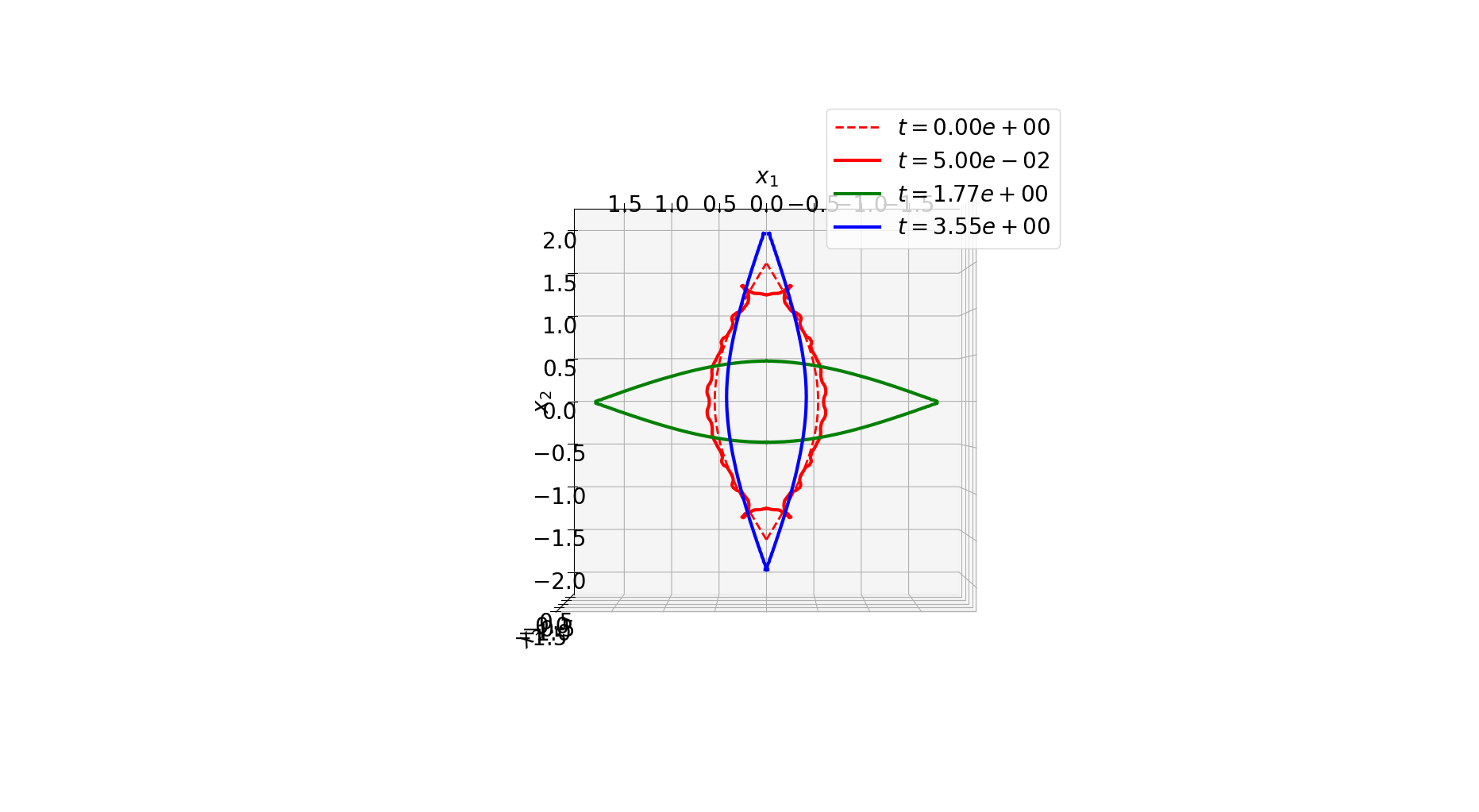}
    \caption{The eye-shaped vortex~\eqref{eq:eye_vortex} with $\theta = \pi / 6$ and $b = 0.4$ at different time moments}
    \label{fig:eye_shaped_vortex}
\end{figure}

The evolution of the eye-shape vortex~\eqref{eq:eye_vortex} with $\theta = \pi / 6$ and $b = 0.4$ is shown in the figure~\ref{fig:eye_shaped_vortex}. 
It is possible to see that the movement is quasi-periodic with period $T = 3.55$ since at that time we see that the vortex has again the eye-shape with
the same orientation (blue line) but a slightly different parameters than at the initial time (red dashed line). We can also see
that at a half-period time the vortex also has an eye-shape but is rotated (green line) similarly to the polygonal vortex~\cite{delahoz2018}. At a rational
fraction $p / q$ of the quasi-period $T$ we can also see a non planar curve with $q$ or $2 q$ corners for even and odd values of $q$
respectively. This behavior also coincides with the one of a polygon. 

\begin{figure}
    \centering
    \begin{subfigure}[c]{0.7\textwidth}
        \centering
        \includegraphics[width=\textwidth, trim={1cm 0cm 3cm 1cm}, clip]{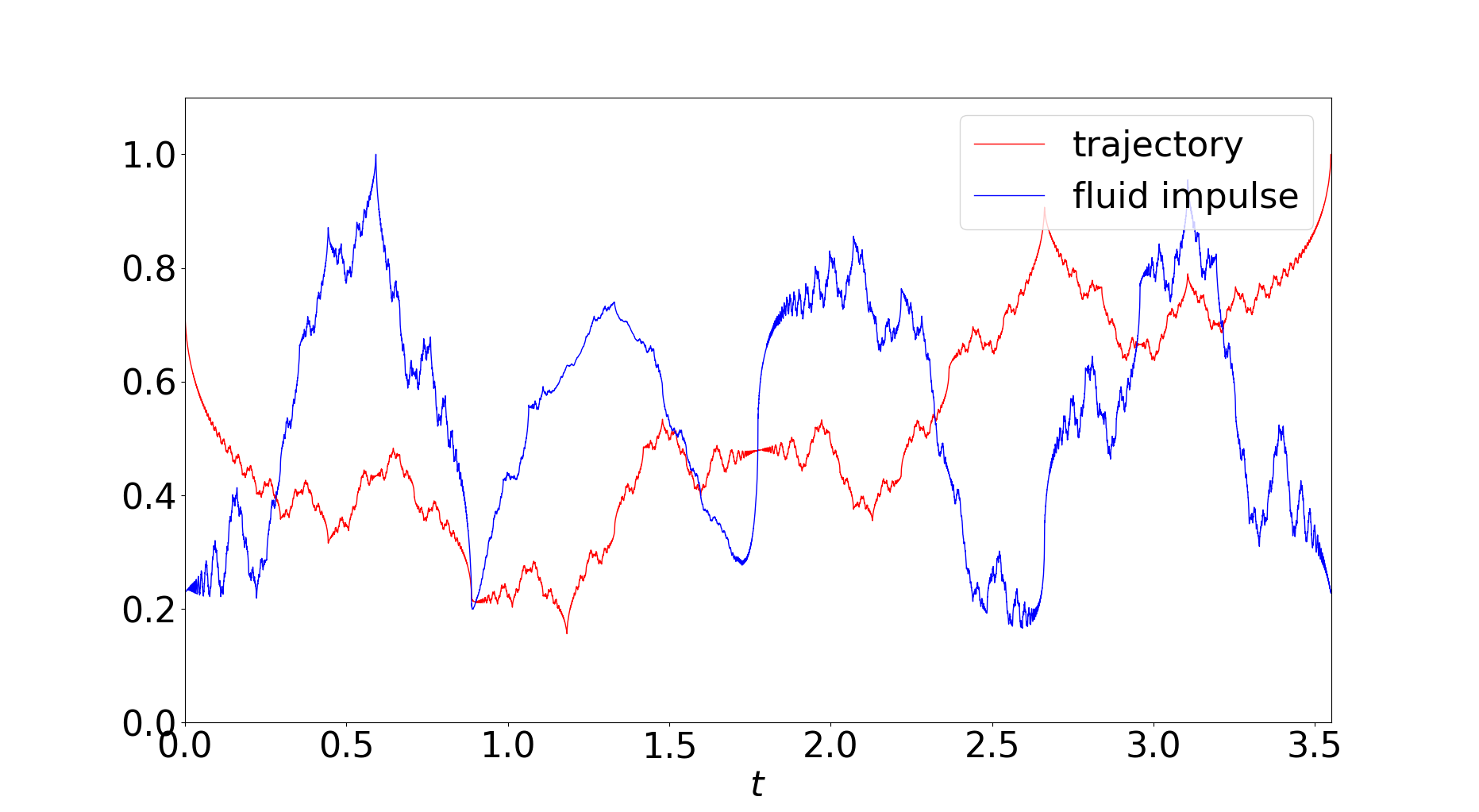}
        \caption{Trajectory of the corner $\mathbf{X}(0,t)$ and fluid impulse around the 
                 corner~\eqref{eq:fluid_impulse}}\label{fig:traj_fi}    
    \end{subfigure}
    \begin{subfigure}[c]{0.7\textwidth}
        \centering
        \includegraphics[width=\textwidth, trim={7mm 0cm 3cm 1cm}, clip]{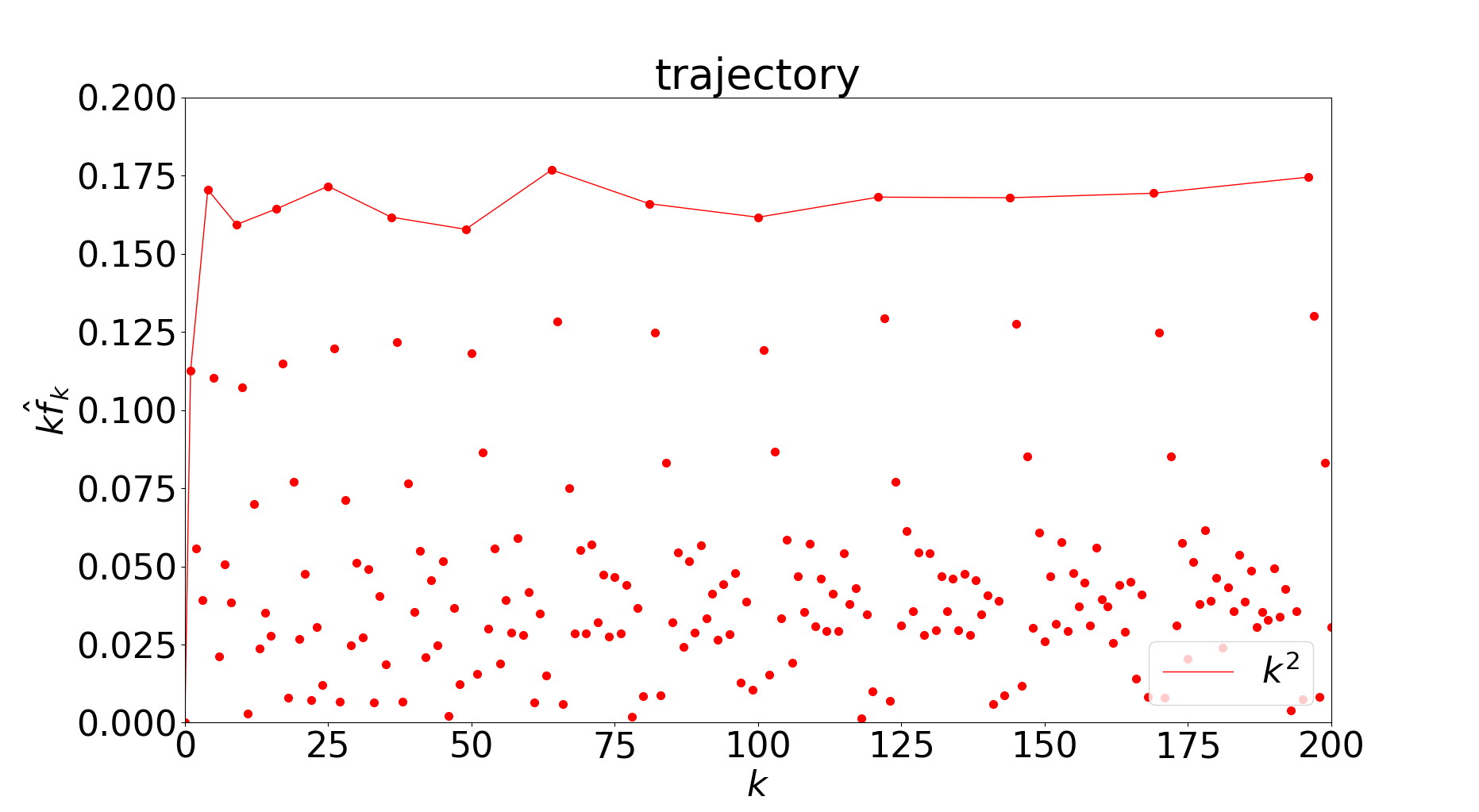}
        \caption{Fourier coefficients of the trajectory}\label{fig:traj_fourier}    
    \end{subfigure} 
    \begin{subfigure}[c]{0.7\textwidth}
        \centering   
        \includegraphics[width=\textwidth, trim={1cm 0cm 3cm 1cm}, clip]{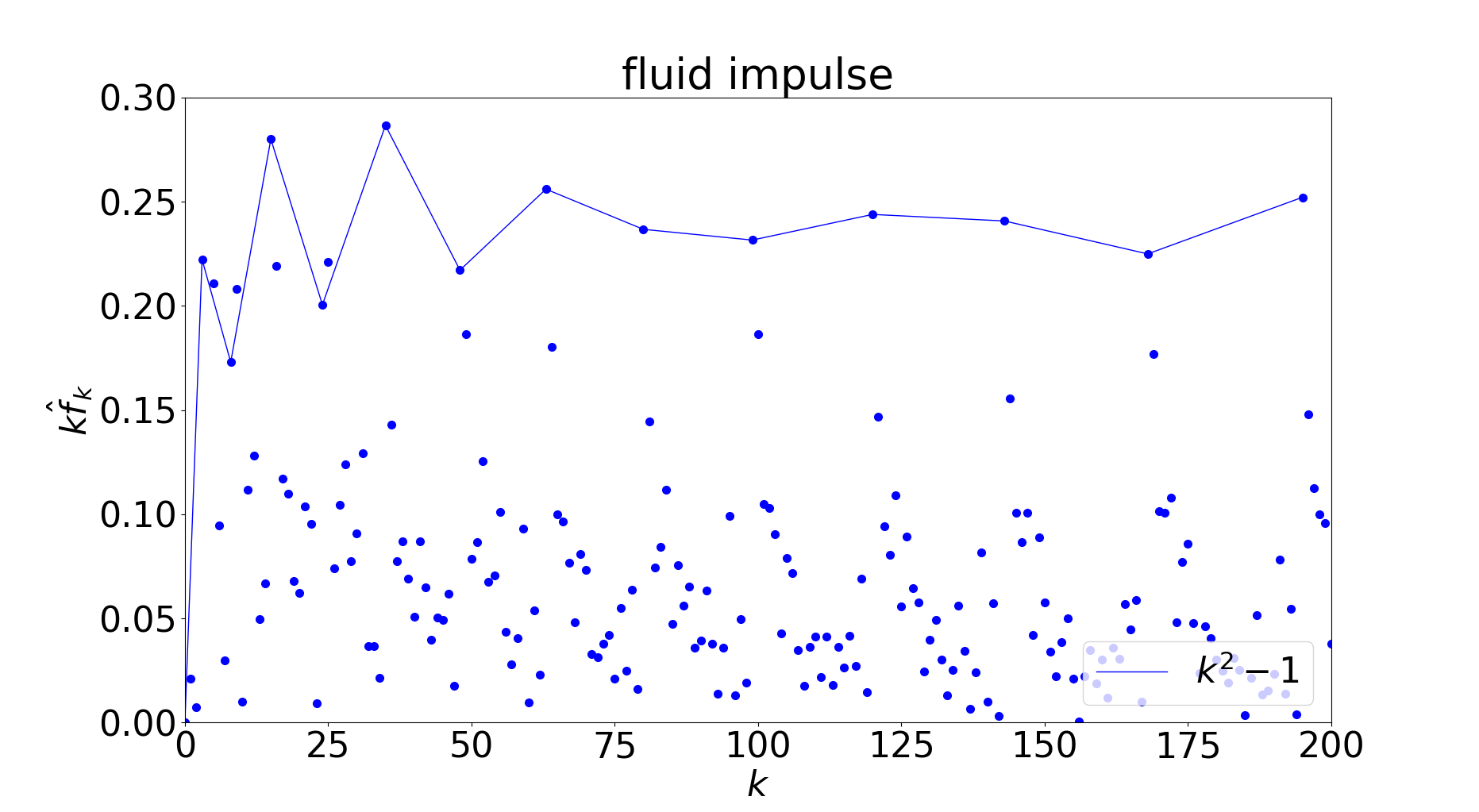}
        \caption{Fourier coefficients of the fluid impulse}\label{fig:fi_fourier}    
    \end{subfigure}     
    \caption{Trajectory, fluid impulse, and their Fourier coefficients for the eye-shaped vortex 
             with  $\theta = \pi / 6$ and $b = 0.4$, the lines connect frequencies corresponding to squares of integers.}\label{fig:traj_fi_fourier}
\end{figure}

In figure~\ref{fig:traj_fi} the trajectory of the corner $\mathbf{X}(0,t)$ and the fluid impulse around the corner are depicted.
Analysing the Fourier coefficients (figures~\ref{fig:traj_fourier} and~\ref{fig:fi_fourier}) we can see that for both the trajectory and the fluid impulse the frequencies
corresponding to squares of integers are dominating similarly to what happens in the case of regular polygons~\cite{delahoz2018}. 
It makes the behavior similar to the RNDF~\eqref{eq:rndf}.

    \printbibliography

@article{crow1970,
% Pioneric paper of the vortex reconnection of Crow instability
    author = "Crow, S. C.",
    title = "{Stability theory for a pair of trailing vortices}",
    journal = "AIAA Journal",
    volume = "8",
    number = "12",
    pages = "2172 -- 2179",
    year = "1970",
    DOI = "https://doi.org/10.2514/3.6083",
}

@article{klein1995,
    author = {Klein R. and Majda A. J. and Damodaran K.},
    title = "{Simplified equations for the interaction of nearly parallel vortex filaments}",
    journal = " J . Fluid Mech.",
    volume = "288",
    pages = "201 -- 248",
    year = "1995",
    DOI = "https://doi.org/10.1017/S0022112095001121",
}

@article{klein1991a,
    author = {Klein R. and Majda A. J.},
    title = "{Self-stretching of a perturbed vortex filament I. The asymptotic equation for deviations from a straight line}",
    journal = "Phys. D: Nonlinear Phenomena",
    volume = "49",
    number = 3,
    pages = "323 -- 352",
    year = "1991",
    DOI = "https://doi.org/10.1016/0167-2789(91)90151-X",
}

@article{hussain2020,
    author = {Yao J. and Hussain F.},
    title = "{A physical model of turbulence cascade via vortex reconnection sequence and avalanche}",
    journal = " J . Fluid Mech.",
    volume = "883",
    year = "2020",
    DOI = "https://doi.org/10.1017/jfm.2019.905",
}

@article{krstulovic2017,
    author = {Villois A. and Proment D. and Krstulovic G.},
    title = "{Universal and nonuniversal aspects of vortex reconnections in superfluids}",
    journal = "Phys. Rev. Fluids",
    volume = "2",
    year = "2017",
    DOI = "https://doi.org/10.1103/PhysRevFluids.2.044701",
}

@article{how1998,
    author = {Hou T.Y. and Klapper I. and Si H.},
    title = "{Removing the Stiffness of Curvature in Computing 3-D Filaments}",
    journal = "J. Comp. Physics",
    volume = "143",
    year = "1998",
    pages = "628 -- 664",
    DOI = "https://doi.org/10.1006/jcph.1998.5977",
}

@book{shaffman1992,
	author = "Shaffman P.G.",
	title = "{Vortex Dynamics}", 
	year = "1992",
	isbn={0-521-42058-X},
    publisher={Cambridge University Press},    
}

@article{pontin2018,
    author = {McGavin P. and Pontin D.I.},
    title = "{Vortex line topology during vortex tube reconnection}",
    journal = "Phys. Rev. Fluids",
    volume = "3",
    year = "2018",
    DOI = "https://doi.org/10.1103/PhysRevFluids.3.054701",
}

@article{lim1992,
% Experiments of the reconnection of two vortex rings
  author          = {Lim T.T. and Nickels T.B.},
  journal         = {NATURE},
  title           = {Instability and reconnect in the head-on collision of two vortex rings},
  volume          = {357},
  year            = {1992}  
}

@article{nemirovskii2020,
  author          = {Nemirovskii S.K.},
  journal         = {Journal of Engineering Thermophysics},
  number          = {1},
  title           = {Statistical Signature of Vortex Filaments in Classical Turbulence: Dog or Tail?},
  volume          = {29},
  year            = {2020}
}

@article{laporte2002,
  author          = {Laporte F. and Leweke T.},
  journal         = {AIAA JOURNAL},
  number          = {12},
  title           = {Elliptic Instability of Counter-Rotating Vortices: Experiment and Direct Numerical Simulation},
  volume          = {40},
  year            = {2002}
}

@article{ledizes2005,
  author          = {Le Dizès S. and Lacaze L.},
  journal         = {Journal of Fluid Mechanics},
  number          = {},
  title           = {An asymptotic description of vortex Kelvin modes},
  volume          = {542},
  year            = {2005},
  pages           = {69 -- 96},
  DOI             = "doi:10.1017/S0022112005005185"
}

@article{hussain2011,
  author          = {Hussain F. and Duraisamy K.},
  journal         = {Physics of Fluids},
  number          = {23},
  title           = {Mechanics of viscous vortex reconnection},
  DOI             = {https://doi.org/10.1063/1.3532039},
  year            = {2011}
}

@book{butcher,
  author         = {Butcher J.C.},
  editor         = {II},
  publisher      = {John Wiley \& Sons Ltd},
  title          = {Numerical Methods for Ordinary Differential Equations},
  year           = {2008}
}

@phdthesis{buttke,
  author      = {Buttke T.F.},
  school      = {University of California, Berkeley},
  title       = {A numerical study of superfluid turbulence in the self-induction approximation},
  year        = {1986}
}

@article{delahoz2014,
  author          = {Fracisco de la Hoz and Luis Vega},
  journal         = {Nonlinearity},
  number          = {27},
  title           = {Vortex Filament Equation for a Regular Polygon},
  volume          = {12},
  year            = {2014}
}

@article{schwarz1985,
  author          = {Schwarz K.W.},
  journal         = {PHYSICAL REVIEW 8},
  number          = {9},
  title           = {Three-dimensional vortex dynamics in superfluid {$^4He$}: Line-line and line-boundary interactions},
  volume          = {31},
  year            = {1985}
}

@article{vega2003,
  author          = {Gutiérrez S. and Rivas J. and Vega L.},
  journal         = {Comm. in PDEs.},
  number          = {28},
  title           = {Formation of Singularities and Self-Similar Vortex Motion Under the Localized Induction Approximation},  
  year            = {2003}
}

@article{delahoz2018,
  author          = {Francisco de la Hoz and Luis Vega},
  journal         = {J. Nonlinear Sci},
  number          = {28},
  title           = {On the Relationship Between the One-Corner Problem and the M-Corner Problem for the Vortex Filament Equation},
  doi             = {https://doi.org/10.1007/s00332-018-9477-7},
  year            = {2018}
}

@article{rosenhead1930,
  author          = {Rosenhead L.},
  journal         = {Proc. R. Soc. Lond. A},
  number          = {127},
  title           = {The Spread of Vorticity in the Wake behind a Cylinder},
  pages          = {590--612},
  year            = {1930},
  doi             = {https://doi.org/10.1098/rspa.1930.0078}
}

@article{banica2022,
  author          = {Banica V. and Vega L.},
  journal         = {Arch Rational Mech Anal},
  number          = {244},
  title           = {Riemann's Non-differentiable Function and the Binormal Curvature Flow},  
  year            = {2022},
  pages           = {501--540},
  doi             = {https://doi.org/10.1007/s00205-022-01769-1}
}

@article{fehlberg1969,
  author          = {Fehlberg E.},
  journal         = {National aeronautics and space administration},
  title           = {Low-order classical Runge-Kutta formulas with stepsize control and their application to some heat transfer problems},
  volume          = {315},
  year            = {1969}
}

@article{jaffard1996,
  author          = {Jaffard S.},
  journal         = {Revista Matematica Iberoamericana},
  number          = {2},
  title           = {The spectrum of singularities of Riemann's function},
  volume          = {12},
  year            = {1996}
}

@article{wendt2007,
  author          = {Wendt H. and Abry P. and Jaffard S.},
  journal         = {IEEE Signal Processing Magazine},  
  title           = {Bootstrap for Empirical Multifractal Analysis},
  volume          = {38},
  year            = {2007},
  doi             = {https://doi.org/10.1109/MSP.2007.4286563}
}

@article{turiel2006,
  author          = {Turiel A. and Pérez-Vicente C.J. and Grazzini J.},
  journal         = {Journal of Computational Physics},  
  title           = {Numerical methods for the estimation of multifractal singularity spectra on sampled data: A comparative study},
  volume          = {216},
  year            = {2006},
  doi             = {doi:10.1016/j.jcp.2005.12.004}
}

@article{kumar2022,
  author          = {Kumar S. and Ponce-Vanegas F. and Roncal L. and Vega L.},
  journal         = {arXiv preprint arXiv:2202.06645},
  title           = {The Frisch-Parisi formalism for fluctuations of the Schr\" odinger equation},
  year            = {2022}
}

@article{muzy1993,
  author          = {Muzy J.F. and Barcy E. and Arneodo A.},
  journal         = {Int. J. of Bifurcation and Chaos},
  number          = {2},
  title           = {The Multifractal Formalism revisited with Wavelets},
  volume          = {4},
  year            = {1993},
  pages           = {245-302}
}

@article{kida1988,
  author          = {Kida S. and Takaoka M. and Hussain F.},
  journal         = {Phys. Fluids A},
  number          = {4},
  title           = {Reconnection of two vortex rings},
  volume          = {1},
  year            = {1988}
}

@article{bewley2008,
  author          = {Bewley G.P. and Paoletti M.S. and Sreenivasan K.R. and Lathrop D.P.},
  journal         = {PNAS},
  number          = {37},
  title           = {Characterization of reconnecting vortices in superfluid helium},
  volume          = {105},
  year            = {2008},
  doi             = {www.pnas.org/cgi/doi/10.1073/pnas.0806002105}
}

@article{leweke2016,
  author          = {Leweke T. and Le Dizes S. and Williamson C.H.K.},
  journal         = {Annu. Rev. Fluid Mech.},
  number          = {48},
  title           = {Dynamics and Instabilities of Vortex Pairs},
  pages           = {1 -- 35},
  year            = {2016},
  doi             = {10.1146/annurev-fluid-000000-000000}
}

@article{han2000,
  author          = {Han J. and Lin Y.L. and Schowalter D.G and Arya S. P. and Proctor F.H.},
  journal         = {AIAA Journal},
  number          = {2},
  title           = {Large Eddy Simulation of Aircraft Wake Vortices within homogeneous Turbulence : Crow Instability},
  volume          = {38},
  year            = {2000},
  doi             = {10.2514/2.956}
}

@article{ortega2003,
  author          = {Ortega J.M. and Bristol R.L. and Savas \"O},
  journal         = {J. Fluid Mech.},
  pages           = {35 -- 84},
  title           = {Experimental study of the instability of unequal-strength counter-rotating vortex pairs},
  volume          = {474},
  year            = {2003},
  doi             = {10.1017/S0022112002002446}
}

@article{fonda2019,
  author          = {Fonda E. and Sreenivasan K.R. and Lathrop D.P.},
  journal         = {PNAS},
  number          = {6},
  title           = {Reconnection scaling in quantum fluids},
  volume          = {116},
  year            = {2019},
  pages           = {1924 -- 1928},
  doi             = {https://doi.org/10.1073/pnas.1816403116}
}

@article{lipniacki2003,
  author          = {Lipniacki T.},
  journal         = { J. Fluid Mech.},
  number          = {447},
  title           = {Quasi-static solutions for quantum vortex motion under the localized induction approximation},
  year            = {2003}
}

@article{jerrard2015,
  author          = {Jerrard R.L. and Smets D.},
  journal         = {Jour. Eur. Math. Soc.},
  number          = {6},
  title           = {On the motion of a curve by its binormal curvature},
  volume          = {17},
  year            = {2015},
  pages = {1487 -- 1515}
}

@article{banica2017,
  author          = {Banica V. and Faou E. and Miot E.},
  journal         = {Pure and Appl. Math.},
  title           = {Collision of Almost Parallel Vortex Filaments},
  volume          = {70},
  year            = {2017},
  doi = {https://doi.org/10.1002/cpa.21637}
}

@article{hussain2022,
  author          = {Yao J. and Hussain F.},
  journal         = {Annual Review of Fluid Mechanics},
  pages           = {317 -- 347},
  title           = {Vortex Reconnection and Turbulence Cascade},
  volume          = {54},
  year            = {2022},
  doi             = {https://doi.org/10.1146/annurev-fluid-030121-125143}
}

@article{hussain1995,
  author          = {Jeong J. and Hussian F.},
  journal         = {J. Fluid Mech},
  title           = {On the identification of a vortex},
  volume          = {285},
  year            = {1995},
  pages           = {69 -- 94},
  doi             = {10.1017/S0022112095000462}
}

@article{melander1988,
  author          = {Melander M.V. and Hussain F.},
  journal         = {Proceedings of the 2nd Summer Program of the Center of Turbulence Research, Stanford, CA: Cent. Turbul. Res.},
  title           = {Cut-and-connect of two antiparallel vortex tubes},
  pages           = {257 -- 286},
  year            = {1988}
}

@article{brenner2016,
  author          = {Brenner M.P. and Hormoz S. and Pumir A.},
  journal         = {Phys. Rev. Fluids},
  number          = {8},
  title           = {Potential singularity mechanism for the Euler equations},
  volume          = {1},
  year            = {2016},
  doi             = {https://doi.org/10.1103/PhysRevFluids.1.084503}
}

\end{document}